\documentclass[10pt,leqno]{article}

\usepackage[%
  tmargin=1.1in,bmargin=1.1in,%
  lmargin=1.4in,rmargin=1.4in,%
  marginparsep=0.2in, marginparwidth=1.0in%
]{geometry}
\AtBeginDocument{%
  \fontsize{9.6}{12}\selectfont
  \setlength{\abovedisplayskip}{1.5ex plus 0.3ex minus 0.3ex}%
  \setlength{\abovedisplayshortskip}{1.0ex plus 0.3ex minus 0.3ex}%
  \setlength{\belowdisplayskip}{1.5ex plus 0.3ex minus 0.3ex}%
  \setlength{\belowdisplayshortskip}{1.0ex plus 0.3ex minus 0.3ex}%
}

\usepackage{datetime2}
\usepackage{fancyhdr}
\usepackage[hyphens]{url}
\usepackage[style=numeric,backend=bibtex8,maxbibnames=100]{biblatex}
\usepackage[newparttoc]{titlesec}
\usepackage[dotinlabels]{titletoc}
\usepackage{appendix}
\usepackage[letterspace=50]{microtype}
\usepackage[bottom]{footmisc}
\usepackage[inline]{enumitem}
\usepackage{xspace}
\usepackage{calc}
\usepackage{etoolbox}
\usepackage{ifthen}
\usepackage{tikz}
\usepackage{tikz-cd}
\usepackage{textcase}

\definecolor{cite}{rgb}{0.45,0.5,0.45}
\definecolor{link}{rgb}{0.45,0.5,0.45}
\definecolor{url}{rgb}{0.45,0.45,0.5}
\usepackage[%
  hyperfootnotes=false,
  colorlinks,%
  linkcolor=link,%
  citecolor=cite,%
  urlcolor=url,%
]{hyperref}
\DeclareCiteCommand{\citebare}
  {\usebibmacro{prenote}}
  {\usebibmacro{citeindex}%
   \usebibmacro{cite}}
  {\multicitedelim}
  {\usebibmacro{postnote}}

\renewbibmacro{in:}{}
\DeclareFieldFormat*{title}{\mkbibquote{#1}}
\DeclareFieldFormat*{booktitle}{\mkbibitalic{#1}}
\DeclareFieldFormat*{journaltitle}{#1\isdot}
\DeclareFieldFormat*{volume}{#1}
\DeclareFieldFormat*[online,misc,incollection,inproceedings]{date}{\mkbibparens{#1}}
\DeclareFieldInputHandler{pages}{}

\newtoggle{bbx:series}
\newtoggle{bbx:volume}

\DeclareBiblatexOption{global,type,entry}[boolean]{series}[true]{%
  \settoggle{bbx:series}{#1}}
\DeclareBiblatexOption{global,type,entry}[boolean]{volume}[true]{%
  \settoggle{bbx:volume}{#1}}

\newbibmacro*{addendum+pubstate+date}{%
  \printfield{addendum}%
  \newunit\newblock
  \printfield{pubstate}
  \usebibmacro{date}}

\newbibmacro*{doi+eprint+url+date}{%
  \iftoggle{bbx:doi}
  {\printfield{doi}}
  {}%
  \newunit\newblock
  \iftoggle{bbx:eprint}
  {\usebibmacro{eprint}}
  {}%
  \iftoggle{bbx:url}
  {\usebibmacro{url+urldate}}
  {}
  \usebibmacro{date}}

\DeclareBibliographyDriver{incollection}{%
  \usebibmacro{bibindex}%
  \usebibmacro{begentry}%
  \usebibmacro{author}%
  \newunit
  \usebibmacro{title}%
  \newunit
  \printfield{booktitle}%
  \newunit
  \iftoggle{bbx:series}{\printfield{series}}%
  {}
  \iftoggle{bbx:volume}{\printfield{volume}}%
  {}
  \usebibmacro{date}%
  \usebibmacro{finentry}}

\DeclareBibliographyDriver{inproceedings}{%
  \usebibmacro{bibindex}%
  \usebibmacro{begentry}%
  \usebibmacro{author}%
  \newunit
  \usebibmacro{title}%
  \newunit
  \printfield{booktitle}%
  \newunit
  \iftoggle{bbx:series}{\printfield{series}}%
  {}
  \iftoggle{bbx:volume}{\printfield{volume}}%
  {}
  \usebibmacro{date}%
  \usebibmacro{finentry}}

\DeclareBibliographyDriver{online}{%
  \usebibmacro{bibindex}%
  \usebibmacro{begentry}%
  \usebibmacro{author}%
  \newunit
  \usebibmacro{title}%
  \newunit
  \usebibmacro{doi+eprint+url+date}%
  \usebibmacro{finentry}}

\DeclareBibliographyDriver{misc}{%
  \usebibmacro{bibindex}%
  \usebibmacro{begentry}%
  \usebibmacro{author}%
  \newunit
  \usebibmacro{title}%
  \newunit
  \printfield{howpublished}%
  \newunit
  \printfield{type}%
  \newunit
  \printfield{version}%
  \newunit
  \printfield{note}%
  \newunit
  \usebibmacro{doi+eprint+url}%
  \newunit
  \usebibmacro{addendum+pubstate+date}%
  \usebibmacro{finentry}}

\DeclareFieldFormat{eprint:ktheory}{%
  K-theory preprint archives: \href{https://faculty.math.illinois.edu/K-theory/#1}{ #1}%
}

\usepackage{amsmath,amsthm}

\usepackage[T1]{fontenc}
\usepackage{lmodern}
\usepackage{MnSymbol}
\usepackage[bbgreekl]{mathbbol}
\DeclareSymbolFontAlphabet{\mathbbl}{bbold}

\makeatletter
\renewcommand\Mn@Text@With@MathVersion[1]{%
  \textnormal{\ifx\Mn@Bold\math@version\bfseries\fi#1}%
}
\makeatother

\usepackage[%
  cal=euler,  calscaled=1.02,%
  bb=boondox, bbscaled=1.05,%
  frak=euler, frakscaled=1.0,%
]{mathalfa}  

\makeatletter
\g@addto@macro\bfseries{\boldmath}
\makeatother

\usepackage[textsize=footnotesize,backgroundcolor=orange!50]{todonotes}
\usepackage[compress,capitalize]{cleveref}

\let\theoldbibliography\thebibliography
\renewcommand{\thebibliography}[1]{%
  \theoldbibliography{#1}%
  \setlength{\parskip}{0ex}
  \setlength{\itemsep}{0.5ex plus 0.2ex minus 0.2ex}
  \small
}
\apptocmd{\thebibliography}{\raggedright}{}{}

\renewcommand{\title}[1]{\newcommand{\thetitle}{#1}}
\renewcommand{\author}[1]{\newcommand{\theauthor}{#1}}
\renewcommand{\date}[1]{\newcommand{\thedate}{#1}}

\renewcommand{\maketitle}{%
  \begin{center}
    {\linespread{1.15}%
      \bfseries\MakeTextUppercase%
      \thetitle\par} \bigskip
    {\footnotesize\MakeUppercase \theauthor\par} \bigskip
  \end{center}
  \thispagestyle{fancy}
}

\renewenvironment{abstract}{\noindent\begin{center}\begin{minipage}{0.8\linewidth}\footnotesize{\scshape Abstract.}}{\end{minipage}\end{center}\medskip}

\newlength{\tagsep}
\titleformat{\part}{\centering\titlerule\vspace{1.0ex}\bfseries}{Part \thepart.}{1.5\tagsep}{}[\vspace{1.0ex}\titlerule]
\titleformat{\section}{\centering}{\textsection\thesection.}{1.5\tagsep}{\scshape}
\titleformat{\subsection}[runin]{}{\fontseries{b}\selectfont\textsection\bfseries\thesubsection.}{1.5\tagsep}{\bfseries}[.]
\titleformat{\subsubsection}[runin]{}{\textsection\thesubsubsection.}{\tagsep}{\itshape}[.]

\titlespacing*{\part}{0pt}{0ex}{4ex}
\titlespacing*{\section}{0pt}{4ex}{\medskipamount}
\titlespacing*{\subsection}{0pt}{\bigskipamount}{0.5em}
\titlespacing*{\subsubsection}{0pt}{\medskipamount}{0.5em}

\titlecontents{part}[0em]{\vspace{1.25ex}}{Part \thecontentslabel.\hspace{0.5em}}{}{}
\dottedcontents{section}[1.6em]{\vspace{0.25ex}}{1.6em}{0.7em}
\dottedcontents{subsection}[2.6em]{}{2.6em}{0.7em}

\newcommand{\crefeqfmt}[1]{
  \crefformat{#1}{(##2##1##3)}
  \Crefformat{#1}{(##2##1##3)}
  \crefrangeformat{#1}{(##3##1##4--##5##2##6)}
  \Crefrangeformat{#1}{(##3##1##4--##5##2##6)}
  \crefmultiformat{#1}{(##2##1##3}{, ##2##1##3)}{, ##2##1##3}{, ##2##1##3)}
  \Crefmultiformat{#1}{(##2##1##3}{, ##2##1##3)}{, ##2##1##3}{, ##2##1##3)}
  \crefrangemultiformat{#1}{(##3##1##4--##5##2##6}{, ##3##1##4--##5##2##6)}{, ##3##1##4--##5##2##6}{, ##3##1##4--##5##2##6)}
  \Crefrangemultiformat{#1}{(##3##1##4--##5##2##6}{, ##3##1##4--##5##2##6)}{, ##3##1##4--##5##2##6}{, ##3##1##4--##5##2##6)}
}
\newcommand{\crefsecfmt}[1]{%
  \crefformat{#1}{\textsection##2##1##3}
  \Crefformat{#1}{\textsection##2##1##3}
  \crefrangeformat{#1}{\textsection\textsection##3##1##4--##5##2##6}
  \Crefrangeformat{#1}{\textsection\textsection##3##1##4--##5##2##6}
  \crefmultiformat{#1}{\textsection\textsection##2##1##3}{--##2##1##3}{, ##2##1##3}{ and~##2##1##3}
  \Crefmultiformat{#1}{\textsection\textsection##2##1##3}{--##2##1##3}{, ##2##1##3}{ and~##2##1##3}
  \crefrangemultiformat{#1}{\textsection\textsection##3##1##4--##5##2##6}{ and~##3##1##4--##5##2##6}{, ##3##1##4--##5##2##6}{ and~##3##1##4--##5##2##6}
  \Crefrangemultiformat{#1}{\textsection\textsectionXS##3##1##4--##5##2##6}{ and~##3##1##4--##5##2##6}{, ##3##1##4--##5##2##6}{ and~##3##1##4--##5##2##6}
}

\crefeqfmt{equation}
\crefeqfmt{enumi}
\crefeqfmt{enumii}
\crefsecfmt{section}
\crefsecfmt{subsection}
\crefsecfmt{appendix}
\crefname{part}{Part}{Parts}
\crefname{chapter}{Chapter}{Chapters}
\crefname{figure}{Figure}{Figures}

\makeatletter

\newcommand{\blocknumfont}{\bfseries}
\newcommand{\blockheadfont}{\bfseries}
\newcommand{\blocknotefont}{\normalfont}
\newcommand{\blockspecialfont}{\itshape}
\newcommand{\blockhorizspace}{0.4em}
\newcommand{\blocknotespace}{0.4em}
\newcommand{\blocknumsep}{}
\newcommand{\blocksep}{.}
\newcommand{\blockvertspace}{\medskipamount}

\newtheoremstyle{block}%
  {\blockvertspace}
  {\blockvertspace}
  {}
  {}
  {} 
  {}
  {0em}
  {\thmname{{\blockheadfont#1}}%
    \@ifnotempty{#1}{\@ifnotempty{#2}{ }}%
    \thmnumber{{\blocknumfont #2\blocknumsep}}%
    \@ifnotempty{#1#2}{{\blockheadfont\blocksep}}%
    \thmnote{\hspace{\blocknotespace}[{\blocknotefont#3}]}%
    \hspace{\blockhorizspace}%
  }

\newtheoremstyle{blockspecial}%
  {\blockvertspace}
  {\blockvertspace}
  {\blockspecialfont}
  {}
  {} 
  {}
  {0em}
  {\thmname{{\blockheadfont#1}}%
    \@ifnotempty{#1}{\@ifnotempty{#2}{ }}%
    \thmnumber{{\blocknumfont #2\blocknumsep}}%
    \@ifnotempty{#1#2}{{\blockheadfont\blocksep}}%
    \thmnote{\hspace{\blocknotespace}[{\blocknotefont#3}]}%
    \hspace{\blockhorizspace}%
  }

\newtheoremstyle{blocknamed}%
  {\blockvertspace}
  {\blockvertspace}
  {}
  {}
  {} 
  {}
  {0em}
  {\@ifempty{#3}{\thmname{\blockheadfont#1}}{\thmname{\blockheadfont#3}}%
    \@ifnotempty{#1#3}{\blockheadfont\blocksep}\hspace{\blockhorizspace}%
  }
  
\theoremstyle{blocknamed}
\newtheorem*{pf}{Proof}

\renewenvironment{proof}[1][]{%
  \pushQED{\qed}%
  \begin{pf}[#1]
}{%
  \popQED\endtrivlist\@endpefalse%
  \end{pf}
}
  
\newcommand{\defblock}[2]{%
  \theoremstyle{block}
  \newtheorem{#1}[block]{#2}%
  \Crefname{#1}{#2}{{#2}s}
  \newtheorem*{#1*}{#2}%
}

\newcommand{\defblockspecial}[2]{%
  \theoremstyle{blockspecial}
  \newtheorem{#1}[block]{#2}%
  \Crefname{#1}{#2}{{#2}s}
  \newtheorem*{#1*}{#2}%
}

\makeatother

\setlist{%
  parsep=0ex, listparindent=\parindent,%
  itemsep=0.5ex, topsep=0.5ex,%
  leftmargin=2.3em,%
}

\setlist[enumerate, 1]{%
  label={\upshape(\arabic*)},%
  ref={\arabic*},%
  widest=9,
}
\setlist[enumerate, 2]{%
  leftmargin=*,%
  label={\upshape(\roman*)},%
  ref=\roman*,%
  widest=iii,
}
\setlist[itemize, 1]{%
  label={\textbf{--}},%
}
\setlist[itemize, 2]{%
  label=--,%
}

\usetikzlibrary{calc,decorations.pathmorphing,shapes,arrows}
\tikzcdset{
  arrow style=tikz,
  diagrams={>={stealth}},
}
  
\newcommand{\from}{\leftarrow}

\newcommand{\lblto}[1]{\xrightarrow{#1}}
\newcommand{\lblfrom}[1]{\xleftarrow{#1}}
\newcommand{\isoto}{\lblto{\raisebox{-0.3ex}[0pt]{$\scriptstyle \sim$}}}
\newcommand{\inj}{\hookrightarrow}

\newcommand{\iso}{\simeq}
\newcommand{\lbliso}[1]{\overset{#1}{\simeq}}

\usepackage{xfrac}

\defblock{claim}{Claim}
\defblock{construction}{Construction}
\defblockspecial{corollary}{Corollary}
\defblock{definition}{Definition}
\defblock{example}{Example}
\defblock{fact}{Fact}
\defblock{hypothesis}{Hypothesis}
\defblockspecial{lemma}{Lemma}
\defblock{notation}{Notation}
\defblock{nothing}{}
\defblockspecial{proposition}{Proposition}
\defblock{question}{Question}
\defblock{remark}{Remark}
\defblockspecial{specification}{Specification}
\defblockspecial{theorem}{Theorem}
\defblock{warning}{Warning}
\defblock{variant}{Variant}

\numberwithin{equation}{subsection}
\numberwithin{block}{subsection}
\makeatletter
\let\c@equation\c@block
\makeatother

\newcommand{\itemref}[2]{\cref{#1}.\ref{#1--#2}}
\newcommand{\itemitemref}[3]{\cref{#1}.\ref{#1--#2}.\ref{#1--#2--#3}}

\bibliography{ref}

\definecolor{arpon}{rgb}{0.7,0.45,0.7}

\newcommand{\cat}[1]{\mathcal{#1}}
\newcommand{\num}[1]{\mathbb{#1}}
\newcommand{\twocat}[1]{\mathbb{#1}}

\newcommand{\A}{\cat{A}}
\renewcommand{\AA}{\num{A}}
\newcommand{\Aff}{\mathrm{Aff}}
\newcommand{\Alg}{\mathrm{Alg}}
\newcommand{\B}{\cat{B}}
\newcommand{\C}{\cat{C}}
\newcommand{\D}{\mathrm{D}}
\newcommand{\Ch}{\mathrm{C}}
\newcommand{\CAlg}{\mathrm{CAlg}}
\newcommand{\Cat}{\mathrm{Cat}}
\newcommand{\Ct}{\widehat\Ch}
\newcommand{\Fun}{\mathrm{Fun}}
\newcommand{\G}{\num{G}}
\newcommand{\Ind}{\mathrm{Ind}}
\newcommand{\K}{\mathrm{K}}
\newcommand{\LMod}{\mathrm{LMod}}
\newcommand{\Loc}{\mathrm{Loc}}
\newcommand{\Map}{\mathrm{Map}}
\newcommand{\Mod}{\mathrm{Mod}}
\newcommand{\Nm}{\mathrm{Nm}}
\renewcommand{\P}{\cat{P}}
\newcommand{\PStk}{\mathrm{PStk}}
\newcommand{\PrL}{\mathrm{Pr^L}}
\newcommand{\PrLst}{\mathrm{Pr^L_{st}}}
\newcommand{\QCoh}{\mathrm{QCoh}}
\renewcommand{\S}{\cat{S}}
\renewcommand{\SS}{\num{S}}
\newcommand{\Span}{\mathrm{Span}}
\newcommand{\Spc}{\mathrm{Spc}}
\newcommand{\Spec}{\mathrm{Spec}}
\newcommand{\Spt}{\mathrm{Spt}}
\newcommand{\X}{\cat{X}}
\newcommand{\Y}{\cat{Y}}
\newcommand{\Z}{\cat{Z}}
\newcommand{\ZZ}{\num{Z}}

\newcommand{\acc}{\mathrm{acc}}
\newcommand{\acck}[1]{{#1\text{-}\mathrm{acc}}}
\newcommand{\bc}{\beta}
\renewcommand{\c}{\mathrm{c}}
\newcommand{\cCAlg}{\mathrm{cCAlg}}
\newcommand{\can}{\mathrm{can}}
\newcommand{\cir}{\mathrm{S}^1}
\newcommand{\clspc}[1]{\mathrm{B}#1}
\newcommand{\cofib}{\operatorname*{cofib}}
\newcommand{\colim}{\operatorname*{colim}}
\newcommand{\comp}{\mu}
\newcommand{\coun}{\epsilon}
\newcommand{\diag}[1]{\Delta_{#1}}
\renewcommand{\diamond}{\diamondsuit}
\newcommand{\dstar}{{2*}}
\newcommand{\dual}[1]{#1^\vee}
\newcommand{\ex}{\xi}
\newcommand{\fib}{\operatorname*{fib}}
\newcommand{\h}{\mathrm{h}}

\newcommand{\id}{\mathrm{id}}
\newcommand{\lp}[1]{p_{0,#1}}

\newcommand{\laxdash}[1]{{#1\text{-}\mathrm{lax}}}

\newcommand{\nil}[1]{{#1\text{-}\mathrm{nil}}}
\renewcommand{\o}{\overline}
\newcommand{\op}{\mathrm{op}}
\newcommand{\pr}{\pi}
\newcommand{\rp}[1]{p_{1,#1}}
\newcommand{\suavetwist}{\omega}
\renewcommand{\t}{\mathrm{t}}
\newcommand{\til}{\widetilde}
\newcommand{\poincaretwist}{\nu}
\newcommand{\primtwist}{\delta}
\renewcommand{\u}{\underline}
\newcommand{\uMap}{\u{\smash{\Map}}}
\newcommand{\un}{\eta}
\newcommand{\unit}{\mathbb{1}}
\newcommand{\zstar}{{0*}}

\title{Notes on Tate cohomology}
\author{Arpon Raksit}
\date{\today}

\begin{document}
\maketitle

\begin{abstract}
  We formulate a definition of Tate cohomology in the context of three functor formalisms, and we establish basic monoidality and functoriality properties of it in this context. Our approach to these properties is based on the treatment of Nikolaus--Scholze in the setting of local systems of spectra on spaces. We discuss a couple of other specific settings of interest that are accommodated by our generalization.
\end{abstract}

{\setcounter{tocdepth}{2}
  \footnotesize\tableofcontents\bigskip}


\addtocounter{section}{-1}
\section{Introduction}
\label{i}

In his study of class field theory, Tate \cite{tate--cft} introduced a certain amalgamation of the homology and cohomology of finite groups. This construction, which now goes by the name of \emph{Tate cohomology}, has turned out to have wider scope and relevance. For instance, regarding its scope, generalizations were developed for groups of finite virtual cohomological dimension by Farrell \cite{farrell--tate}, for compact Lie groups by Greenlees--May \cite{greenlees-may--tate}, and then for all topological groups by Klein \cite{klein--dualizing}. Regarding its relevance, let us mention as one example that the case of the circle group $\cir$ plays a crucial role in the theory of (topological) cyclic homology, going back to work of Jones \cite{jones--cyclic} and B\"okstedt--Madsen \cite{bokstedt-madsen--TCZ}, and more recently put into focus by work of Nikolaus--Scholze \cite{nikolaus-scholze--TC} and Bhatt--Morrow--Scholze \cite{bms2}.

This paper gives an exposition of the definition of Tate cohomology and some of its fundamental features in an abstract context, generalizing those mentioned above and accommodating at least a couple of other specific situations that have arisen in other work.

In \cref{i--bg} below, we will review how the theory looks in the topological setting of Klein, highlighting some aspects of interest. In \cref{i--ov}, we will briefly motivate and outline our account.


\subsection{Background}
\label{i--bg}

Let $\Spc$ and $\Spt$ denote the $\infty$-categories of spaces and of spectra, respectively. Let $T$ be a space, and let $L$ be a local system of spectra on $T$, i.e. a functor $T \to \Spt$. This diagram in $\Spt$ admits a limit and a colimit: we set $\Ch^*(T;L) := \lim_T L$ and $\Ch_*(T;L) := \colim_T L$, and we refer to these as the \emph{cohomology} and \emph{homology} of $T$ with coefficients in $L$. When $L$ is the constant local system with value $E \in \Spt$, these recover the spectra $E^T \iso \uMap_\Spt(\Sigma^\infty_+T,E)$ and $E[T] \iso E \otimes \Sigma^\infty_+T$, whose homotopy groups are the cohomology and homology groups of $T$ with coefficients in $E$.

By work of Klein \cite{klein--dualizing}, later revisited by Nikolaus--Scholze \cite{nikolaus-scholze--TC}, these two constructions are related by a canonical \emph{norm map},
\begin{equation}
  \label{i--bg--norm}
  \Nm_T : \Ch_*(T; L \otimes \poincaretwist_T) \to \Ch^*(T;L),
\end{equation}
where $\poincaretwist_T$ is a certain local system of spectra on $T$ (independent of $L$) and $\otimes$ denotes the pointwise tensor product of local systems of spectra on $T$.

\medskip

This general relation specializes to the classical relation of Poincar\'e duality between cohomology and homology, as follows. Following Spivak \cite{spivak--poincare} and Wall \cite{wall--poincare}, the space $T$ is said to be \emph{Poincar\'e} if the following two conditions are satisfied:
\begin{enumerate}
\item $T$ is compact (as an object of $\Spc$, equivalently finitely dominated);
\item $\poincaretwist_T$ is invertible with respect to tensor product (equivalently, the value of $\poincaretwist_T$ at each point $t \in T$ is equivalent to a shift of the sphere spectrum $\SS$).
\end{enumerate}
The first condition guarantees that $\Nm_T$ is an equivalence, allowing us to express the cohomology of $S$, with arbitrary coefficients $L$, in terms of homology. The second condition allows us to reverse the flow: we may replace $L$ by $L \otimes \poincaretwist_T^{-1}$ to obtain a map $\Nm_T : \Ch_*(T;L) \to \Ch^*(T;L \otimes \poincaretwist_T^{-1})$. When both conditions are satisfied, the homology and cohomology of $T$ can each be recovered in terms of the other; this is Poincar\'e duality, in an abstract form.

This abstract formulation of Poincar\'e duality is connected to the theory of manifolds by the following assertion: the space underlying any compact topological manifold is Poincar\'e. Using this fact and the discussion above in the case that $L$ is the constant local system with value $\ZZ$, we recover classical Poincar\'e duality for such manifolds ($L \otimes \poincaretwist_T$ recovering the classical orientation local system in this case).

\medskip

The norm map \cref{i--bg--norm} is also of interest in situations in which it is not an equivalence: we let $\Ct^*(T;L)$ denote its cofiber, and refer to this as the \emph{Tate cohomology} of $T$ with coefficients in $L$.

The terminology here is rooted in the special case where $T$ is the classifying space $\clspc{G}$ of a finite group $G$; let's discuss this case for a bit. Let $t \in \clspc{G}$ denote the canonical basepoint and let $M$ be the value of the local system $L$ at $t$. Then $L$ is precisely the data of an action of $G$ on $M$, and the cohomology $\Ch^*(\clspc{G};L)$ and homology $\Ch_*(\clspc{G};L)$ are alternatively referred to as the \emph{homotopy fixed points} $M^{\h G}$ and \emph{homotopy orbits} $M_{\h G}$.

Note that $\clspc{G}$ is not compact when $G$ is nontrivial, and indeed $\Nm_{\clspc{G}}$ is not in general an equivalence. On the other hand, $\poincaretwist_{\clspc{G}}$ is invertible: in fact, it is canonically trivial, i.e. equivalent to the constant local system with value $\SS$. The norm map thus takes the form $\Nm_{\clspc{G}} : M_{\h G} \to M^{\h G}$, and, as its name suggests, it is induced by the map $\sum_{g \in G} g : M \to M$. The cofiber of the norm map, $\Ct^*(\clspc{G};L)$, is alternatively denoted $M^{\t G}$.

In the case that $M$ is a discrete abelian group, the homotopy groups of $M_{\h G}$ and $M^{\h G}$ are the group homology and cohomology groups of $G$ with coefficients in $M$, and the homotopy groups of $M^{\t G}$ are the original Tate cohomology groups of $G$ with coefficients in $M$.

\medskip

We now return to discussing a general space $T$. As we said at the outset, the object $\Ct^*(T;L)$ is defined as a kind of amalgamation of homology and cohomology, so what inclines us to call it a type of cohomology? To answer this, let's recall two distinctions between cohomology and homology:
\begin{enumerate}
\item \emph{Functoriality}: In the variable $T$, cohomology exhibits a contravariant functoriality while homology exhibits a covariant functoriality: in particular, for $g : U \to T$ a map of spaces, letting $g^*(L) := L \circ g$ be the pullback of the local system $L$ to $U$, we have canonical maps
  \[
    \Ch^*(T;L) \to \Ch^*(U;g^*(L)), \qquad
    \Ch_*(U;g^*(L)) \to \Ch_*(T;L).
  \]
\item \emph{Monoidality}: In the variable $L$, cohomology is lax symmetric monoidal, while homology is oplax symmetric monoidal: in particular, given another local system $L' : T \to \Spt$, we have canonical maps
  \[
    \Ch^*(T;L) \otimes \Ch^*(T;L') \to \Ch^*(T; L \otimes L'), \qquad
    \Ch_*(T;L \otimes L') \to \Ch_*(T;L) \otimes \Ch_*(T;L').
  \]
  This is the source of the multiplicativity of cohomology, i.e. the fact that a ring structure on $L$ induces a ring structure on $\Ch^*(T;L)$.
\end{enumerate}

In light of these discrepancies, the definition of $\Ct^*(T;L)$ does not make it clear either what type of functoriality it may exhibit with respect to the variable $T$ or what type of monoidality it may exhibit with respect to the variable $L$. What ends up happening, though, is that, in some natural situations, it behaves like cohomology in both of these respects, and moreover the canonical map $\Ch^*(T;L) \to \Ct^*(T;L)$ respects these structures.

These features of Tate cohomology can be found in particular when considering the case that $T = \clspc{G}$ for $G$ a finite group. Indeed, it goes back to the origins of the theory that, in this case, the contravariant functoriality of cohomology with respect to \emph{injective} maps of finite groups $H \to G$ and the multiplicativity of cohomology with respect to the coefficient system both extend to Tate cohomology. An account of this in the language we are using here was given by Nikolaus--Scholze \cite{nikolaus-scholze--TC}.

Nikolaus--Scholze explain the presence of these features of Tate cohomology by characterizing it, together with the norm map, by a certain universal property. A very similar characterization was established previously by Klein \cite{klein--axioms}. Klein's original work \cite{klein--dualizing} also shows that the key condition required in Nikolaus--Scholze's explanation of the monoidality of Tate cohomology is satisfied more generally when $T = \clspc{G}$ for $G$ the group space underlying a compact Lie group. As mentioned earlier, the case $G = \cir$ is also important in \cite{nikolaus-scholze--TC}.

In fact, a closer examination of Klein's analysis shows that what is actually relevant for the aforementioned condition to be satisfied is that $G$ is Poincar\'e (just as a space, in the sense discussed above). That is, one kind of special behavior of the Tate cohomology of $G$ is tied to another kind of special behavior of the Tate cohomology of $\clspc{G}$. And, as one final point of interest, it turns out that any compact group space is Poincar\'e; this result, for which only a very weak form of group structure is really needed, is due to Browder \cite{browder--torsion}.\footnote{It's furthermore true that any compact group space admits a smooth manifold structure \cite{bknp--loop} (albeit not necessarily a Lie group structure \cite{hilton-roitberg--bundles}, even up to rational equivalence \cite{abgp--rational-lie}), but this type of result is not our subject here.}


\subsection{Overview}
\label{i--ov}

Our goal in this paper is to abstract the theory described in \cref{i--bg}, so that we may apply the same ideas in other, formally similar situations, like the following.

\begin{example}
  \label{i--ov--l}
  We might consider local systems on spaces valued in an $\infty$-category $\C$ other than $\Spt$. A case of interest is where $\C$ is a chromatic localization of $\Spt$: in this setting, there is the phenomenon of ``ambidexterity''---developed in work of Greenlees--Sadofsky \cite{greenlees-sadofsky--tate}, Hovey--Sadofsky \cite{hovey-sadofsky--blueshift}, Kuhn \cite{kuhn--blueshift}, Hopkins--Lurie \cite{hopkins-lurie--ambi}, and Carmeli--Schlank--Yanovski \cite{csy--ambi}---which we may regard as supplying new examples of Poincar\'e duality. We would like to derive from this new examples of functoriality in Tate cohomology; for instance, this is relevant in an analysis of the topological cyclic homology of $\ZZ$ carried out in joint work with Devalapurkar \cite{devalapurkar-r--THHZ}.
\end{example}

\begin{example}
  \label{i--ov--q}
  We might consider quasicoherent sheaves on certain kinds of stacks, in the sense of algebraic geometry. For instance, in place of representations of a group (equivalent to local systems over its classifying space), we may be interested in representations of a group scheme or group stack (equivalent to quasicoherent sheaves on its classifying stack). A situation of this type arises in the study of the ``filtered circle'' in \cite{antieau-riggenbach--synthetic,raksit--hhdr}, and, relatedly, in the study of certain variants of Hochschild homology introduced by Moulinos--Robalo--To\"en \cite{mrt--hkr}.
\end{example}

The setting for our abstraction will be that of a \emph{three functor formalism}, as defined by Mann \cite{mann--p-adic-six,heyer-mann--six}. We will review this notion in \cref{t--f}, in particular indicating in what sense it is a framework for reasoning about homology and cohomology of topological or geometric objects. The remainder of \cref{t} is devoted firstly to generalizing the definition of the norm map and Tate cohomology to this setting, and then to developing the features of these that were highlighted in \cref{i--bg}, involving Poincar\'e duality, functoriality, and monoidality. Then, in \cref{e}, we will explain how \cref{i--ov--l,i--ov--q} are accommodated by this framework.


\subsection{Acknowledgements}
\label{i--ac}

First and foremost, I thank Lars Hesselholt, for explaining to me the relevance of the recent developments around three/six functor formalisms, especially concerning $(\infty,2)$-categories of kernels, to the study of norm maps. I am also thankful that all of this is now thoroughly documented, particularly in the work of Heyer--Mann \cite{heyer-mann--six}; for consistency and convenience, we will stick to referencing this source for these matters, but see for instance the introduction there for more on the history of these ideas. Perhaps this a good moment to stress the expository nature of this paper: we will largely be reviewing and lightly mixing ideas explained in \cite{heyer-mann--six} and the work of Nikolaus--Scholze \cite{nikolaus-scholze--TC}.

In addition, I am grateful to Sanath Devalapurkar, for the collaboration resulting in \cite{devalapurkar-r--THHZ}, alongside which the present work developed; and to Ben Antieau, Jacob Lurie, Maxime Ramzi, and Allen Yuan, for helpful conversations related to this material.

This work was supported by the Bell Systems Fellowship at IAS, NSF grant DMS-2103152, and the Simons Collaboration on Perfection. It took me too long and I appreciated the time.


\subsection{Conventions}
\label{i--cv}

\begin{enumerate}[leftmargin=*]
\item \label{i--cv--cat}
  We let $\Cat_\infty$ denote the extra large $\infty$-category of possibly large $\infty$-categories, $\PrL$ the subcategory of $\Cat_\infty$ spanned by the presentable $\infty$-categories and colimit preserving functors, and $\PrLst$ the full subcategory of $\PrL$ spanned by the stable, presentable $\infty$-categories.
\item \label{i--cv--2cat}
  By an $(\infty,2)$-category we mean a (possibly large) $\Cat_\infty$-enriched $\infty$-category. We refer to \cite[Appendices C and D]{heyer-mann--six} for a convenient overview of the theories of enriched $\infty$-categories and $(\infty,2)$-categories. We will use the following notation in the context of an $(\infty,2)$-category $\AA$:
  \begin{enumerate}
  \item \label{i--cv--2cat--map}
    For $\A,\B \in \AA$, we let $\uMap_\AA(\A,\B)$ denote the associated mapping $\infty$-category.
  \item \label{i--cv--2cat--op}
    We let $\twocat{A}^\op$ and $\twocat{A}^\c$ denote the opposite and conjugate $(\infty,2)$-categories of $\twocat{A}$, respectively: these have the same objects as $\twocat{A}$ and mapping $\infty$-categories given by
    \[
      \uMap_{\twocat{A}^\op}(\A,\B) \iso \uMap_{\twocat{A}}(\B,\A), \quad
      \uMap_{\twocat{A}^\c}(\A,\B) \iso \uMap_{\twocat{A}}(\A,\B)^\op.
    \]
  \end{enumerate}
\end{enumerate}


\section{Generalities}
\label{t}

Throughout this section, we let $\D$ be a three functor formalism on a geometric setup $(\S,\S_!)$, in the sense of Heyer-Mann's work \cite{heyer-mann--six}. What this means will be reviewed at more length in \cref{t--f}; for the moment, let us just say that $\S$ is an $\infty$-category whose objects we think of as ``spaces'', and that $\D$ assigns to each $S \in \S$ an $\infty$-category $\D(S)$ whose objects we think of as ``coefficient systems'' for homology and cohomology. Our purpose in this section is to formulate and study a notion of Tate cohomology in this setting.

As was indicated in \cref{i--bg}, Tate cohomology is closely tied to the phenomenon of Poincar\'e duality. One of the key points of \cite{heyer-mann--six} is to express Poincar\'e duality in the context of a three functor formalism in terms of the abstract notion of an adjunction in an $(\infty,2)$-category; more specifically, they consider this general notion in the case of an \emph{$(\infty,2)$-category of kernels} $\K_\D(S)$, which is constructed for each $S \in \S$ using $\D$ (also to be reviewed in \cref{t--f}). What we will spell out here, in \cref{t--d}, is that the same ideas may be applied to define norm maps and Tate cohomology in this context. For this, what is relevant is not quite the notion of adjunction, but a certain weakening thereof; we will discuss this general $(\infty,2)$-categorical notion in \cref{t--a}. Then, in \cref{t--p}, we will review the description of Poincar\'e duality from \cite{heyer-mann--six}, slightly rephrasing some things in terms of our norm maps, and in \cref{t--h}, we will prove in this context a version of the result of Browder mentioned in \cref{i--bg}.

Following that, the main aim of the section is to follow Nikolaus--Scholze's work \cite{nikolaus-scholze--TC} to establish, in certain circumstances, a universal property characterizing our construction of Tate cohomology, and to then deduce functoriality and monoidality properties of the type described in \cref{i--bg}. We will come to these results in \cref{t--u} and \cref{t--n}. Before that, in \cref{t--b} and \cref{t--c}, we will need to articulate the compatibility of norm maps with base change and composition in $\S$.


\subsection{Three functor formalisms}
\label{t--f}

As stated above, we have fixed a geometric setup $(\S,\S_!)$ and a three functor formalism $\D$ on it. This subsection is devoted to reviewing the features of this general situation that will be relevant for our discussion, and setting notation and terminology.

\begin{enumerate}[leftmargin=*]
\item \label{t--f--S}
  $\S$ is an $\infty$-category, whose objects we think of as some kind of ``spaces''. $\S_!$ is a wide subcategory of $\S$, and we say that a map in $\S$ is \emph{$\D$-$!$-able} if it lies in $\S_!$. The $\D$-$!$-able maps are assumed to admit and be closed under base changes along arbitrary maps in $\S$, and to be closed under formation of diagonals. In this subsection, to simplify the exposition, we will assume that $\S$ admits finite products; but this assumption will not be needed in the remainder of the section.
\item \label{t--f--span}
  There is an associated $\infty$-category $\Span(\S,\S_!)$, called an \emph{$\infty$-category of spans} (or alternatively \emph{of correspondences}):
  \begin{enumerate}
  \item \label{t--f--span--obj}
    The objects of $\Span(\S,\S_!)$ are the objects of $\S$.
  \item \label{t--f--span--map}
    For $S,T \in \S$, a map from $S$ to $T$ in $\Span(\S,\S_!)$ is given by a span $\smash{S \lblfrom{f} U \lblto{g} T}$ in $\S$ in which $g$ is $\D$-$!$-able. Composition of spans involves forming fiber products.
  \item \label{t--f--span--iota}
    We have functors $\iota_0 : \S^\op \to \Span(\S,\S_!)$ and $\iota_1 : \S_! \to \Span(\S,\S_!)$. On objects, both are given by the identity map. On maps, $\iota_0$ sends a map $f : U \to S$ in $\S$ to the span $\smash{S \lblfrom{f} U \lblto{\id_U} U}$ in $\S$, and $\iota_1$ sends a $\D$-$!$-able map $g : U \to T$ in $\S$ to the span $\smash{U \lblfrom{\id_U} U \lblto{g} T}$ in $\S$. A general map $\smash{S \lblfrom{f} U \lblto{g} T}$ in $\Span(\S,\S_!)$ can then be written as the composition $\iota_1(g)\iota_0(f)$.
  \item \label{t--f--span--algebra}
    There is a canonical symmetric monoidal structure on $\Span(\S,\S_!)$, together with one on the functor $\iota_0 : \S^\op \to \Span(\S,\S_!)$, where $\S^\op$ is equipped with its cocartesian symmetric monoidal structure. In particular, the tensor product of $\Span(\S,\S_!)$ is given on objects by cartesian product in $\S$, and the unit object is the final object $* \in \S$. Furthermore, each $S \in \Span(\S,\S_!)$ carries a canonical commutative algebra structure, the multiplication map being $\iota_0(\diag{S})$, i.e. the span $\smash{S \times S \lblfrom{\diag{S}} S \lblto{\id_S} S}$, where $\diag{S}$ is the diagonal map, and the unit map being $\iota_0(p)$, i.e. the span $\smash{* \lblfrom{p} S \lblto{\id_S} S}$, where $p$ is the unique such map.
  \item \label{t--f--span--frobenius}
    Let $S \in \S$ and let $p : S \to *$ be the unique such map. If $p$ is $\D$-$!$-able, then $S$ is moreover a commutative Frobenius algebra in $\Span(\S,\S_!)$: the map $\iota_1(p)$, i.e. the span $\smash{S \lblfrom{\id_S} S \lblto{p} *}$, is nondegenerate, in the sense that the composition $\iota_1(p)\iota_0(\diag{S})$, i.e. the span $\smash{S \times S \lblfrom{\diag{S}} S \lblto{p} *}$, is a duality datum in $\Span(\S,\S_!)$.
  \end{enumerate}
\item \label{t--f--D}
  $\D$ is a lax symmetric monoidal functor $\Span(\S,\S_!) \to \Cat_\infty$:
  \begin{enumerate}
  \item \label{t--f--D--D}
    To each object $S \in \S$ is assigned a symmetric monoidal $\infty$-category $\D(S)$. We denote the tensor product operation in $\D(S)$ simply by $- \otimes - $ and we denote the unit object by $\unit_S$. We think of the objects of $\D(S)$ as some kind of ``coefficient systems'' or ``sheaves'' on $S$. The main results below concern the case where $\D(S)$ is in fact a stable, presentable symmetric monoidal $\infty$-category, but many of the preliminary pieces make sense more generally, and it will be good to retain the symmetry of passing to opposites where possible; we will state additional hypotheses when they are necessary (cf. \cref{t--f--D--presentable} below).
  \item \label{t--f--D--pullback}
    To each map $f : T \to S$ in $\S$ is assigned a symmetric monoidal functor $f^* : \D(S) \to \D(T)$, namely the image under $\D$ of $\iota_0(f)$; we refer to $f^*$ as \emph{$*$-pullback} along $f$.

    We say that $f$ is \emph{$\D$-$*$-able} if $f^*$ admits a right adjoint, in which case we denote the right adjoint by $f_* : \D(T) \to \D(S)$ and refer to this as \emph{$*$-pushforward} or \emph{cohomology} along $f$. We say that $f$ is \emph{$\D$-$\sharp$-able} if $f^*$ admits a left adjoint, in which case we denote the left adjoint by $f_\sharp : \D(T) \to \D(S)$ and refer to this as \emph{$\sharp$-pushforward} or \emph{homology} along $f$.

    If $f$ is $\D$-$*$-able, then the functor $f_*$ inherits a lax symmetric monoidal structure, and hence a lax $\D(S)$-linear structure, from $f^*$. For instance, we have in this case a natural \emph{projection map}
    \[
      \pr^*_f : f_*(M) \otimes L \to f_*(M \otimes f^*(L))
    \]
    for $M \in \D(T)$ and $L \in \D(S)$; we will say that $f$ satisfies the \emph{$\D$-$*$-projection formula} if this is a natural isomorphism, i.e. if $f_*$ is in fact strictly $\D(S)$-linear. Similarly, if $f$ is $\D$-$\sharp$-able, then the functor $f_\sharp$ inherits an oplax symmetric monoidal structure and oplax $\D(S)$-linear structure; again we have a natural \emph{projection map}
    \[
      \pr^\sharp_f : f_\sharp(M \otimes f^*(L)) \to f_\sharp(M) \otimes L
    \]
    for $M \in \D(T)$ and $L \in \D(S)$, and we will say that $f$ satisfies the \emph{$\D$-$\sharp$-projection formula} if this is a natural isomorphism, i.e. if $f_\sharp$ is strictly $\D(S)$-linear.

    A commutative square $\sigma :$
    \[
      \begin{tikzcd}
        V \ar[r, "q"] \ar[d, "g", swap] &
        T \ar[d, "f"] \\
        U \ar[r, "p"] &
        S,
      \end{tikzcd}
    \]
    in $\S$ induces an isomorphism $q^*f^* \iso g^*p^*$ of functors from $\D(S)$ to $\D(V)$. This then induces \emph{base change maps}
    \[
      \bc^*_\sigma : p^*f_* \to g_*q^*, \quad\quad
      \bc^\sharp_\sigma : g_\sharp q^* \to p^* f_\sharp,
    \]
    assuming that $f$ and $g$ are both $\D$-$*$-able or $\D$-$\sharp$-able, respectively: namely, the compositions
    \[
      p^*f_* \lblto{\un^*_g} g_*g^*p^*f_* \iso g_*q^*f^*f_* \lblto{\coun^*_f} g_*q^*, \quad\quad
      g_\sharp q^* \lblto{\un^\sharp_f} g_\sharp q^* f^* f_\sharp \iso g_\sharp g^*p^* f_\sharp \lblto{\coun^\sharp_g} p^* f_\sharp,
    \]
    where the labelled maps are the unit and counit maps for the relevant adjunctions. We will only consider these base change maps in the case that $\sigma$ is cartesian. If $f$ is $\D$-$*$-able, we will say that $f$ satisfies \emph{$\D$-$*$--base change} if for every cartesian square $\sigma$ as above, $g$ is also $\D$-$*$-able and the base change map $\bc^*_\sigma$ is an isomorphism. Similarly, if $f$ is $\D$-$\sharp$-able, we will say that $f$ satisfies \emph{$\D$-$\sharp$--base change} if for every cartesian square $\sigma$ as above, $g$ is also $\D$-$\sharp$-able and the base change map $\bc^\sharp_\sigma$ is an isomorphism.
  \item \label{t--f--D--shriek}
    To each $\D$-$!$-able map $f : T \to S$ in $\S$ is assigned a functor $f_! : \D(T) \to \D(S)$, namely the image under $\D$ of $\iota_1(f)$; we refer to this as \emph{$!$-pushforward} along $f$. This comes equipped with a canonical $\D(T)$-linear structure, including a natural \emph{projection isomorphism}
    \[
      \pr^!_f : f_!(L) \otimes M \iso f_!(L \otimes f^*(M))
    \]
    for $L \in \D(T)$ and $M \in \D(S)$. In addition, given a cartesian square $\sigma :$
    \[
      \begin{tikzcd}
        V \ar[r, "q"] \ar[d, "g", swap] &
        U \ar[d, "f"] \\
        T \ar[r, "p"] &
        S
      \end{tikzcd}
    \]
    in $\S$, there is a canonical \emph{base change isomorphism}
    \[
      \bc^!_\sigma : p^*f_! \iso g_!q^*.
    \]
    This then induces \emph{exchange maps}
    \[
      \ex^*_\sigma : f_!q_* \to p_*g_!, \quad\quad
      \ex^\sharp_\sigma : p_\sharp g_! \to f_! q_\sharp,
    \]
    assuming that $p$ and $q$ are both $\D$-$*$-able or $\D$-$\sharp$-able, respectively: namely, the compositions
    \[
      f_!q_* \lblto{\un^*_p} p_*p^*f_!q_* \lbliso{\bc^!_\sigma} p_*g_!q^*q_* \lblto{\coun^*_q} g_*q^*, \quad\quad
      p_\sharp g_! \lblto{\un^\sharp_q} p_\sharp g_!q^*q_\sharp \lbliso{\bc^!_\sigma} p_\sharp p^*f_!q_\sharp \lblto{\coun^\sharp_p} f_!q_\sharp.
    \]
  \item \label{t--f--D--presentable}
    We say that $\D$ is presentable (resp. stable and presentable) if it promotes to a lax symmetric monoidal functor from $\Span(\S,\S_!)$ to $\PrL$ (resp. $\PrLst$). This is equivalent to the following properties being satisfied: for each $S \in \S$, the $\infty$-category $\D(S)$ is presentable (resp. stable and presentable) and its tensor product preserves colimits in each variable; for each map $f : T \to S$ in $\S$, the functor $f^* : \D(S) \to \D(T)$ preserves colimits, which by the adjoint functor theorem implies that it admits a right adjoint, i.e. that $f$ is $*$-able; for each $\D$-$!$-able map $f : T \to S$ in $\S$, the functor $f_! : \D(T) \to \D(S)$ preserves colimits, hence admits a right adjoint, which we denote by $f^! : \D(S) \to \D(T)$ and refer to as \emph{$!$-pullback} along $f$. (In this case, $\D$ is more often referred to as a six functor formalism.)
  \item \label{t--f--D--cover}
    Let $S \in \S$, and let $\cat{U}$ be a collection of maps in $\S$ with target $S$. Let $\cat{U}'$ be the sieve generated by $\cat{U}$. We say that $\cat{U}$ is a \emph{$\D$-$*$-cover} if, for every map $s : S' \to S$ in $\S$, the composite functor
    \[
      \S^\op \lblto{\iota_0} \Span(\S,\S_!) \lblto{\D}  \Cat_\infty
    \]
    satisfies descent along the pullback sieve $s^*(\cat{U}')$, i.e. the functor
    \[
      \D(S') \to \lim_{U \in s^*(\cat{U}')^\op} \D(U)
    \]
    induced by the $*$-pullback functoriality of $\D$ is an equivalence. Assuming that $\D$ is presentable, we say that $\cat{U}$ is a \emph{$\D$-$!$-cover} if, for every map $s : S' \to S$ in $\S$, the composite functor
    \[
      \S \lblto{\iota_1} \Span(\S,\S_!) \lblto{\D}  \PrL
    \]
    satisfies codescent along the pullback sieve $s^*(\cat{U}')$, i.e. the functor
    \[
      \colim_{U \in s^*(\cat{U}')^\op} \D(U) \to \D(S')
    \]
    induced by the $!$-pushforward functoriality of $\D$ is an equivalence; here the colimit is taken in $\PrL$, so the condition is equivalent to the functor
    \[
      \D(S') \to \lim_{U \in s^*(\cat{U}')^\op} \D(U)
    \]
    induced by the $!$-pullback functoriality of $\D$ being an equivalence.
  \end{enumerate}
\item \label{t--f--K}
  For each $S \in \S$, there is an associated $(\infty,2)$-category $\K_\D(S)$, called an \emph{$(\infty,2)$-category of kernels}:
  \begin{enumerate}
  \item \label{t--f--K--obj}
    The objects of $\K_\D(S)$ are the objects of $(\S_!)_{/S}$, i.e. $\D$-$!$-able maps $f : T \to S$ in $\S$. Here, when we regard such a map as an object of $\K_\D(S)$, we will denote it by $[T]$ or $[T]_S$ (allowing the standard abuse of notation of leaving the map $f$ implicit).
  \item \label{t--f--K--map}
    For $T_0,T_1 \in (\S_!)_{/S}$, the mapping $\infty$-category $\uMap_{\K_\D(S)}([T_0],[T_1])$ is given by $\D(T_0 \times_S T_1)$. Similar to the previous point, for $L_{01} \in \D(T_0 \times_S T_1)$, we will denote the corresponding map in $\K_\D(S)$ by $[L_{01}] : [T_0] \to [T_1]$ or $[L_{01}]_S : [T_0]_S \to [T_1]_S$. Composition in $\K_\D(S)$ is given by ``convolution'': for $T_0,T_1,T_2 \in (\S_!)_{/S}$, $L_{01} \in \D(T_0 \times_S T_1)$, and $L_{12} \in \D(T_1 \times_S T_2)$, we have a natural isomorphism
    \[
      [L_{12}] \circ [L_{01}] \iso [(p_{02})_!(p_{01}^*(L_{01}) \otimes p_{12}^*(L_{12}))],
    \]
    where $p_{ij}$ denotes the projection map $T_0 \times_S T_1 \times_S T_2 \to T_i \times_S T_j$.
  \item \label{t--f--K--id}
    For $T \in (\S_!)_{/S}$, the identity map $\id_{[T]} : [T] \to [T]$ in $\K_\D(S)$ is given by $[(\diag{T})_!(\unit_T)]$, where $\diag{T} : T \to T \times_S T$ is the diagonal map.
  \item \label{t--f--K--Psi}
    We let $\Psi_{\D,S} : \K_\D(S) \to \Cat_\infty$ denote the mapping functor $\smash{\uMap_{\K_\D(S)}([S],-)}$. This is a functor of $(\infty,2)$-categories. For $T \in (\S_!)_{/S}$, it sends the object $[T] \in \K_\D(S)$ to $\D(T) \in \Cat_\infty$; and for $T_0,T_1 \in (\S_!)_{/S}$ and $L_{01} \in \D(T_0 \times_S T_1)$, it sends the map $[L_{01}] : [T_0] \to [T_1]$ in $\K_\D(S)$ to the functor
    \[
      \Psi_{\D,S}([L_{01}]) : \D(T_0) \to \D(T_1), \quad L_0 \mapsto (p_1)_!(p_0^*(L_0) \otimes L_{01}),
    \]
    where $p_i$ denotes the projection map $T_0 \times_S T_1 \to T_i$.

    We note also that $\Psi_{\D,S}$ factors canonically through the $(\infty,2)$-category $\Mod_{\D(S)}(\Cat_\infty)$ of $\D(S)$-linear $\infty$-categories; this follows from the characterization of $\Psi_{\D,S}$ appearing in the proof of \cite[Proposition 4.1.5]{heyer-mann--six}, and the fact that $\D$ induces a lax symmetric monoidal functor $\Span((\S_!)_{/S}) \to \Mod_{\D(S)}(\Cat_\infty)$. In particular, the functor $\Psi_{\D,S}([L_{01}])$ displayed above is canonically $\D(S)$-linear, and any map $L_{01} \to L_{01}'$ in $\D(T_0 \times_S T_1)$ induces a $\D(S)$-linear natural transformation $\Psi_{\D,S}([L_{01}]) \to \Psi_{\D,S}([L_{01}'])$.
  \item \label{t--f--K--op}
    The opposite and conjugate $(\infty,2)$-categories $\K_\D(S)^\op$ and $\K_\D(S)^\c$ (see \itemitemref{i--cv}{2cat}{op}) may be described as follows: we have an equivalence $\K_\D(S)^\op \iso \K_\D(S)$, which acts by the identity on objects and by the canonical equivalence
    \[
      \uMap_{\K_\D(S)}([T_0],[T_1]) \iso \D(T_0 \times_S T_1) \iso \D(T_1 \times_S T_0) \iso \uMap_{\K_\D(S)}([T_1],[T_0])
    \]
    on mapping $\infty$-categories \cite[Proposition 4.1.4]{heyer-mann--six}; and by definition we have an equivalence $\K_\D(S)^\c \iso \K_{\D^\op}(S)$, where $\D^\op$ denotes the opposite three functor formalism of $\D$, i.e. the composition of $\D$ with $(-)^\op : \Cat_\infty \to \Cat_\infty$.
  \item \label{t--f--K--Dprime}
    Suppose given another three functor formalism $\D' : \Span(\S,\S_!) \to \Cat_\infty$ and a map (i.e. lax symmetric monoidal natural transformation) $\alpha : \D \to \D'$. Then there is an induced functor of $(\infty,2)$-categories $\phi_\alpha : \K_\D(S) \to \K_{\D'}(S)$ which acts by the identity on objects and by $\alpha$ on mapping $\infty$-categories \cite[Proposition 4.2.1]{heyer-mann--six}.
  \end{enumerate}
\item \label{t--f--K-monoidal}
  For each $S \in \S$, there is a canonical symmetric monoidal structure on $\K_\D(S)$, together with a symmetric monoidal functor $\Phi_{\D,S} : \Span((\S_!)_{/S}) \to \K_\D(S)$ \cite[Theorem 4.2.4]{heyer-mann--six}:
  \begin{enumerate}
  \item \label{t--f--K-monoidal--span}
    Here $\Span((\S_!)_{/S})$ is an $\infty$-category of spans defined similarly to $\Span(\S,\S_!)$, but beginning with the $\infty$-category $(\S_!)_{/S}$ in place of $\S$, and with no extra restrictions needed to be placed on the legs of the spans (note that, for $T_0,T_1 \in (\S_!)_{/S}$, any map $T_0 \to T_1$ in $\S_{/S}$ is $\D$-$!$-able, i.e. is in fact a map in $(\S_!)_{/S}$).
  \item \label{t--f--K-monoidal--Phi}
    The functor $\Phi_{\D,S}$ acts by the identity on objects, i.e. sends $T \in \Span((\S_!)_{/S})$ to $[T] \in \K_\D(S)$. On maps, $\Phi_{\D,S}$ sends a span $\smash{T_0 \lblfrom{f_0} U \lblto{f_1} T_1}$ in $(\S_!)_{/S}$ to the map $[(f_{01})_!(\unit_U)] : [T_0] \to [T_1]$ in $\K_\D(S)$, where $f_{01} : U \to T_0 \times_S T_1$ is the map induced by $f_0$ and $f_1$.
  \item \label{t--f--K-monoidal--formulas}
    The symmetric monoidality of $\Phi_{\D,S}$ tells us that the tensor product of $\K_\D(S)$ is given on objects by $[T_0] \otimes [T_1] \iso [T_0 \times_S T_1]$, and that the unit object is $[S]$. In addition, it follows from \cref{t--f--span}\cref{t--f--span--algebra,t--f--span--frobenius} that each object $[T] \in \K_\D(S)$ carries a canonical commutative Frobenius algebra structure, with multiplication map $\mu : [T \times_S T] \to [T]$ given by $[(\diag{f}')_!(\unit_T)]$, where $\diag{f}' : T \to T \times_S T \times_S T$ is the diagonal map, and unit map $\eta : [S] \to [T]$ and nondegenerate map $\lambda : [T] \to [S]$ each being given by $[\unit_T]$.
  \end{enumerate}
\item \label{t--f--K-functoriality}
  Similar to $\D$, the assignment $S \mapsto \K_\D(S)$ extends canonically to a lax symmetric monoidal functor $\Span(\S,\S_!) \to \Cat_{(\infty,2)}$ \cite[Theorem 4.2.4]{heyer-mann--six}:
  \begin{enumerate}
  \item \label{t--f--K-functoriality--pull}
    For each map $f : T \to S$ in $\S$, the image of the map $\iota_0(f)$ under this functor is a symmetric monoidal functor $f^\dstar : \K_\D(S) \to \K_\D(T)$. On objects, $f^\dstar$ is given by base change, i.e. it sends $[U]$ to $[U \times_S T]$. It sends a map $[L_{01}] : [U_0] \to [U_1]$ in $\K_\D(S)$ to the map $[q^*(L_{01})] : [U_0 \times_S T] \to [U_1 \times_S T]$ in $\K_\D(T)$, where $q : (U_0 \times_S T) \times_T (U_1 \times_S T) \to U_0 \times_S U_1$ is the projection map.
  \item \label{t--f--K-functoriality--push}
    For each $\D$-$!$-able map $f : T \to S$ in $\S$, the image of $\iota_1(f)$ under this functor is a functor $f_\dstar : \K_\D(T) \to \K_\D(S)$. On objects, $f_\dstar$ is given by composition with $f$, i.e. sends $[U]_T$ to $[U]_S$. It sends a map $[L_{01}] : [U_0]_T \to [U_1]_T$ in $\K_\D(T)$ to the map $[r_!(L_{01})] : [U_0]_S \to [U_1]_S$, where $r : U_0 \times_T U_1 \to U_0 \times_S U_1$ is the canonical such map.
  \item \label{t--f--K-functoriality--adjunction}
    For $f : T \to S$ as in the previous point, $f_\dstar$ is both right and left adjoint to $f^\dstar$ \cite[Lemma 4.2.7]{heyer-mann--six}. One consequence of the right adjointness is that $f_\dstar$ is lax symmetric monoidal. Another consequence, given the identification $f^\dstar([S]) \iso [T]$, is a natural isomorphism $\Psi_{\D,T} \iso \Psi_{\D,S} \circ f_\dstar$. It follows that the unit map $\id_{\K_\D(S)} \to f_\dstar \circ f^\dstar$ exhibiting $f_\dstar$ as right adjoint to $f^\dstar$ induces a natural transformation $\Psi_{\D,S} \to \Psi_{\D,T} \circ f^\dstar$, and similarly the counit map $f_\dstar \circ f^\dstar \to \id_{\K_\D(S)}$ exhibiting $f_\dstar$ as left adjoint to $f^\dstar$ induces a natural transformation $\Psi_{\D,T} \circ f^\dstar \to \Psi_{\D,S}$; unwinding definitions shows that the former is given on $[U] \in \K_\D(S)$ by the functor $q^* : \D(U) \to \D(U \times_S T)$ and the latter by the functor $q_! : \D(U \times_S T) \to \D(U)$, where $q : U \times_S T \to U$ is the projection map.
  \end{enumerate}
\end{enumerate}


\subsection{Weak adjoints}
\label{t--a}

Recall that the notion of adjunction makes sense within any $(\infty,2)$-category. As was indicated in the opening to this section, it is shown in \cite{heyer-mann--six} that considering this notion within the $(\infty,2)$-categories of kernels $\K_\D(S)$ (\itemref{t--f}{K}) is quite fruitful for understanding Poincar\'e duality in the context of the three functor formalism $\D$. In the next subsection, we will elaborate on that idea, constructing norm maps and Tate cohomology in this setting by considering a certain weakening of the notion of adjunction in the $(\infty,2)$-categories $\K_\D(S)$. This weaker notion also makes sense in a general $(\infty,2)$-category, and the purpose of this subsection is to introduce it and to collect the abstract results about it that will be invoked later in the section.

\begin{definition}
  \label{t--a--weak-adjoint}
  Let $\twocat{A}$ be an $(\infty,2)$-category, and let $G : \A \to \B$ be a map in $\twocat{A}$. A \emph{weak right adjoint} to $G$ is a map $H : \B \to \A$ in $\twocat{A}$ together with a natural isomorphism
  \[
    \alpha : \Map_{\uMap_{\twocat{A}}(\B,\A)}(J,H) \iso \Map_{\uMap_{\twocat{A}}(\B,\B)}(GJ,\id_\B)
  \]
  for $J \in \uMap_{\twocat{A}}(\B,\A)$. A \emph{weak left adjoint} to $G$ is a map $F : \B \to \A$ in $\twocat{A}$ together with a natural isomorphism
  \[
    \alpha : \Map_{\uMap_{\twocat{A}}(\B,\A)}(F,J) \iso \Map_{\uMap_{\twocat{A}}(\B,\B)}(\id_\B, GJ)
  \]
  for $J \in \uMap_{\twocat{A}}(\B,\A)$. We will often leave the natural isomorphism $\alpha$ implicit when referring to weak adjoints.
\end{definition}

\begin{remark}
  \label{t--a--weak-adjoint-props}
  In the situation of \cref{t--a--weak-adjoint}, we note the following:
  \begin{enumerate}
  \item \label{t--a--weak-adjoint-props--adjoint}
    A right adjoint of $G$ is in particular a weak right adjoint of $G$, and a left adjoint of $G$ is in particular a weak left adjoint of $G$.
  \item \label{t--a--weak-adjoint-props--unique}
    The space of weak right adjoints to $G$ is either empty or contractible, and the same for the space of weak left adjoints. We therefore allow ourselves to refer to \emph{the} weak right adjoint or weak left adjoint to $G$ when they exist.
  \item \label{t--a--weak-adjoint-props--symmetry}
    Weak adjunctions share only part of the symmetry of honest adjunctions. On the negative side, passing between $\twocat{A}$ and its opposite $\twocat{A}^\op$ does not exchange weak right and left adjoints: that is, $H$ being weakly right adjoint to $G$ does not imply that $G$ is weakly left adjoint to $H$, and $F$ being weakly left adjoint to $G$ does not imply that $G$ is weakly right adjoint to $F$. On the other hand, passing between $\twocat{A}$ and its conjugate $\twocat{A}^\c$ does exchange weak right and left adjoints; this is what is meant below whenever we say that a result involving weak left adjoints is ``formally dual'' to another result involving weak right adjoints.
  \item \label{t--a--weak-adjoint-props--presentable}
    Let $g : \uMap_{\twocat{A}}(\B,\A) \to \uMap_{\twocat{A}}(\B,\B)$ be the functor sending $J \mapsto GJ$. If $g$ admits a right adjoint $h$ (for example, if $\twocat{A}$ is $\PrL$-enriched), then $G$ admits a weak right adjoint, namely the map $H := h(\id_B)$. Similarly, if $g$ admits a left adjoint $f$, then $F := f(\id_B)$ is a weak left adjoint to $G$.
  \item \label{t--a--weak-adjoint-props--counit}
    As with right adjoints, a weak right adjoint $H$ to $G$ is characterized by an associated \emph{counit map} $\coun : G H \to \id_\B$, namely the image of the identity map $\id_G$ under the natural isomorphism $\alpha$ of the weak adjunction. That is, $\alpha$ may be recovered from $\coun$, as the composition
    \[
      \Map_{\uMap_{\twocat{A}}(\B,\A)}(J,H) \lblto{G \circ -} \Map_{\uMap_{\twocat{A}}(\B,\B)}(GJ,GH) \lblto{\coun} \Map_{\uMap_{\twocat{A}}(\B,\B)}(GJ,\id_\B).
    \]
    Similarly, a weak left adjoint $F$ to $G$ is characterized by an associated \emph{unit map} $\un : \id_\B \to G F$.
  \end{enumerate}
\end{remark}

\begin{example}
  \label{t--a--PrL}
  Fix a presentable monoidal $\infty$-category $\C \in \Alg(\PrL)$, and consider the $(\infty,2)$-category $\AA = \LMod_\C(\PrL)$ of $\C$-linear presentable $\infty$-categories (recall that the $\LMod_\C(\PrL)$ is enriched over $\PrL$ by means of its canonical tensoring over $\PrL$). Let $G : \A \to \B$ be a map in $\LMod_\C(\PrL)$. As a map in $\PrL$, the functor $G$ admits a right adjoint functor $H : \B \to \A$. The functor $H$ inherits a lax $\C$-linear structure from $G$ \cite[Example 7.3.2.8]{lurie--HA}, but this may or may not be a strict $\C$-linear structure, and $H$ may or may not preserve colimits; $H$ promotes to a right adjoint in $\LMod_\C(\PrL)$ if and only if these conditions do hold (see \cite[Remark 7.3.2.9]{lurie--HA}). On the other hand, $G$ always admits a weak right adjoint $\til{H} : \B \to \A$ in $\LMod_\C(\PrL)$: this follows from \itemref{t--a--weak-adjoint-props}{presentable}, since $\LMod_\C(\PrL)$ is $\PrL$-enriched. For further discussion of this situation, see \cref{t--a--assembly} below.
\end{example}

The remainder of this subsection consists of a few general constructions and observations involving weak adjoints. We will only state these for weak right adjoints, but there are formally dual statements for weak left adjoints (which will also be relevant later in the section). We begin with some comments on functoriality.

\begin{construction}  
  \label{t--a--functoriality}
  Let $\X : \AA \to \AA'$ be a functor of $(\infty,2)$-categories, let $G : \A \to \B$ be a map in $\AA$, and suppose given a weak right adjoint $H : \B \to \A$ to $G$ and a weak right adjoint $H' : \X(\B) \to \X(\A)$ to $\X(G)$. Then there is a canonical comparison map $\gamma : \X(H) \to H'$: namely, letting $\coun : GH \to \id_\B$ be the counit map, $\gamma$ is the map carried to $\X(\coun)$ under the composite equivalence
  \begin{align*}
    \Map_{\uMap_{\AA'}(\X(\B),\X(\A))}(\X(H),H')
    &\iso \Map_{\uMap_{\AA'}(\X(\B),\X(\B))}(\X(G)\X(H), \id_{\X(\B)}) \\
    &\iso \Map_{\uMap_{\AA'}(\X(\B),\X(\B))}(\X(GH), \X(\id_\B)).
  \end{align*}
\end{construction}

\begin{remark}
  \label{t--a--strict-functoriality}
  In the situation of \cref{t--a--functoriality}, if $H$ is in fact right adjoint to $G$, then the map $\X(H) \to H'$ is an isomorphism. This is because functors of $(\infty,2)$-categories preserve adjunctions.
\end{remark}

\begin{example}
  \label{t--a--assembly}
  Let us return to the situation of \cref{t--a--PrL}: $G : \A \to \B$ is a map in $\LMod_\C(\PrL)$, $H : \B \to \A$ is its right adjoint in $\Cat_\infty$, and $\til{H} : \B \to \A$ is its weak right adjoint in $\LMod_\C(\PrL)$. Applying \cref{t--a--functoriality} in the context of the forgetful functor $\LMod_\C(\PrL) \to \Cat_\infty$, we obtain a comparison map $\gamma : \til{H} \to H$. By \cref{t--a--strict-functoriality}, $\gamma$ is an isomorphism if $\til{H}$ is a right adjoint in $\LMod_\C(\PrL)$. In fact, the discussion in \cref{t--a--PrL} shows that the converse holds as well: both conditions are equivalent to $H$ being strictly $\C$-linear and preserving colimits. In general, the map $\gamma$ may be understood as an ``assembly map'': it exhibits $\til{H}$ as final among colimit preserving $\C$-linear functors from $\B$ to $\A$ equipped with a lax $\C$-linear natural transformation to $H$; this is a rephrasing of the universal property of the weak right adjoint.
\end{example}

\begin{remark}
  \label{t--a--functoriality-composition}
  Suppose given a commutative diagram of $(\infty,2)$-categories
  \[
    \begin{tikzcd}
      \AA \ar[rr, "\X"] \ar[dr, "\Z", swap] &
      &
      \AA' \ar[dl, "\Y"] \\
      &
      \AA''.
    \end{tikzcd}
  \]
  Let $G : \A \to \B$ be a map in $\AA$, and suppose given weak right adjoints $H : \B \to \A$ of $G$, $H' : \X(\B) \to \X(\A)$ of $\X(G)$, and $H'' : \Z(\B) \to \Z(\A)$ of $\Z(G)$. By \cref{t--a--functoriality}, we have canonical maps $\gamma_\X : \X(H) \to H'$, $\gamma_\Y : \Y(H') \to H''$, and $\gamma_\Z : \Z(H) \to H''$. These fit into a commutative diagram
  \[
    \begin{tikzcd}
      \Y\X(H) \ar[r, "\gamma_\X"] \ar[d, "\sim", swap] &
      \Y(H') \ar[d, "\gamma_\Y"] \\
      \Z(H) \ar[r, "\gamma_\Z"] &
      H''.
    \end{tikzcd}
  \]
\end{remark}

\begin{remark}
  \label{t--a--functoriality-transformation}
  Suppose given functors of $(\infty,2)$-categories $\o\X,\X : \AA \to \AA'$ and a natural transformation $\tau : \o\X \to \X$. Let $G : \A \to \B$ be a map in $\AA$, and suppose given weak right adjoints $H : \B \to \A$ of $G$ and $\o{H}' : \o\X(\B) \to \o\X(\A)$ of $\o\X(G)$, and a right adjoint $H' : \X(\B) \to \X(\A)$ of $\X(G)$. We have naturality squares
  \[
    \begin{tikzcd}
      \o\X(A) \ar[r, "\tau_A"] \ar[d, "\o\X(G)", swap] &
      \X(A) \ar[d, "\X(G)"] \\
      \o\X(B) \ar[r, "\tau_B"] &
      \X(B),
    \end{tikzcd}
    \quad\quad\quad
    \begin{tikzcd}
      \o\X(A) \ar[r, "\tau_A"] &
      \X(A) \\
      \o\X(B) \ar[r, "\tau_B"] \ar[u, "\o\X(H)"] &
      \X(B). \ar[u, "\X(H)", swap]
    \end{tikzcd}
  \]
  The first of these induces a map $\bc : \tau_A \circ \o{H}' \to H' \circ \tau_B$, namely the composition
  \[
    \tau_A \circ \o{H}' \lblto{\un} H' \circ \X(G) \circ \tau_A \circ \o{H}' \iso H' \circ \tau_B \circ \o\X(G) \circ \o{H}' \lblto{\coun} H' \circ \tau_B,
  \]
  where $\un$ is the unit for the adjunction between $\X(G)$ and $H'$ and $\coun$ is the counit for the weak adjunction between $\o\X(G)$ and $\o{H}'$. And by \cref{t--a--functoriality}, we have canonical maps $\o\gamma : \o\X(H) \to \o{H}'$ and $\gamma : \X(H) \to H'$. These fit into a commutative diagram
  \[
    \begin{tikzcd}
      \tau_A \circ \o\X(H) \ar[r, "\sim"] \ar[d, "\o\gamma", swap] &
      \X(H) \circ \tau_B \ar[d, "\gamma"] \\
      \tau_A \circ \o{H}' \ar[r, "\bc"] &
      H' \circ \tau_B,
    \end{tikzcd}
  \]
  where the upper equivalence comes from the second naturality square above.
\end{remark}

We next address the interaction of weak adjoints with composition.

\begin{construction}
  \label{t--c--categorical-comp-right}
  Let $\AA$ be an $(\infty,2)$-category. Suppose given maps $G' : \A \to \B$ and $G'' : \B \to \C$ in $\AA$, as well as maps $H' : \B \to \A$ and $H'' : \C \to \B$ and maps $\coun' : G'H' \to \id_\B$ and $\coun'' : G''H'' \to \id_\C$. Then we may form the composition
  \[
    \smash{\coun''' : GH'H'' \iso G''G'H'H'' \lblto{\coun'} G''H'' \lblto{\coun''} \id_\C.}
  \]
  Now set $G := G'' \circ G' : \A \to \C$, and suppose that $G$ admits a weak right adjoint $H : \C \to \A$ in $\AA$; let $\coun : GH \to \id_\C$ be the associated counit map. Then we have a unique map $\comp : H'H'' \to H$ making the following diagram commute:
  \[
    \begin{tikzcd}
      GH'H'' \ar[r, "\comp"] \ar[dr, "\coun'''", swap] &
      GH \ar[d, "\coun"] \\
      &
      \id_\C.
    \end{tikzcd}
  \]
\end{construction}

\begin{proposition}
  \label{t--c--categorical-equiv-right}
  In the situation of \cref{t--c--categorical-comp-right}, suppose that $\coun'$ and $\coun''$ exhibit $H'$ and $H''$ as weak right adjoints of $G'$ and $G''$, and that one of the following further conditions is satisfied:
  \begin{enumerate}
  \item \label{t--c--categorical-equiv-right--post}
    $\coun'$ exhibits $H'$ as right adjoint to $G'$;
  \item \label{t--c--categorical-equiv-right--pre}
    $H''$ admits a right adjoint $I'' : \B \to \C$.
  \end{enumerate}
  Then the map $\comp : H'H'' \to H$ is an isomorphism.
\end{proposition}

\begin{proof}
  Let $J : \C \to \A$ be any map in $\AA$. Suppose first that \cref{t--c--categorical-equiv-right--post} is satisfied. Then we have equivalences
  \begin{align*}
    \Map_{\uMap_{\AA}(\C,\A)}(J, H'H'')
    \iso{} & \Map_{\uMap_{\AA}(\C,\B)}(G' J, H'') \\
    \iso{} & \Map_{\uMap_{\AA}(\C,\C)}(G''G' J, \id_\C) \\
    \iso{} & \Map_{\uMap_{\AA}(\C,\C)}(G J, \id_\C) \\
    \iso{} & \Map_{\uMap_{\AA}(\C,\A)}(J, H)
  \end{align*}
  (the first using that $H'$ is right adjoint to $G'$, the second using that $H''$ is weakly right adjoint to $G''$, the third using the definition of $G$, and the fourth using that $H$ is weakly right adjoint to $G$). It is straightforward to verify that the composite equivalence carries, in the case $J = H'H''$, the identity map $\id_{H'H''}$ to the map $\comp$, so we deduce that $\comp$ is an equivalence.

  Now suppose that \cref{t--c--categorical-equiv-right--pre} is satisfied. Then we have equivalences
  \begin{align*}
    \Map_{\uMap_{\AA}(\C,\A)}(J, H'H'')
    \iso{} & \Map_{\uMap_{\AA}(\B,\A)}(JI'', H') \\
    \iso{} & \Map_{\uMap_{\AA}(\B,\B)}(G'JI'', \id_\B) \\
    \iso{} & \Map_{\uMap_{\AA}(\C,\B)}(G'J, H'') \\
    \iso{} & \Map_{\uMap_{\AA}(\C,\C)}(G''G' J, \id_\C) \\
    \iso{} & \Map_{\uMap_{\AA}(\C,\C)}(G J, \id_\C) \\
    \iso{} & \Map_{\uMap_{\AA}(\C,\A)}(J, H)
  \end{align*}
  (the first and third using that $I''$ is right adjoint to $H''$, the second using that $H'$ is weakly right adjoint to $G'$, and the others being the same as in the first case above). Again, the composite equivalence carries $\id_{H'H''}$ to $\comp$, proving that $\comp$ is an equivalence.
\end{proof}

Finally, we close this subsection by recalling the following result from \cite{heyer-mann--six}. It is not stated directly in terms of weak adjoints, but there is a connection via \itemref{t--a--weak-adjoint-props}{presentable}, and we will see later in the section that it gives a useful criterion for determining when a weak adjoint is actually an adjoint.

\begin{proposition}[{\citebare[Proposition D.2.8]{heyer-mann--six}}]
  \label{t--p--right-adjoint-criterion}
  Let $\AA$ be an $(\infty,2)$-category, and let $G : \A \to \B$ be a map in $\AA$. Then the following are equivalent:
  \begin{enumerate}
  \item \label{t--p--right-adjoint-criterion--adjoint}
    $G$ admits a right adjoint $H : \B \to \A$;
  \item \label{t--p--right-adjoint-criterion--criterion}
    the functors $g_\B : \uMap_{\AA}(\B,\A) \to \uMap_{\AA}(\B,\B)$ and $g_\A : \uMap_{\AA}(\A,\A) \to \uMap_{\AA}(\A,\B)$ given by composition with $G$ admit right adjoint functors $h_\B : \uMap_{\AA}(\B,\B) \to \uMap_{\AA}(\B,\A)$ and $h_\A : \uMap_{\AA}(\A,\B) \to \uMap_{\AA}(\A,\A)$, and the canonical map
    \[
      h_\B(\id_\B) \circ G \to h_\A(G)
    \]
    in $\uMap_{\AA}(\A,\A)$ (adjoint to the map $G \circ h_\B(\id_\B) \circ G \to G$ in $\uMap_{\AA}(\A,\B)$ induced by the counit map $G \circ h_\B(\id_\B) \to \id_\B$) becomes an isomorphism upon applying $\Map_{\uMap_{\AA}(\A,\A)}(\id_\A,-)$.
  \end{enumerate}
\end{proposition}

\begin{proof}
  If \cref{t--p--right-adjoint-criterion--adjoint} holds, then the functors $g_\B$ and $g_\A$ admit right adjoints given by composition with $H$, and hence \cref{t--p--right-adjoint-criterion--criterion} holds. Conversely, assume that \cref{t--p--right-adjoint-criterion--criterion} holds. We claim that then $H := h_\B(\id_\B) : \B \to \A$ is right adjoint to $G$. By hypothesis, we have an equivalence
  \[
    \Map_{\uMap_{\AA}(\A,\A)}(\id_\A,HG) \iso     \Map_{\uMap_{\AA}(A,A)}(\id_\A,h_\A(G));
  \]
  let $\un : \id_\A \to HG$ be the map corresponding under this equivalence to the unit map for the adjunction $g_\A \dashv h_\A$. Let $\coun : GH \to \id_\B$ be the map given by the counit for the adjunction $g_\B \dashv h_\B$. One may verify that $\un$ and $\coun$ satisfy the triangle identities to complete the argument.
\end{proof}


\subsection{Norm maps and Tate cohomology}
\label{t--d}

In this subsection, we define norm maps and Tate cohomology in the setting of our three functor formalism $\D$. We will do so by considering weak adjoints (\cref{t--a--weak-adjoint}) in the $(\infty,2)$-categories of kernels $\K_\D(S)$ (\itemref{t--f}{K}).

\begin{notation}
  \label{t--d--product}
  Let $f : T \to S$ be a $\D$-$!$-able map in $\S$. Then $f$ admits a base change along itself, so we have a commutative diagram
  \[
    \begin{tikzcd}
      T \ar[dr, "\diag{f}"] \ar[drr, "\id", bend left=20] \ar[ddr, "\id", bend right=20, swap] &
      &
      \\
      &
      T \times_S T \ar[r, "\rp{f}"] \ar[d, "\lp{f}", swap] &
      T \ar[d, "f"] \\
      &
      T \ar[r, "f"] &
      S;
    \end{tikzcd}
  \]
  here, and below, $\lp{f}$ denotes projection onto the left factor, $\rp{f}$ projection onto the right factor, and $\diag{f}$ the diagonal map. Assuming that $\lp{f}$ is $\D$-$*$-able, we define
  \[
    \primtwist_f := (\lp{f})_*(\diag{f})_!(\unit_T) \in \D(T),
  \]
  and assuming that $\lp{f}$ is $\D$-$\sharp$-able, we define
  \[
    \suavetwist_f := (\lp{f})_\sharp(\diag{f})_!(\unit_T) \in \D(T);
  \]
  when there is danger of confusion, we will write $\primtwist_f^\D$ and $\suavetwist_f^\D$ rather than $\primtwist_f$ and $\suavetwist_f$.
\end{notation}

\begin{proposition}
  \label{t--d--kernel-weak-adjoints}
  Let $f : T \to S$ be a $\D$-$!$-able map in $\S$. Then:
  \begin{enumerate}
  \item \label{t--d--kernel-weak-adjoints--prim}
    assuming that $\lp{f}$ is $\D$-$*$-able, the map $[\unit_T] : [S] \to [T]$ in $\K_\D(S)$ has weak right adjoint given by $[\primtwist_f] : [T] \to [S]$;
  \item \label{t--d--kernel-weak-adjoints--suave}
    assuming that $\lp{f}$ is $\D$-$\sharp$-able, the map $[\unit_T] : [S] \to [T]$ in $\K_\D(S)$ has weak left adjoint given by $[\suavetwist_f] : [T] \to [S]$.
  \end{enumerate}
\end{proposition}

\begin{proof}
  We have $\id_{[T]} \iso [(\diag{f})_!(\unit_T)]$, and, for $L \in \D(T)$ with associated map $[L] : [T] \to [S]$ in $\K_\D(S)$, we have a natural identification
  \[
    [\unit_T] \circ [L] \iso [\lp{f}^*(L)].
  \]
  The two claims (which are formally dual) now follow from \itemref{t--a--weak-adjoint-props}{presentable}.
\end{proof}

\begin{remark}
  \label{t--d--Dprime}
  Suppose given another three functor formalism $\D' : \Span(\S,\S_!) \to \Cat_\infty$ and a map $\alpha : \D \to \D'$. Let $f : T \to S$ be a $\D$-$!$-able (hence also $\D'$-$!$-able) map in $\S$. Assuming that $\lp{f}$ is $\D$-$*$-able and $\D'$-$*$-able, the functor $\phi_\alpha : \K_\D(S) \to \K_{\D'}(S)$ of \itemitemref{t--f}{K}{Dprime} induces, by \cref{t--d--kernel-weak-adjoints,t--a--functoriality}, a canonical comparison map $\smash{\alpha(\primtwist_f^\D) \to \primtwist_f^{\D'}}$ in $\D'(T)$. Dually, if $\lp{f}$ is $\D$-$\sharp$-able and $\D'$-$\sharp$-able, we have a comparison map $\smash{\suavetwist_f^{\D'} \to \alpha(\suavetwist_f^\D)}$ in $\D'(T)$.
\end{remark}

The following terminology is motivated by the constructions which follow.

\begin{definition}
  \label{t--d--normed}
  Let $f : T \to S$ be a $\D$-$!$-able map in $\S$. We say that $f$ is:
  \begin{itemize}
  \item \emph{$\D$-$*$-normed} if both $f$ and $\lp{f}$ are $\D$-$*$-able;
  \item \emph{$\D$-$\sharp$-normed} if both $f$ and $\lp{f}$ are $\D$-$\sharp$-able;
  \item \emph{$\D$-normed} if it is $\D$-$\sharp$-normed and $\D$-$*$-normed;
  \item \emph{$\D$-Tate} if it is $\D$-normed and $\suavetwist_f$ is an invertible object of $\D(T)$.
  \end{itemize}
\end{definition}

\begin{construction}
  \label{t--d--star-norm}
  Let $f : T \to S$ be a $\D$-$*$-normed map in $\S$. The functor $\Psi_{\D,S} : \K_\D(S) \to \Cat_\infty$ (\itemitemref{t--f}{K}{Psi}) carries the map $[\unit_T] : [S] \to [T]$ to the functor $f^* : \D(S) \to \D(T)$, and it carries the map $[\primtwist_f] : [T] \to [S]$ to the functor $f_!(- \otimes \primtwist_f) : \D(T) \to \D(S)$. Since $f^*$ admits a right adjoint $f_*$, it follows from \itemref{t--d--kernel-weak-adjoints}{prim} and \cref{t--a--functoriality} that we have a canonical natural transformation
  \[
    \Nm^*_f : f_!(- \otimes \primtwist_f) \to f_*
  \]
  of functors from $\D(T)$ to $\D(S)$; we refer to this as the \emph{$*$-norm map} associated to $f$.
\end{construction}

\begin{construction}
  \label{t--d--sharp-norm}
  Let $f : T \to S$ be a $\D$-$\sharp$-normed map in $\S$. The functor $\Psi_{\D,S} : \K_\D(S) \to \Cat_\infty$ carries the map $[\unit_T] : [S] \to [T]$ to the functor $f^* : \D(S) \to \D(T)$ and the map $[\suavetwist_f] : [T] \to [S]$ to the functor $f_!(- \otimes \suavetwist_f) : \D(T) \to \D(S)$. Since $f^*$ admits a left adjoint $f_\sharp$, it follows from \itemref{t--d--kernel-weak-adjoints}{suave} and \cref{t--a--functoriality} that we have a canonical natural transformation
  \[
    \Nm^\sharp_f : f_\sharp \to f_!(- \otimes \suavetwist_f)
  \]
  of functors from $\D(T)$ to $\D(S)$; we refer to this as the \emph{$\sharp$-norm map} associated to $f$.
\end{construction}

\begin{construction}
  \label{t--d--sharpstar-norm}
  let $f : T \to S$ be a $\D$-normed map in $\S$. Composing the $\sharp$-norm map of \cref{t--d--sharp-norm} with the $*$-norm map of \cref{t--d--star-norm}, we obtain a natural transformation
  \[
    \Nm_f : f_\sharp(- \otimes \primtwist_f) \to f_*(- \otimes \suavetwist_f)
  \]
  of functors from $\D(T)$ to $\D(S)$, which we refer to as the \emph{norm map} associated to $f$.
\end{construction}

\begin{construction}
  \label{t--d--norm}
  Let $f : T \to S$ be a $\D$-Tate map in $\S$. Then we define
  \[
    \poincaretwist_f := \suavetwist_f^{-1} \otimes \primtwist_f \in \D(T),
  \]
  so that we may rewrite the norm map of \cref{t--d--sharpstar-norm} equivalently as $\Nm_f : f_\sharp(- \otimes \poincaretwist_f) \to f_*$. Assuming that $\D(S)$ is pointed and admits cofibers, we then define the functor $f_\t : \D(T) \to \D(S)$ and the \emph{canonical map} $\can_f : f_* \to f_\t$ by forming a cofiber sequence
  \[
    f_\sharp(- \otimes \poincaretwist_f) \lblto{\Nm_f} f_* \lblto{\can_f} f_\t;
  \]
  we refer to $f_\t$ as \emph{Tate pushforward} or \emph{Tate cohomology} along $f$.
\end{construction}

\begin{remark}
  \label{t--d--norm-linearity}
  In the situation of \cref{t--d--star-norm}, the counit map $\coun : [\unit_T][\primtwist_f] \to \id_{[T]}$ induces a $\D(S)$-linear natural transformation $f^*f_!(- \otimes \primtwist_f) \to \id_{\D(T)}$ (\itemitemref{t--f}{K}{Psi}), from which it follows that the $*$-norm map $\Nm^*_f : f_!(- \otimes \primtwist_f) \to f_*$ is canonically lax $\D(S)$-linear. Dually, in the situation of \cref{t--d--sharp-norm}, the unit map $\un : \id_{[T]} \to [\unit_T][\suavetwist_f]$ induces a $\D(S)$-linear natural transformation $\id_{\D(T)} \to f^*f_!(- \otimes \suavetwist_f)$, from which it follows that the $\sharp$-norm map $\Nm^\sharp_f : f_\sharp \to f_!(- \otimes \suavetwist_f)$ is canonically oplax $\D(S)$-linear.
\end{remark}

To end this subsection, we note that, in some cases, the above norm maps may be thought of in terms of the notion of ``assembly map'' discussed in \cref{t--a--assembly}.

\begin{remark}
  \label{t--d--norm-assembly}
  Suppose that $\D$ is presentable (\itemitemref{t--f}{D}{presentable}). Then, for $S \in \S$, the functor $\Psi_{\D,S}$ admits a factorization
  \[
    \begin{tikzcd}[row sep=3ex]
      &
      \Mod_{\D(S)}(\PrL) \ar[d] \\
      \K_\D(S) \ar[ur, "\Psi'_{\D,S}"] \ar[r, "\Psi_{\D,S}"] &
      \Cat_\infty,
    \end{tikzcd}
  \]
  where the downward arrow is the forgetful functor (cf. \itemitemref{t--f}{K}{Psi}). Letting $f : T \to S$ be a $\D$-$!$-able map in $\S$, the functor $\Psi'_{\D,S}$ induces functors
  \begin{align*}
    \D(T) \iso \uMap_{\K_\D(S)}([S],[T]) \to \uMap_{\Mod_{\D(S)}(\PrL)}(\D(S),\D(T)), \quad L \mapsto f^*(-) \otimes L \\
    \D(T) \iso \uMap_{\K_\D(S)}([T],[S]) \to \uMap_{\Mod_{\D(S)}(\PrL)}(\D(T),\D(S)), \quad L \mapsto f_!(- \otimes L).
  \end{align*}
  The first of these is always an equivalence. Sometimes the second is as well; a sufficient condition for this to hold is that the canonical maps
  \[
    \D(T) \otimes_{\D(S)} \D(T) \to \D(T \times_S T), \quad \D(T) \otimes_{\D(S)} \D(T) \otimes_{\D(S)} \D(T) \to \D(T \times_S T \times_S T)
  \]
  in $\CAlg(\PrL)$, induced by pullback along the projection maps, be equivalences. Indeed, then $\D$ preserves the self-duality of $T$ in $\Span((\S_!)_{/S})$ (see \itemitemref{t--f}{span}{frobenius}): that is, the composition
  \[
    \D(T) \otimes_{\D(S)} \D(T) \lblto{\otimes} \D(T) \lblto{f_!} \D(S)
  \]
  is a duality datum in $\Mod_{\D(S)}(\PrL)$.

  Suppose that the second functor is in fact an equivalence. Then it is immediate from the definition of weak adjoints that the functor $\Psi'_{\D,S}$ must carry a weak right (resp. left) adjoint of the map $[\unit_T] : [S] \to [T]$ in $\K_\D(S)$ to a weak right (resp. left) adjoint of the map $\Psi'_{\D,S}([\unit_T]) \iso f^* : \D(S) \to \D(T)$ in $\Mod_{\D(S)}(\PrL)$. The presentability hypothesis guarantees that the weak right adjoints exist, and it then follows from \cref{t--a--assembly} that the $*$-norm map $\Nm_f^* : f_!(- \otimes \primtwist_f) \to f_*$ is a $\D(S)$-linear assembly map: that is, it exhibits the source as the final colimit preserving $\D(S)$-linear functor equipped with a lax $\D(S)$-linear map to $f_*$. Dually, assuming that $f$ is $\sharp$-normed, the $\sharp$-norm map $\Nm_f^\sharp : f_\sharp \to f_!(- \otimes \suavetwist_f)$ is a $\D(S)$-linear coassembly map: that is, it exhibits the target as the initial colimit preserving $\D(S)$-linear functor equipped with an oplax $\D(S)$-linear map from $f_\sharp$.
\end{remark}


\subsection{Poincar\'e duality}
\label{t--p}

In the last subsection, we defined the class of $\D$-normed maps $f : T \to S$ in $\S$, and associated to each such $f$ a norm map comparing homology and cohomology along it. As was indicated in \cref{i--bg}, one matter of interest is whether or not this norm map is an isomorphism: when it is, we may interpret this as a form of Poincar\'e duality. In \cite{heyer-mann--six}, it is shown that this form of Poincar\'e duality is closely related to the weak adjoints of \cref{t--d--kernel-weak-adjoints} actually being adjoints, as we now review.

\begin{proposition}[{\citebare[Lemma 4.5.5]{heyer-mann--six}}]
  \label{t--p--prim-equiv}
  Let $f : T \to S$ be a $\D$-$!$-able map in $\S$. Then the following conditions are equivalent:
  \begin{enumerate}
  \item \label{t--p--prim-equiv--adjoint}
    $[\unit_T] : [S] \to [T]$ admits a right adjoint in $\K_\D(S)$;
  \item \label{t--p--prim-equiv--iso}
    $f$ is $\D$-$*$-normed and the natural transformation $\Nm^*_f : f_!(- \otimes \primtwist_f) \to f_*$ is an isomorphism;
  \item \label{t--p--prim-equiv--unit}
    $f$ is $\D$-$*$-normed and the map
    \[
      \Map_{\D(S)}(\unit_S,f_!(\primtwist_f)) \to \Map_{\D(S)}(\unit_S,f_*(\unit_T)) \iso \Map_{\D(T)}(\unit_T,\unit_T)
    \]
    induced by $\Nm^*_f$ is an equivalence.
  \end{enumerate}
\end{proposition}

\begin{proof}
  The functors
  \[
    \uMap_{\K_\D(S)}([T],[S]) \to \uMap_{\K_\D(S)}([T],[T]), \quad
    \uMap_{\K_\D(S)}([S],[S]) \to \uMap_{\K_\D(S)}([S],[T])
  \]
  given by composition with $[\unit_T]$ identify with the functors
  \[
    \lp{f}^* : \D(T) \to \D(T \times_S T), \quad
    f^* : \D(S) \to \D(T).
  \]
  Thus, if \cref{t--p--prim-equiv--adjoint} holds, then $f$ is $\D$-$*$-normed by the forward implication in \cref{t--p--right-adjoint-criterion}; and in this case $\Nm_f^*$ is an isomorphism by \cref{t--a--strict-functoriality}, proving that \cref{t--p--prim-equiv--iso} holds. That \cref{t--p--prim-equiv--iso} implies \cref{t--p--prim-equiv--unit} is obvious. Finally, \cref{t--p--prim-equiv--unit} implies \cref{t--p--prim-equiv--adjoint} by the reverse implication of \cref{t--p--right-adjoint-criterion}.
\end{proof}

\begin{proposition}[{\citebare[Lemma 4.5.6]{heyer-mann--six}}]
  \label{t--p--suave-equiv}
  Let $f : T \to S$ be a $\D$-$!$-able map in $\S$. Then the following conditions are equivalent:
  \begin{enumerate}
  \item \label{t--p--suave-equiv--adjoint}
    $[\unit_T] : [S] \to [T]$ admits a left adjoint in $\K_\D(S)$;
  \item \label{t--p--suave-equiv--iso}
    $f$ is $\D$-$\sharp$-normed and the natural transformation $\Nm^\sharp_f : f_\sharp \to f_!(- \otimes \suavetwist_f)$ is an isomorphism;
  \item \label{t--p--suave-equiv--unit}
    $f$ is $\D$-$\sharp$-normed and the map
    \[
      \Map_{\D(T)}(\unit_T,\unit_T) \iso \Map_{\D(S)}(f_\sharp(\unit_T),\unit_S) \to \Map_{\D(S)}(f_!(\suavetwist_f),\unit_S)
    \]
    induced by $\Nm^\sharp_f$ is an equivalence.
  \end{enumerate}
\end{proposition}

\begin{proof}
  This is formally dual to \cref{t--p--prim-equiv}.
\end{proof}

To facilitate later discussion involving these duality phenomena, we introduce the following terminology (three fifths of which comes from \cite{heyer-mann--six}).

\begin{definition}
  \label{t--p--terms}
  Let $f : T \to S$ be a $\D$-$!$-able map in $\S$. Then we say that $f$ is:
  \begin{itemize}
  \item \emph{$\D$-prim} if the equivalent conditions of \cref{t--p--prim-equiv} hold;
  \item \emph{$\D$-terse} if it is $\D$-prim and $\primtwist_f$ is an invertible object of $\D(T)$;
  \item \emph{$\D$-suave} if the equivalent conditions of \cref{t--p--suave-equiv} hold;
  \item \emph{$\D$-smooth} if it is $\D$-suave and $\suavetwist_f$ is an invertible object of $\D(T)$;
  \item \emph{$\D$-Poincar\'e} if it is $\D$-terse and $\D$-smooth.
  \end{itemize}
\end{definition}

\begin{remark}
  \label{t--p--poincare-duality}
  Let $f : T \to S$ be a $\D$-Poincar\'e map in $\S$. Then $f$ is in particular $\D$-Tate, and moreover the norm map $\Nm_f : f_\sharp(- \otimes \poincaretwist_f) \to f_*$ is an isomorphism and $\poincaretwist_f$ is an invertible object of $\D(T)$. We think of the combination of these assertions as a statement of Poincar\'e duality for homology and cohomology along $f$.
\end{remark}

As indicated by \cref{t--p--poincare-duality}, our interest here will primarily be in the second conditions in each of \cref{t--p--prim-equiv,t--p--suave-equiv}, i.e. the maps $\Nm_f^*$ and $\Nm_f^\sharp$ being isomorphisms. The fact that these conditions are equivalent to the a priori stronger first characterizations there, in terms of adjunctions in $\K_\D(S)$, implies that they are more robust than one might expect\footnote{This pleasant surprise was explained to me by Hesselholt.}; see for instance \cref{t--b--closure,t--b--descent,t--c--closure}.

Before moving on, let us record a few basic observations about these conditions.

\begin{remark}
  \label{t--p--pr-bc}
  Let $f : T \to S$ be a map in $\S$. If $f$ is $\D$-prim, then $f$ satisfies the $\D$-$*$-projection formula and $\D$-$*$--base change (see \itemref{t--f}{D}). Dually, if $f$ is $\D$-suave, then $f$ satisfies the $\D$-$\sharp$-projection formula and $\D$-$\sharp$--base change.

  Let us justify in the prim case. For the projection formula, note that the $*$-norm map $\Nm^*_f : f_!(- \otimes \primtwist_f) \to f_*$ is a map, in this case an isomorphism, of lax $\D(S)$-linear functors (\cref{t--d--norm-linearity}), and that its source is strictly $\D(S)$-linear (see \itemitemref{t--f}{K}{Psi}). For base change, see \cite[Lemma 4.5.13]{heyer-mann--six}, and note that primness is stable under base change (\cref{t--b--closure}). 
\end{remark}

\begin{remark}
  \label{t--p--assembly}
  As in \cref{t--d--norm-assembly}, suppose that $\D$ is presentable, and let $f : T \to S$ be a $\D$-$!$-able map such that the functor
  \[
    \D(T) \to \uMap_{\Mod_{\D(S)}(\PrL)}(\D(T),\D(S)), \quad L \mapsto f_!(- \otimes L)
  \]
  is an equivalence. It follows from the discussion in loc. cit. that $f$ is $\D$-prim if and only if the functor $f_* : \D(T) \to \D(S)$ is $\D(S)$-linear (i.e. $f$ satisfies the $\D$-$*$-projection formula) and preserves colimits, and that $f$ is $\D$-suave if and only if it is $\sharp$-normed and satisfies the $\D$-$\sharp$-projection formula.
\end{remark}

\begin{remark}
  \label{t--p--Dprime}
  Suppose given another three functor formalism $\D' : \Span(\S,\S_!) \to \Cat_\infty$ and a map $\alpha : \D \to \D'$. Let $f : T \to S$ be a $\D$-$!$-able (hence also $\D'$-$!$-able) map in $\S$. If $f$ is $\D$-prim, then, by \cref{t--a--strict-functoriality}, $f$ is also $\D'$-prim and the comparison map $\smash{\alpha(\primtwist_f^\D) \to \primtwist_f^{\D'}}$ of \cref{t--d--Dprime} is an isomorphism. Dually, if $f$ is $\D$-suave, then $f$ is $\D'$-suave and we have an isomorphism $\smash{\suavetwist_f^{\D'} \isoto \alpha(\suavetwist_f^\D)}$.
\end{remark}

\begin{remark}
  \label{t--p--op}
  Let $f : T \to S$ be a $\D$-$!$-able map in $\S$. The canonical equivalence $\K_\D(S) \iso \K_\D(S)^\op$ (\itemitemref{t--f}{K}{op}) carries a right (resp. left) adjoint of the map $[\unit_T] : [S] \to [T]$ in $\K_\D(S)$ to a left (resp. right) adjoint of the map $[\unit_T] : [T] \to [S]$ in $\K_\D(S)$. This gives additional characterizations and consequences of primness and suaveness. For instance, if $f$ is $\D$-prim, then $[\primtwist_f] : [S] \to [T]$ is left adjoint to $[\unit_T] : [T] \to [S]$, and by applying the functor $\Psi_{\D,S} : \K_\D(S) \to \Cat_\infty$ we deduce that the functor $f^*(-) \otimes \primtwist_f : \D(S) \to \D(T)$ is left adjoint to the functor $f_! : \D(T) \to \D(S)$. Dually, if $f$ is $\D$-suave, then $[\suavetwist_f] : [S] \to [T]$ is right adjoint to $[\unit_T] : [T] \to [S]$, implying that $f^*(-) \otimes \suavetwist_f : \D(S) \to \D(T)$ is right adjoint to $f_! : \D(T) \to \D(S)$.
\end{remark}

\begin{remark}[{\citebare[Remark 4.5.12]{heyer-mann--six}}]
  \label{t--p--dualizable}
  Let $f : T \to S$ be a map in $\S$. If $f$ is $\D$-prim and $\primtwist_f$ is a dualizable object of $\D(T)$, then $\primtwist_f$ is in fact invertible, i.e. $f$ is $\D$-terse. Dually, if $f$ is $\D$-suave and $\suavetwist_f$ is dualizable, then $f$ is $\D$-smooth.

  For the prim case, recall from \cref{t--p--op} that $f_! : \D(T) \to \D(S)$ has left adjoint $f^*(-) \otimes \primtwist_f : \D(S) \to \D(T)$. On the other hand, our dualizability assumption gives us that $- \otimes \primtwist_f : \D(T) \to \D(T)$ has left adjoint $- \otimes \dual{\primtwist}_f : \D(T) \to \D(T)$. Thus, the $*$-norm map $\Nm^*_f : f_!(- \otimes \primtwist_f) \to f_*$, which in this case is an isomorphism, induces upon passing to left adjoints a natural isomorphism $f^* \iso f^*(-) \otimes \primtwist_f \otimes \dual{\primtwist}_f$. We deduce that $\primtwist_f \otimes \dual{\primtwist}_f \iso \unit_T$, as desired.
\end{remark}

To close this subsection, we recall that there are special classes of $\D$-prim maps and $\D$-suave maps for which the respective twisting object $\primtwist_f$ or $\suavetwist_f$ is canonically trivialized (cf. \cite[\textsection 4.6]{heyer-mann--six}).

\begin{definition}
  \label{t--p--truncated}
  Let $f : T \to S$ be a map in $\S$. The property that $f$ be \emph{$n$-truncated} is defined inductively for integers $n \ge -2$ as follows:
  \begin{itemize}
  \item we say that $f$ is $(-2)$-truncated if it is an isomorphism;
  \item for $n \ge -1$, we say that $f$ is $n$-truncated if it admits a base change along itself and the diagonal map $\diag{f} : T \to T \times_S T$ is $(n-1)$-truncated.
  \end{itemize}
  We say that $f$ is \emph{truncated} if it is $n$-truncated for some integer $n \ge -2$.
\end{definition}

\begin{definition}
  \label{t--p--etale}
  Let $f : T \to S$ be a truncated, $\D$-$!$-able map in $\S$. The property that $f$ be \emph{$\D$-\'etale} is defined by the following stipulations: if $f$ is an isomorphism, then $f$ is $\D$-\'etale; and $f$ is $\D$-\'etale if and only if it is $\D$-suave and $\diag{f}$ is $\D$-\'etale.
\end{definition}

\begin{proposition}
  \label{t--p--etale-twist}
  Let $f : T \to S$ be an $\D$-\'etale map in $\S$. Then there is a canonical isomorphism $\suavetwist_f \iso \unit_T$ in $\D(T)$; in particular, $f$ is $\D$-smooth.
\end{proposition}

\begin{proof}
  Choose $n \in \ZZ_{\ge-2}$ minimally such that $f$ is $n$-truncated. We proceed by induction on $n$. Suppose first that $n=-2$, meaning that $f$ is an isomorphism. Then $\lp{f}$ and $\diag{f}$ are inverse isomorphisms. It follows that we have an identification $(\diag{f})_! \iso \lp{f}^*$ and that $(\lp{f})_\sharp$ is inverse to $\lp{f}^*$. The composition
  \[
    \suavetwist_f = (\lp{f})_\sharp(\diag{f})_!(\unit_T) \iso (\lp{f})_\sharp(\lp{f})^*(\unit_T) \iso \unit_T
  \]
  gives the desired isomorphism.
  
  Suppose now that $n \ge -1$. Then $\diag{f}$ is an $(n-1)$-truncated $\D$-\'etale map, so the inductive hypothesis grants us an isomorphism $\suavetwist_{\diag{f}} \iso \unit_T$ in $\D(T)$. Since \'etaleness by definition implies suaveness, it follows that we have a natural isomorphism $\Nm^\sharp_{\diag{f}} : (\diag{f})_\sharp \isoto (\diag{f})_!$. The composition
  \[
    \suavetwist_f = (\lp{f})_\sharp(\diag{f})_!(\unit_T) \iso (\lp{f})_\sharp(\diag{f})_\sharp(\unit_T) \iso (\lp{f}\diag{f})_\sharp(\unit_T) \iso (\id_T)_\sharp(\unit_T) \iso \unit_T
  \]
  gives the desired isomorphism.
\end{proof}

\begin{definition}
  \label{t--p--proper}
  Let $f : T \to S$ be a truncated, $\D$-$!$-able map in $\S$. The property that $f$ be \emph{$\D$-proper} is defined by the following stipulations: if $f$ is an isomorphism, then $f$ is $\D$-proper; and $f$ is $\D$-proper if and only if it is $\D$-prim and $\diag{f}$ is $\D$-proper.
\end{definition}

\begin{proposition}
  \label{t--p--proper-twist}
  Let $f : T \to S$ be a $\D$-proper map in $\S$. Then there is a canonical isomorphism $\primtwist_f \iso \unit_T$ in $\D(T)$; in particular, $f$ is $\D$-terse.
\end{proposition}

\begin{proof}
  This is formally dual to \cref{t--p--etale-twist}.
\end{proof}


\subsection{Poincar\'e duality for H-spaces}
\label{t--h}

As alluded to at the end of \cref{i--bg}, it is a result of Browder \cite{browder--torsion} that any compact, connected H-space is Poincar\'e. In this subsection, we will prove a generalization of this result in our more abstract context (\cref{t--h--main}), following Klein's argument in \cite{klein--dualizing}.

\begin{definition}
  \label{t--h--H}
  Let $S \in \S$. By an \emph{$\AA_2$-monoid} over $S$, we mean a map $f : T \to S$ in $\S$ that admits a base change along itself, together with maps $\un_f : S \to T$ and $\comp_f : T \times_S T \to T$ and commutative diagrams
  \[
    \begin{tikzcd}[row sep=small]
      &
      T \ar[dr, "f"] \\
      S \ar[ur, "\un_f"] \ar[rr, "\id_S"] &
      &
      S,
    \end{tikzcd}
    \quad\quad\quad
    \begin{tikzcd}[row sep=small]
      &
      T \times_S T \ar[dr, "\comp_f"] \\
      T \ar[ur, "{(\id_T,\un_f \circ f)}"] \ar[rr, "\id_T"] &
      &
      T
    \end{tikzcd}
  \]
  (i.e. the structure of a right unital multiplication on $T$ over $S$). We say that an $\AA_2$-monoid $f : T \to S$ is \emph{left divisible} if the map $(\lp{f},\comp_f) : T \times_S T \to T \times_S T$ is an isomorphism.
\end{definition}

\begin{remark}
  \label{t--h--group}
  For $S \in \S$, a group object in $\S_{/S}$ is in particular a left divisible $\AA_2$-monoid over $S$.
\end{remark}

\begin{proposition}
  \label{t--h--constant}
  Let $S \in \S$, let $f : T \to S$ be a left divisible $\AA_2$-monoid over $S$, and let $\sigma$ denote the cartesian square
  \[
    \begin{tikzcd}
      T \times_S T \ar[r, "\rp{f}"] \ar[d, "\lp{f}", swap] &
      T \ar[d, "f"] \\
      T \ar[r, "f"] &
      S.
    \end{tikzcd}
  \]
  Then:
  \begin{enumerate}
  \item \label{t--h--constant--star}
    if $f$ is $\D$-$*$-normed and the base change map $\bc^*_\sigma : f^*f_* \to (\lp{f})_*(\rp{f})^*$ is an isomorphism, there is a canonical isomorphism $\primtwist_f \iso f^*(\o\primtwist_f)$ for $\o\primtwist_f := f_*(\eta_f)_!(\unit_S) \in \D(S)$.
  \item \label{t--h--constant--sharp}
    if $f$ is $\D$-$\sharp$-normed and the base change map $\bc^\sharp_\sigma : (\lp{f})_\sharp(\rp{f})^* \to f^*f_\sharp$ is an isomorphism, there is a canonical isomorphism $\suavetwist_f \iso f^*(\o\suavetwist_f)$ for $\o\suavetwist_f := f_\sharp(\eta_f)_!(\unit_S) \in \D(S)$.
  \end{enumerate}
\end{proposition}

\begin{proof}
  The two claims are formally dual; let us prove \cref{t--h--constant--star}. So we assume that $f$ is $\D$-$*$-normed and $\bc^*_\sigma$ is an isomorphism. Consider the commutative diagram
  \[
    \begin{tikzcd}
      &
      T \ar[dl, "{(\id_T, \eta_f \circ f)}", swap] \ar[dr, "\diag{f}"] \\
      T \times_S T \ar[rr, "{(\lp{f},\mu_f)}"] \ar[dr, "\lp{f}", swap] &
      &
      T \times_S T \ar[dl, "\lp{f}"] \\
      &
      T.
    \end{tikzcd}
  \]
  By hypothesis, the middle horizontal arrow $(\lp{f},\mu_f)$ is an isomorphism. In particular, this map is $\D$-proper, giving us a canonical identification $(\lp{f},\mu_f)_! \iso (\lp{f},\mu_f)_*$ (\cref{t--p--proper-twist}). We thus have a natural identification
  \begin{align}
    \label{t--h--constant--iso}
    (\lp{f})_*(\diag{f})_!
    &\iso (\lp{f})_*(\lp{f},\mu_f)_!(\id_T,\eta_f \circ f)_! \\
    &\iso (\lp{f})_*(\lp{f},\mu_f)_*(\id_T,\eta_f \circ f)_! \notag \\
    &\iso (\lp{f})_*(\id_T,\eta_f \circ f)_! \notag.
  \end{align}
  We next consider the cartesian squares
  \[
    \begin{tikzcd}
      T \ar[r, "f"] \ar[d, "{(\id_T,\eta_f \circ f)}", swap] &
      S \ar[d, "\eta_f"] \\
      T \times_S T \ar[r, "\rp{f}"] \ar[d, "\lp{f}", swap] &
      T \ar[d, "f"] \\
      T \ar[r, "f"] &
      S.
    \end{tikzcd}
  \]
  Let us denote the upper square by $\tau$; the lower one is the square $\sigma$ in the statement. Then we have
  \begin{align*}
    \primtwist_f
    ={} & (\lp{f})_*(\diag{f})_!(\unit_T) \\[0.8ex]
    \smash{\lbliso{\cref{t--h--constant--iso}}}{} & (\lp{f})_*(\id_T,\eta_f \circ f)_!(\unit_T) \\
    \iso{} & (\lp{f})_*(\id_T,\eta_f \circ f)_!f^*(\unit_S) \\[0.5ex]
    \smash{\lbliso{\bc^!_\tau}}{} & (\lp{f})_*(\rp{f})^*(\eta_f)_!(\unit_S) \\[0.5ex]
    \smash{\lbliso{\bc^*_\sigma}}{} & f^*f_*(\eta_f)_!(\unit_S). \qedhere
  \end{align*}
\end{proof}

\begin{lemma}
  \label{t--h--dualizable}
  Let $f : T \to S$ be a map in $\S$. Suppose that $f$ is $\D$-$*$-able and $\D$-$\sharp$-able and satisfies both the $\D$-$*$- and $\D$-$\sharp$-projection formulas. Then $f_\sharp(\unit_T)$ is a dualizable object of $\D(S)$, with dual $f_*(\unit_T)$.
\end{lemma}

\begin{proof}
  The projection formulas imply that the functors $f_\sharp f^* : \D(S) \to \D(S)$ and $f_*f^* : \D(S) \to \D(S)$ are $\D(S)$-linear, and hence given by tensoring with the objects $f_\sharp(\unit_T)$ and $f_*(\unit_T)$, respectively. We also have that the adjunction between these functors is $\D(S)$-linear, which translates into the claimed duality between $f_\sharp(\unit_T)$ and $f_*(\unit_T)$.
\end{proof}

\begin{theorem}
  \label{t--h--main}
  Let $S \in \S$ and let $f : T \to S$ be a left divisible $\AA_2$-monoid over $S$. Suppose that $f$ is $\D$-suave and $\D$-prim and that $\D(S)$ is idempotent complete. Then, if $\D$-smooth or $\D$-terse, it is $\D$-Poincar\'e.
\end{theorem}

\begin{proof}
  The two cases are formally dual; let us make the argument in the case that $f$ is $\D$-smooth. Then it remains to show that  $\primtwist_f$ is invertible. By \cref{t--p--dualizable}, this is equivalent to $\primtwist_f$ being dualizable, which is furthermore equivalent to $\poincaretwist_f = \suavetwist_f^{-1} \otimes \primtwist_f$ being dualizable. Since $f$ is $\D$-smooth and $\D$-prim, it satisfies the $\D$-$*$- and $\D$-$\sharp$-projection formulas and $\D$-$*$-- and $\D$-$\sharp$--base change (\cref{t--p--pr-bc}). Base change allows us to apply \cref{t--h--constant}, and we find that $\poincaretwist_f \iso f^*(\o\poincaretwist_f)$ for some $\o\poincaretwist_f \in \D(S)$; now it suffices to show that $\o\poincaretwist_f$ is dualizable in $\D(S)$. The projection formulas on the other hand allow us to apply \cref{t--h--dualizable}; in combination with the isomorphism $\Nm_f : f_\sharp(\poincaretwist_f) \iso f_*(\unit_T)$, we find that $f_\sharp(\poincaretwist_f)$ is dualizable in $\D(S)$. Finally, we have
  \[
    f_\sharp(\poincaretwist_f) \iso f_\sharp f^*(\o\poincaretwist_f) \smash{\lbliso{\pr^\sharp_f}} f_\sharp(\unit_T) \otimes \o\poincaretwist_f,
  \]
  and the maps $S \lblto{\eta_f} T \lblto{f} S$ induce a retraction $\unit_S \to f_\sharp(\unit_T) \to \unit_S$. Since dualizability in $\D(S)$ is closed under retracts by our idempotent completeness hypothesis, we deduce that $\o\poincaretwist_f$ is dualizable, as desired .
\end{proof}


\subsection{Norm and base change}
\label{t--b}

This subsection and the next concern the behavior of the norm maps defined in \cref{t--d} with respect to the fundamental geometric operations on maps in $\S$: base change and composition. We begin here with base change; until further notice, we suppose given a cartesian square $\sigma :$
\[
  \begin{tikzcd}
    V \ar[r, "q"] \ar[d, "g", swap] &
    T \ar[d, "f"] \\
    U \ar[r, "p"] &
    S
  \end{tikzcd}
\]
in $\S$.

\begin{construction}
  \label{t--b--star}
  Suppose that the maps $f$ and $g$ are $\D$-$*$-normed. By \itemref{t--d--kernel-weak-adjoints}{prim}, the map $[\unit_T]_S : [S]_S \to [T]_S$ in $\K_\D(S)$ has weak right adjoint $[\primtwist_f]_S : [T]_S \to [S]_S$, and the map $[\unit_V]_U : [U]_U \to [V]_U$ in $\K_\D(U)$ has weak right adjoint $[\primtwist_g]_U : [V]_U \to [U]_U$. We have a functor of $(\infty,2)$-categories $p^\dstar : \K_\D(S) \to \K_\D(U)$ (\itemitemref{t--f}{K-functoriality}{pull}), which carries the map $[\unit_T]_S$ to the map $[\unit_V]_U$ and carries the map $[\primtwist_f]_S$ to the map $[q^*(\primtwist_f)]_U : [V]_U \to [U]_U$. Applying \cref{t--a--functoriality}, we obtain a canonical map
  \[
    \bc^\primtwist_\sigma : q^*(\primtwist_f) \to \primtwist_g
  \]
  in $\D(V)$. Furthermore, \cref{t--a--functoriality-composition,t--a--functoriality-transformation} give us the following commutative diagrams involving this map:
    \[
    \begin{tikzcd}
      p^*f_!(- \otimes \primtwist_f) \ar[r, "\bc^!_\sigma"] \ar[d, "\Nm^*_f", swap] &
      g_!q^*(- \otimes \primtwist_f) \ar[r, "\sim"] &
      g_!(q^*(-) \otimes q^*(\primtwist_f)) \ar[r, "\bc^\primtwist_\sigma"] &
      g_!(q^*(-) \otimes \primtwist_g) \ar[d, "\Nm^*_g"] \\
      p^*f_* \ar[rrr, "\bc^*_\sigma"] &
      &
      &
      g_*q^*
    \end{tikzcd}
  \]
  \[
    \begin{tikzcd}
      p_!g_!(- \otimes \primtwist_g) \ar[d, "\Nm^*_g", swap] &
      p_!g_!(- \otimes q^*(\primtwist_f))  \ar[l, "\bc^\primtwist_\sigma", swap] \ar[r, "\sim"] &
      f_!q_!(- \otimes q^*(\primtwist_f)) \ar[r, "\pr^!_q"] &
      f_!(q_!(-) \otimes \primtwist_f) \ar[d, "\Nm^*_f"] \\
      p_!g_* \ar[rrr, "\ex^*_{\sigma^\top}"] &
      &
      &
      f_*q_!
    \end{tikzcd}
  \]
  (here $\bc^!_\sigma$ and $\bc^*_\sigma$ are the base change maps associated to $\sigma$, $\pr^!_q$ is the projection isomorphism associated to $q$, and $\ex^*_{\sigma^\top}$ is the exchange map associated to the transpose of $\sigma$; see \itemref{t--f}{D}). More specifically, for both diagrams we apply \cref{t--a--functoriality-composition} to the composition of $p^\dstar : \K_\D(S) \to \K_\D(U)$ and $\Psi_{\D,U} : \K_\D(U) \to \Cat_\infty$, and for the upper diagram (resp. lower diagram) we apply \cref{t--a--functoriality-transformation} to the natural transformation $\Psi_{\D,S} \to \Psi_{\D,U} \circ p^\dstar$ (resp. $\Psi_{\D,U} \circ p^\dstar \to \Psi_{\D,S}$) described in \itemref{t--f}{K-functoriality}.
\end{construction}

\begin{remark}
  \label{t--b--star-equiv}
  In the situation of \cref{t--b--star}, if $f$ is $\D$-prim, then the map $\bc^\primtwist_\sigma$ is an isomorphism, by \cref{t--a--strict-functoriality}.
\end{remark}

\begin{construction}
  \label{t--b--sharp}
  Suppose that the maps $f$ and $g$ are $\D$-$\sharp$-normed. Formally dualizing \cref{t--b--star}, we obtain a canonical map
  \[
    \bc^\suavetwist_\sigma : \suavetwist_g \to q^*(\suavetwist_f)
  \]
  in $\D(V)$, making the following diagrams commute:
  \[
    \begin{tikzcd}
      p^*f_\sharp  \ar[d, "\Nm^\sharp_f", swap] &
      &
      &
      g_\sharp q^* \ar[d, "\Nm^\sharp_g"] \ar[lll, "\bc^\sharp_\sigma", swap] \\
      p^*f_!(- \otimes \suavetwist_f) &
      g_!q^*(- \otimes \suavetwist_f) \ar[l, "\bc^!_\sigma", swap] &
      g_!(q^*(-) \otimes q^*(\suavetwist_f)) \ar[l, "\sim", swap]  &
      g_!(q^*(-) \otimes \suavetwist_g) \ar[l, "\bc^\suavetwist_\sigma", swap]
    \end{tikzcd}
  \]
  \[
    \begin{tikzcd}
      p_!g_\sharp  \ar[d, "\Nm^*_g", swap] &
      &
      &
      f_\sharp q_! \ar[lll, "\ex^\sharp_{\sigma^\top}", swap] \ar[d, "\Nm^*_f"] \\
      p_!g_!(- \otimes \suavetwist_g) \ar[r, "\bc^\suavetwist_\sigma"] &
      p_!g_!(- \otimes q^*(\suavetwist_f)) &
      f_!q_!(- \otimes q^*(\suavetwist_f)) \ar[l, "\sim", swap] &
      f_!(q_!(-) \otimes \suavetwist_f). \ar[l, "\pr^!_q", swap]
    \end{tikzcd}
  \]
\end{construction}

\begin{remark}
  \label{t--b--sharp-equiv}
  In the situation of \cref{t--b--sharp}, if $f$ is $\D$-suave, then the map $\bc^\suavetwist_\sigma$ is an isomorphism, by \cref{t--a--strict-functoriality}.
\end{remark}

\begin{remark}
  \label{t--b--sharpstar}
  Suppose that the maps $f$ and $g$ are $\D$-normed. Combining \cref{t--b--star,t--b--sharp}, we obtain the following commutative diagrams:
  \[
    \begin{tikzcd}
      p^*f_\sharp(- \otimes \primtwist_f) \ar[d, "\Nm_f", swap] &
      g_\sharp q^*(- \otimes \primtwist_f) \ar[r, "\sim"] \ar[l, "\bc^\sharp_\sigma", swap] &
      g_\sharp(q^*(-) \otimes q^*(\primtwist_f)) \ar[r, "\bc^\primtwist_\sigma"] &
      g_\sharp(q^*(-) \otimes \primtwist_g) \ar[d, "\Nm_g"] \\
      p^*f_*(- \otimes \suavetwist_f) \ar[r, "\bc^*_\sigma"] &
      g_* q^*(- \otimes \suavetwist_f) \ar[r, "\sim"] &
      g_*(q^*(-) \otimes q^*(\suavetwist_f))  &
      g_*(q^*(-) \otimes \suavetwist_g) \ar[l, "\bc^\suavetwist_\sigma", swap]
    \end{tikzcd}
  \]
  \[
    \begin{tikzcd}
      p_!g_\sharp(- \otimes \primtwist_g) \ar[d, "\Nm_g", swap] &
      f_\sharp q_!(- \otimes \primtwist_g) \ar[l, "\ex^\sharp_{\sigma^\top}", swap] &
      f_\sharp q_!(- \otimes q^*(\primtwist_f)) \ar[l, "\bc^\primtwist_\sigma", swap] &
      f_\sharp(q_!(-) \otimes \primtwist_f) \ar[l, "\pr^!_q", swap] \ar[d, "\Nm_f"] \\
      p_!g_*(- \otimes \suavetwist_g) \ar[r, "\ex^*_\sigma"] &
      f_*q_!(- \otimes \suavetwist_g) \ar[r, "\bc^\suavetwist_\sigma"] &
      f_*q_!(- \otimes q^*(\suavetwist_f))  \ar[r, "\pr^!_q"] &
      f_*(q_!(-) \otimes \suavetwist_f).
    \end{tikzcd}
  \]
\end{remark}

\begin{construction}
  \label{t--b--tate}
  Suppose that the maps $f$ and $g$ are $\D$-Tate. Then the map $\bc^\suavetwist_\sigma : \suavetwist_g \to q^*(\suavetwist_f)$ of \cref{t--b--sharp} induces a map $\bc^\suavetwist_\sigma : q^*(\suavetwist_f^{-1}) \to \suavetwist_g^{-1}$. Tensoring this with the map $\bc^\primtwist_\sigma : q^*(\primtwist_f) \to \primtwist_g$ of \cref{t--b--star}, we obtain a map
  \[
    \bc^\poincaretwist_\sigma : q^*(\poincaretwist_f) \to \poincaretwist_g.
  \]
  Then the commutative diagrams of \cref{t--b--sharpstar} may be rewritten as the following:
  \[
    \begin{tikzcd}
      p^*f_\sharp(- \otimes \poincaretwist_f)  \ar[d, "\Nm_f", swap] &
      g_\sharp q^*(- \otimes \poincaretwist_f) \ar[l, "\bc^\sharp_\sigma", swap] \ar[r, "\sim"] &
      g_\sharp (q^*(-) \otimes q^*(\poincaretwist_f)) \ar[r, "\bc^\poincaretwist_\sigma"] &
      g_\sharp(q^*(-) \otimes \poincaretwist_g) \ar[d, "\Nm_g"] \\
      p^*f_* \ar[rrr, "\bc^*_\sigma"] &
      &
      &
      g_*q^*
    \end{tikzcd}
  \]
    \[
    \begin{tikzcd}
      p_!g_\sharp(- \otimes \poincaretwist_g) \ar[d, "\Nm_g", swap] &
      p_!g_\sharp(- \otimes q^*(\poincaretwist_f))  \ar[l, "\bc^\poincaretwist_\sigma", swap] \ar[r, "\sim"] &
      f_\sharp q_!(- \otimes q^*(\poincaretwist_f)) \ar[r, "\pr^!_q"] &
      f_\sharp (q_!(-) \otimes \poincaretwist_f) \ar[d, "\Nm_f"] \\
      p_!g_* \ar[rrr, "\ex^*_{\sigma^\top}"] &
      &
      &
      f_*q_!.
    \end{tikzcd}
  \]
\end{construction}

\begin{remark}
  \label{t--b--poincare-equiv}
  In the situation of \cref{t--b--star}, if $f$ is both $\D$-prim and $\D$-suave (e.g. if it is $\D$-Poincar\'e), then the map $\bc^\poincaretwist_\sigma$ is an isomorphism, by combining \cref{t--b--star-equiv,t--b--sharp-equiv}.
\end{remark}

\begin{proposition}[{\citebare[Lemmas 4.5.9 and 4.6.3]{heyer-mann--six}}]
  \label{t--b--closure}
  The classes of $\D$-prim, $\D$-terse, $\D$-proper, $\D$-suave, $\D$-smooth, $\D$-\'etale, and $\D$-Poincar\'e maps in $\S$ are closed under base change along all maps in $\S$.
\end{proposition}

\begin{proof}
  Suppose in our cartesian square
  \[
    \begin{tikzcd}
      V \ar[r, "q"] \ar[d, "g", swap] &
      T \ar[d, "f"] \\
      U \ar[r, "p"] &
      S
    \end{tikzcd}
  \]
  that $f$ is $\D$-prim, i.e. that the map $[\unit_T]_S : [S]_S \to [T]_S$ in $\K_\D(S)$ admits a right adjoint. As noted in \cref{t--b--star}, the functor of $(\infty,2)$-categories $p^\dstar : \K_\D(S) \to \K_\D(U)$ carries $[\unit_T]_S$ to the map $[\unit_V]_U : [U]_U \to [V]_U$. Since functors of $(\infty,2)$-categories preserve adjunctions, it follows that $[\unit_V]_U$ admits a right adjoint, meaning that $g$ is $\D$-prim. Moreover, we have in this case that the map $\bc^\primtwist_\sigma : q^*(\primtwist_f) \to \primtwist_g$ is an isomorphism (\cref{t--b--star-equiv}), which implies that if $f$ is $\D$-terse, then so is $g$. The case of $\D$-suave and $\D$-smooth maps is dual, and then it follows for $\D$-proper, $\D$-\'etale, and $\D$-Poincar\'e maps.
\end{proof}

In \cref{e}, we will also use a descent result from \cite{heyer-mann--six}, converse to \cref{t--b--closure}. We end this subsection by recalling this, now unfixing the cartesian square $\sigma$ from the first part of the subsection.

\begin{proposition}[{\citebare[Lemma 4.5.7]{heyer-mann--six}}]
  \label{t--b--descent}
  The classes of $\D$-prim, $\D$-terse, $\D$-suave, $\D$-smooth, and $\D$-Poincar\'e maps in $\S_!$ are all local with respect to $\D$-$*$-covers (\itemitemref{t--f}{D}{cover}) of the target: that is, for any of these classes $\mathsf{p}$, and for any $\D$-$!$-able map $f : T \to S$ in $\S$ and a $\D$-$*$-cover $\cat{U}$ of $S$, if for every cartesian square
  \[
    \begin{tikzcd}
      V \ar[r, "q"] \ar[d, "g", swap] &
      T \ar[d, "f"] \\
      U \ar[r, "p"] &
      S
    \end{tikzcd}
  \]
  in $\S$ where $p \in \cat{U}$ the map $g$ lies in $\mathsf{p}$, then the map $f$ also lies in $\mathsf{p}$.
\end{proposition}

\begin{proof}
  For primness and suaveness, we refer to the cited result of \cite{heyer-mann--six}. The claim then follows for $\D$-terseness, $\D$-smoothness, and $\D$-Poincar\'eness, thanks to \cref{t--b--star-equiv,t--b--sharp-equiv} and the fact that invertibility in the limit of a diagram of symmetric monoidal $\infty$-categories and symmetric monoidal functors can be checked after projection to each factor in the diagram.
\end{proof}


\subsection{Norm and composition}
\label{t--c}

We now turn to the interaction of our norm maps with composition in $\S$. This will be particularly important in our discussion of the universal property and functoriality of Tate cohomology in the next two subsections. Until further notice we fix a commutative triangle $\tau :$
\[
  \begin{tikzcd}[row sep=small]
    U \ar[rr, "g"] \ar[dr, "h", swap] &
    &
    T \ar[dl, "f"] \\
    &
    S
  \end{tikzcd}
\]
in $\S$.

\begin{construction}
  \label{t--c--prim-comp}
  Suppose that $f$ and $g$ are $\D$-$*$-normed, and note that $h$ is then also $\D$-$*$-normed. Consider the following diagram of functors from $\D(U)$ to $\D(S)$:
  \begin{equation}
    \label{t--c--prim-comp--diagram}
    \begin{tikzcd}
      f_!g_!(- \otimes \primtwist_g \otimes g^*(\primtwist_f)) \ar[r, "\pr^!_g"] \ar[dd, "\sim", swap] &
      f_!(g_!(- \otimes \primtwist_g) \otimes \primtwist_f) \ar[r, "\Nm^*_g"] \ar[d, "\Nm^*_f", swap] &
      f_!(g_*(-) \otimes \primtwist_f) \ar[d, "\Nm^*_f"]
      \\
      &
      f_*(g_!(- \otimes \primtwist_g)) \ar[r, "\Nm^*_g"] &
      f_*g_* \ar[d, "\sim"] \\
      h_!(- \otimes \primtwist_g \otimes g^*(\primtwist_f)) \ar[r, "\comp^\primtwist_{f,g}", dashed]&
      h_!(- \otimes \primtwist_h) \ar[r, "\Nm^*_h"] &
      h_*
    \end{tikzcd}
  \end{equation}
  Recall that $\pr^!_g$ denotes the projection isomorphism (\itemitemref{t--f}{D}{shriek}), and observe that the upper right square commutes by naturality. We will now construct a map $\comp^\primtwist_{f,g} : \primtwist_g \otimes g^*(\primtwist_f) \to \primtwist_h$ in $\D(U)$, defining the dashed arrow above, such that the other piece of the diagram commutes as well.

  Consider the map $[\unit_U]_T : [T]_T \to [U]_T$ in $\K_\D(T)$, i.e. the unit map for the canonical commutative algebra structure of $[U]_T$ (\itemitemref{t--f}{K-monoidal}{formulas}). Applying the lax symmetric monoidal functor $f_\dstar : \K_\D(T) \to \K_\D(S)$ (\itemitemref{t--f}{K-functoriality}{push}), we obtain a map of commutative algebra objects $f_\dstar([\unit_U]_T) : [T]_S \to [U]_S$ in $\K_\D(S)$. This map in particular preserves units, meaning that we have a canonical isomorphism
  \[
    [\unit_U]_S \iso f_\dstar([\unit_U]_T) \circ [\unit_T]_S
  \]
  (which is also straightforward to construct directly). By \itemref{t--d--kernel-weak-adjoints}{prim}, the maps $[\unit_T]_S$, $[\unit_U]_T$, and $[\unit_U]_S$ here have weak right adjoints given by $[\primtwist_f]_S$, $[\primtwist_g]_T$, and $[\primtwist_h]_S$. Let $\coun' : [\unit_T]_S \circ [\primtwist_f]_S \to \id_{[T]_S}$, $\coun'' : [\unit_U]_T \circ [\primtwist_g]_T \to \id_{[U]_T}$, and $\coun : [\unit_U]_S \circ [\primtwist_h]_S \to \id_{[U]_S}$ be the associated counit maps. Noting that there is an induced map
  \[
    f_\dstar(\coun'') : f_\dstar([\unit_U]_T) \circ f_\dstar([\primtwist_g]_T) \to f_\dstar(\id_{[U]_T}) \iso \id_{[U]_S},
  \]
  we may apply \cref{t--c--categorical-comp-right} to obtain a map $\coun''' : [\unit_U]_S \circ [\primtwist_f]_S \circ f_\dstar([\primtwist_g]_T) \to \id_{[U]_S}$ and a unique map $\comp : [\primtwist_f]_S \circ f_\dstar([\primtwist_g]_T) \to [\primtwist_h]_S$ making the diagram
  \[
    \begin{tikzcd}
      {[\unit_U]_S \circ [\primtwist_f]_S \circ f_\dstar([\primtwist_g]_T)} \ar[r, "\comp"] \ar[dr, "\coun'''", swap] &
      {[\unit_U]_S \circ [\primtwist_h]_S} \ar[d, "\coun"] \\
      &
      \id_{[U]_S}
    \end{tikzcd}
  \]
  commute. Unwinding definitions gives an identification $[\primtwist_f]_S \circ f_\dstar([\primtwist_g]_T) \iso [\primtwist_g \otimes g^*(\primtwist_f)]_S$, under which $\mu$ corresponds to the desired map $\comp^\primtwist_{f,g} : \primtwist_g \otimes g^*(\primtwist_f) \to \primtwist_h$ in $\D(U)$. Applying $\Psi_{\D,S}$ to the above commutative triangle gives the claimed commutativity in \cref{t--c--prim-comp--diagram}.
\end{construction}

\begin{proposition}
  \label{t--c--prim-equiv}
  Suppose that $f$ and $g$ are $\D$-$*$-normed, and also that one of the following further conditions is satisfied:
  \begin{enumerate}
  \item \label{t--c--prim-equiv--prim}
    $f$ is $\D$-prim;
  \item \label{t--c--prim-equiv--terse}
    $g$ is $\D$-Poincar\'e.
  \end{enumerate}
  Then the map $\comp^\primtwist_{f,g} : \primtwist_g \otimes g^*(\primtwist_f) \to \primtwist_h$ of \cref{t--c--prim-comp} is an isomorphism.
\end{proposition}

\begin{proof}
  We will apply \cref{t--c--categorical-equiv-right}. If $f$ is $\D$-prim, then \itemref{t--c--categorical-equiv-right}{post} is satisfied. Suppose on the other hand that $g$ is $\D$-Poincar\'e. Then we claim that \itemref{t--c--categorical-equiv-right}{pre} is satisfied, i.e. that $f_\dstar([\primtwist_g]_T) : [U]_S \to [T]_S$ admits a right adjoint in $\K_\D(S)$. It suffices to show that $[\primtwist_g]_T : [U]_T \to [T]_T$ admits a right adjoint in $\K_\D(T)$. For this we will use \cref{t--p--right-adjoint-criterion}. The functors
  \[
    \uMap_{\K_\D(T)}([T],[U]) \to \uMap_{\K_\D(T)}([T],[T]), \quad
    \uMap_{\K_\D(T)}([U],[U]) \to \uMap_{\K_\D(T)}([U],[T])
  \]
  given by composition with $[\primtwist_g]_T$ identify with the functors
  \[
    g_!(- \otimes \primtwist_g) : \D(U) \to \D(T), \quad
    (\lp{g})_!(- \otimes \rp{g}^*(\primtwist_g)) : \D(U \times_T U) \to \D(U).
  \]
  Since $g$ is $\D$-smooth, the first of these identifies with $g_\sharp(- \otimes \poincaretwist_g)$, which has right adjoint given by $g^*(-) \otimes \poincaretwist_g^{-1} : \D(T) \to \D(U)$. By \cref{t--b--closure}, the base change $\lp{g}$ of $g$ is also $\D$-Poincar\'e, and we have isomorphisms $\suavetwist_{\lp{g}} \iso \rp{g}^*(\suavetwist_g)$ and $\primtwist_{\lp{g}} \iso \rp{g}^*(\primtwist_g)$ by \cref{t--b--star-equiv,t--b--sharp-equiv}. Hence, the second functor above similarly has right adjoint given by $\lp{g}^*(-) \otimes \rp{g}^*(\poincaretwist_g^{-1}) : \D(U) \to \D(U \times_T U)$. Finally, the composition
  \[
    \smash{[U]_T \lblto{[\poincaretwist_g^{-1}]_T} [T]_T \lblto{[\primtwist_g]_T} [U]_T}
  \]
  identifies with the map $[\lp{g}^*(\primtwist_g) \otimes \rp{g}^*(\poincaretwist_g^{-1})]_T : [U]_T \to [U]_T$, verifying that \itemref{t--p--right-adjoint-criterion}{criterion} holds and finishing the proof.
\end{proof}

\begin{construction}
  \label{t--c--suave-comp}
  Suppose that $f$ and $g$ are $\D$-$\sharp$-normed, and hence so is $h$. Formally dualizing \cref{t--c--prim-comp} gives a map $\comp^\suavetwist_{f,g} : \suavetwist_h \to \suavetwist_g \otimes g^*(\suavetwist_f)$ in $\D(U)$ and the following commutative diagram of functors from $\D(U)$ to $\D(S)$:
  \begin{equation}
    \label{t--c--suave-comp--diagram}
    \begin{tikzcd}[row sep=3ex]
      h_\sharp \ar[r, "\Nm^\sharp_h"] \ar[d, "\sim", swap] &
      h_!(- \otimes \suavetwist_h) \ar[r, "\comp^\suavetwist_{f,g}"] &
      h_!(- \otimes \suavetwist_g \otimes g^*(\suavetwist_f))  \ar[dd, "\sim"] \\
      f_\sharp g_\sharp  \ar[r, "\Nm^\sharp_g"] \ar[d, "\Nm^\sharp_f", swap] &
      f_\sharp(g_!(- \otimes \suavetwist_g))  \ar[d, "\Nm^\sharp_f"] \\
      f_!(g_\sharp(-) \otimes \suavetwist_f) \ar[r, "\Nm^\sharp_g"] &
      f_!(g_!(- \otimes \suavetwist_g) \otimes \suavetwist_f) \ar[r, "\pr^!_g"] &
      f_!g_!(- \otimes \suavetwist_g \otimes g^*(\suavetwist_f)).
    \end{tikzcd}
  \end{equation}
\end{construction}

\begin{proposition}
  \label{t--c--suave-equiv}
  Suppose that $f$ and $g$ are $\D$-$\sharp$-normed, and also that one of the following further conditions is satisfied:
  \begin{enumerate}
  \item $f$ is $\D$-suave;
  \item $g$ is $\D$-Poincar\'e.
  \end{enumerate}
  Then the map $\comp^\suavetwist_{f,g} : \suavetwist_h \to \suavetwist_g \otimes g^*(\suavetwist_f)$ of \cref{t--c--suave-comp} is an isomorphism.
\end{proposition}

\begin{proof}
  This is formally dual to \cref{t--c--prim-equiv}.
\end{proof}

\begin{remark}
  \label{t--c--norm-comp}
  Suppose that $f$ and $g$ are $\D$-normed, and hence so is $h$. Combining the commutative diagrams \cref{t--c--prim-comp--diagram} and \cref{t--c--suave-comp--diagram}, and using the compatibility of $*$-norm and $\sharp$-norm maps with projection maps (\cref{t--d--norm-linearity}), we obtain the following commutative diagram of functors from $\D(U)$ to $\D(S)$:
  \[
    \begin{tikzcd}[column sep=1.05em]
      h_\sharp(- \otimes \primtwist_g \otimes g^*(\primtwist_f)) \ar[r, "\comp^\primtwist_{f,g}"] \ar[d, "\sim", swap] &
      h_\sharp(- \otimes \primtwist_h) \ar[r, "\Nm_h"] &
      h_*(- \otimes \suavetwist_h) \ar[r, "\comp^\suavetwist_{f,g}"] &
      h_*(- \otimes \suavetwist_g \otimes g^*(\suavetwist_f)) \ar[dd, "\sim"] \\
      f_\sharp g_\sharp(- \otimes \primtwist_g \otimes g^*(\primtwist_f)) \ar[r, "\pr^\sharp_g"] &
      f_\sharp(g_\sharp(- \otimes \primtwist_g) \otimes \primtwist_f) \ar[r, "\Nm_g"] \ar[d, "\Nm_f", swap] &
      f_\sharp(g_*(- \otimes \suavetwist_g) \otimes \primtwist_f) \ar[d, "\Nm_f"] \\
      &
      f_*(g_\sharp(- \otimes \primtwist_g) \otimes \suavetwist_f) \ar[r, "\Nm_g", swap] &
      f_*(g_*(- \otimes \suavetwist_g) \otimes \suavetwist_f) \ar[r, "\pr^*_g"] &
      f_*g_*(- \otimes \suavetwist_g \otimes g^*(\suavetwist_f))
    \end{tikzcd}
  \]
\end{remark}

\begin{construction}
  \label{t--c--poincare-comp}
  Suppose that $f$ and $g$ are $\D$-Tate maps, and suppose that the map $\comp^\suavetwist_{f,g} : \suavetwist_h \to \suavetwist_g \otimes g^*(\suavetwist_f)$ of \cref{t--c--suave-comp} is an isomorphism (e.g. one of the conditions of \cref{t--c--suave-equiv} is satisfied), so that $h$ is also $\D$-Tate, with  $\suavetwist_g^{-1} \otimes g^*(\suavetwist_f^{-1}) \iso \suavetwist_h^{-1}$. Tensoring this last identification with the map $\comp^\primtwist_{f,g} : \primtwist_g \otimes g^*(\primtwist_f) \to \primtwist_h$ of \cref{t--c--prim-comp}, we obtain a map $\comp^\poincaretwist_{f,g} : \poincaretwist_g \otimes g^*(\poincaretwist_f) \to \poincaretwist_h$ in $\D(U)$. Then the commutative diagram of \cref{t--c--norm-comp} induces the following commutative diagram of functors from $\D(U)$ to $\D(S)$:
  \[
    \begin{tikzcd}[row sep=3ex]
      f_\sharp g_\sharp (- \otimes \poincaretwist_g \otimes g^*(\poincaretwist_f)) \ar[r, "\pr^\sharp_g"] \ar[dd, "\sim", swap] &
      f_\sharp(g_\sharp(- \otimes \poincaretwist_g) \otimes \poincaretwist_f) \ar[r, "\Nm_g"] \ar[d, "\Nm_f", swap] &
      f_\sharp(g_*(-) \otimes \poincaretwist_f) \ar[d, "\Nm_f"]
      \\
      &
      f_*(g_\sharp(- \otimes \poincaretwist_g)) \ar[r, "\Nm_g"] &
      f_*g_* \ar[d, "\sim"] \\
      h_\sharp(- \otimes \poincaretwist_g \otimes g^*(\poincaretwist_f)) \ar[r, "\comp^\poincaretwist_{f,g}"]&
      h_\sharp(- \otimes \poincaretwist_h) \ar[r, "\Nm_h"] &
      h_*.
    \end{tikzcd}
  \]
\end{construction}

\begin{proposition}
  \label{t--c--poincare-equiv}
  Suppose that $f$ and $g$ are $\D$-Tate, and also that one of the following further conditions is satisfied:
  \begin{enumerate}
  \item $f$ is $\D$-prim and $\D$-suave;
  \item $g$ is $\D$-Poincar\'e.
  \end{enumerate}
  Then $h$ is $\D$-Tate and the map $\comp^\poincaretwist_{f,g} : \poincaretwist_g \otimes g^*(\poincaretwist_f) \to \poincaretwist_h$ of \cref{t--c--poincare-comp} is an isomorphism.
\end{proposition}

\begin{proof}
  This follows from combining \cref{t--c--prim-equiv,t--c--suave-equiv}.
\end{proof}

\begin{proposition}[{\citebare[Lemmas 4.5.9 and 4.6.3]{heyer-mann--six}}]
  \label{t--c--closure}
  The classes of $\D$-prim, $\D$-terse, $\D$-proper, $\D$-suave, $\D$-smooth, $\D$-\'etale, and $\D$-Poincar\'e maps in $\S$ are closed under composition.
\end{proposition}

\begin{proof}
  Suppose in our commutative triangle
  \[
    \begin{tikzcd}[row sep=small]
      U \ar[rr, "g"] \ar[dr, "h", swap] &
      &
      T \ar[dl, "f"] \\
      &
      S
    \end{tikzcd}
  \]
  that $f$ and $g$ are $\D$-prim, i.e. that the maps $[\unit_T]_S : [S]_S \to [T]_S$ in $\K_\D(S)$ and $[\unit_U]_T : [T]_T \to [U]_T$ in $\K_\D(T)$ admit right adjoints. As noted in \cref{t--c--prim-comp}, the map $[\unit_U]_S : [S]_S \to [T]_S$ identifies with the composition $f_\dstar([\unit_U]_T) \circ [\unit_T]_S$, and so then also admits a right adjoint, meaning that $h$ is $\D$-prim. Moreover, we have in this case by \cref{t--c--prim-equiv} that the map $\comp^\primtwist_{f,g} : \primtwist_g \otimes g^*(\primtwist_f) \to \primtwist_h$ is an isomorphism, implying that if $f$ and $g$ are $\D$-terse, then so is $h$. The case of $\D$-suave and $\D$-smooth maps is dual, and then it follows for $\D$-Poincar\'e maps. The cases of $\D$-proper and $\D$-\'etale maps also follow formally, given that $\diag{h} : U \to U \times_S U$ factors as the composition of $\diag{g} : U \to U \times_T U$ and the base change $\diag{f}' : U \times_T U \to U \times_S U$ of $\diag{f} : T \to T \times_S T$, and the fact that these conditions are stable under base change (\cref{t--b--closure}).
\end{proof}

Parallel to the previous subsection, let us end our discussion here by unfixing the commutative triangle $\tau$ from above and recalling a descent result from \cite{heyer-mann--six} which gives a partial converse to \cref{t--c--closure}, in the suave and smooth cases (these are the cases that will be relevant in \cref{e}).

\begin{proposition}[{\citebare[Lemma 4.5.8]{heyer-mann--six}}]
  \label{t--c--descent}
  The class of $\D$-suave maps in $\S_!$ is local with respect to $\D$-suave $\D$-$*$-covers of the source: that is, for any $\D$-$!$-able map $f : T \to S$ in $\S$ and any $\D$-$*$cover (\itemitemref{t--f}{D}{cover}) $\cat{U}$ of $T$  consisting of $\D$-suave maps, if for every commutative triangle
  \[
    \begin{tikzcd}[row sep=small]
      U \ar[rr, "g"] \ar[dr, "h", swap] &
      &
      T \ar[dl, "f"] \\
      &
      S
    \end{tikzcd}
  \]
  in $\S$ where $g \in \cat{U}$ the map $h$ is $\D$-suave, then $f$ is also in $\D$-suave. The same statement with $\D$-suave everywhere replaced by $\D$-smooth also holds.
\end{proposition}

\begin{proof}
  For the suave case, we refer to the cited result of \cite{heyer-mann--six}. For the smooth case, we have to further check that, if $\suavetwist_h$ is invertible for all such commutative triangles, then so too is $\suavetwist_f$. Since $\cat{U}$ is a $\D$-$*$-cover, it suffices to show that $g^*(\suavetwist_f)$ is invertible for each $g \in \cat{U}$. This follows from the hypotheses and the fact that the map $\comp^\suavetwist_{f,g} : \suavetwist_h \to \suavetwist_g \otimes g^*(\suavetwist_f)$ is an isomorphism for all such $g$, by \cref{t--c--suave-equiv}.
\end{proof}


\subsection{Universal property of Tate cohomology}
\label{t--u}

Throughout this subsection, we assume that our three functor formalism is stable and presentable, i.e. is given by a lax symmetric monoidal functor $\D : \Span(\S,\S_!) \to \PrLst$. The stability in particular means that we have, for each $\D$-Tate map $f : T \to S$ in $\S$, a Tate cohomology functor $f_\t : \D(T) \to \D(S)$ (\cref{t--d--norm}). Our goal here is to show that, under certain hypotheses on $f$, the functor $f_\t$ is characterized by a universal property having to do with its vanishing on certain ``nilpotent'' objects. When it applies, this characterization will endow $f_\t$ with a lax symmetric monoidal structure, compatible with that of the cohomology functor $f_*$. These results are closely adapted from those of Nikolaus--Scholze in \cite[\textsection I.3]{nikolaus-scholze--TC}.

The following two standard definitions will be relevant here.

\begin{definition}
  \label{t--u--thick}
  Let $\X$ be an idempotent complete, stable $\infty$-category. A \emph{thick subcategory} of $\X$ is a full subcategory $\X_0 \subseteq \X$ that is closed under retracts, finite limits, and finite colimits.
\end{definition}

\begin{definition}
  \label{t--u--tensor-ideal}
  Let $\X$ be a symmetric monoidal $\infty$-category. A \emph{tensor ideal} of $\X$ is a full subcategory $\X_0 \subseteq \X$ such that for all $X,X' \in \X$, if $X \in \X_0$, then $X \otimes X' \in \X_0$.
\end{definition}

The relevance of these notions comes from the following construction.

\begin{definition}
  \label{t--u--nilpotent}
  Let $i : V \to T$ be a $\D$-$\sharp$-able map in $\S$. We let $\D(T)_{\nil{i}}$ denote the thick subcategory of $\D(T)$ generated by the objects $i_\sharp(M)$ for $M \in \D(V)$. We say that an object $L \in \D(T)$ is \emph{$i$-nilpotent} if it is contained in $\D(T)_{\nil{i}}$.
\end{definition}

\begin{lemma}
  \label{t--u--nilpotent-ideal}
  Let $i : V \to T$ be a map in $\S$ that satisfies the $\D$-$\sharp$-projection formula. Then $\D(T)_{\nil{i}}$ is a tensor ideal in $\D(T)$.
\end{lemma}

\begin{proof}
  We need to show for $L,L' \in \D(T)$ that, if $L$ is $i$-nilpotent, then $L \otimes L'$ is $i$-nilpotent. Fixing $L'$, the collection of $L$ for which this is satisfied forms a thick subcategory of $\D(T)$, so we need only consider the case where $L \iso i_\sharp(M)$ for $M \in \D(V)$. The projection map
  \[
     i_\sharp(M \otimes i^*(L')) \lblto{\pr^\sharp_i} i_\sharp(M) \otimes L' \iso L \otimes L'
  \]
  is an isomorphism by hypothesis, proving the claim.
\end{proof}

We may now formulate the main result of the subsection.

\begin{definition}
  \label{t--u--cover}
  Let $f : T \to S$ be a $\D$-Tate map in $\S$. We say that a map $i : V \to T$ in $\S$ is a \emph{$\D$-Tate cover} of $f$ if the following conditions hold:
  \begin{enumerate}
  \item \label{t--u--cover--homology}
    $i$ satisfies the $\D$-$\sharp$-projection formula;
  \item \label{t--u--cover--exhaustive}
    the functor $i^* : \D(T) \to \D(V)$ is conservative;
  \item \label{t--u--cover--vanishing}
    the functor $f_\t : \D(T) \to \D(S)$ vanishes on $\D(T)_{\nil{i}}$ (equivalently, the composition $f_\t \circ i_\sharp : \D(V) \to \D(S)$ vanishes).
  \end{enumerate}
\end{definition}

\begin{remark}
  \label{t--u--nilpotent-tate-vanishing-ad}
  See \cref{t--u--nilpotent-tate-vanishing} for a way to check that condition \cref{t--u--cover--vanishing} of \cref{t--u--cover} holds.
\end{remark}

\begin{theorem}
  \label{t--u--main}
  Let $f : T \to S$ be a $\D$-Tate map in $\S$ and let $i : V \to T$ be a $\D$-Tate cover of $f$. Then the following statements hold:
  \begin{enumerate}
  \item \label{t--u--main--initial}
    The canonical map $\can : f_* \to f_\t$ is initial among maps $\gamma : f_* \to F$ of functors from $\D(T)$ to $\D(S)$ such that $F$ is accessible and $F|_{\D(T)_{\nil{i}}} \iso 0$;
  \item \label{t--u--main--fiber-criterion}
    Suppose given a map $\gamma : f_* \to F$ of functors from $\D(T)$ to $\D(S)$ such that $F$ is accessible and $F|_{\D(T)_{\nil{i}}} \iso 0$, and let $\o\gamma : f_\t \to F$ be the unique transformation with an equivalence $\o\gamma \circ \can \iso \gamma$ (as per \cref{t--u--main--initial}). Then $\o\gamma$ is an isomorphism if and only if the functor $\fib(\gamma) : \D(T) \to \D(S)$ preserves colimits.
  \item \label{t--u--main--monoidal}
    There are lax symmetric monoidal structures on the Tate cohomology functor $f_\t : \D(T) \to \D(S)$ and on the canonical map $\can : f_* \to f_\t$ making the latter initial among maps $f_* \to F$ of lax symmetric monoidal functors from $\D(T)$ to $\D(S)$ such that $F$ is accessible and $F|_{\D(T)_{\nil{i}}} \iso 0$.
  \end{enumerate}
\end{theorem}

This will be a consequence of a more general categorical result, \cref{t--u--abstract} below, which is a minor variation of \cite[Theorems I.3.3(iii) and I.3.6(ii)]{nikolaus-scholze--TC}.

\begin{definition}
  \label{t--u--relatively-accessible}
  Let $\B$ be a presentable $\infty$-category and let $\B_0$ be a full subcategory of $\B$. Recall that, for $\kappa$ a regular cardinal, $\B^\kappa$ denotes the full subcategory of $\B$ spanned by the $\kappa$-compact objects, and let us write $\B^\kappa_0 := \B^\kappa \cap \B_0$; we say that $\B_0$ is \emph{relatively $\kappa$-accessible} in $\B$ if every object of $\B_0$ is the colimit in $\B$ of a $\kappa$-filtered diagram in $\B^\kappa_0$. We say that $\B_0$ is \emph{relatively accessible} in $\B$ if, for every regular cardinal $\kappa$, there exists a regular cardinal $\kappa' \ge \kappa$ such that $\B_0$ is relatively $\kappa'$-accessible in $\B$.
\end{definition}

\begin{lemma}
  \label{t--u--relatively-accessible-meaning}
  Let $\A$ and $\B$ be presentable $\infty$-categories, let $\B_0$ be a relatively accessible full subcategory of $\B$, and let $E : \B \to \A$ be an accessible functor. Then there exists a regular cardinal $\kappa$ such that the restriction $E|_{\B_0}$ is left Kan extended from $\B^\kappa_0$ (where $\B^\kappa_0$ is as in \cref{t--u--relatively-accessible}).
\end{lemma}

\begin{proof}
  Choose $\kappa$ so that $E$ is $\kappa$-accessible, $\B$ is $\kappa$-accessible, and $\B_0$ is relatively $\kappa$-accessible in $\B$. The inclusion $\B^\kappa_0 \inj \B$ extends uniquely to a fully faithful, $\kappa$-accessible functor $\Ind_\kappa(\B^\kappa_0) \inj \B$; let us abuse notation and regard $\Ind_\kappa(\B^\kappa_0)$ as a full subcategory of $\B$. Then we have a sequence of inclusions of full subcategories
  \[
    \B^\kappa_0 \inj \B_0 \inj \Ind_\kappa(\B^\kappa_0) \inj \B.
  \]
  Since $E$ is $\kappa$-accessible, its restriction to $\Ind_\kappa(\B^\kappa_0)$ is left Kan extended from $\B^\kappa_0$. It follows that its restriction to $\B_0$ is left Kan extended from $\B^\kappa_0$, as desired.
\end{proof}

\begin{proposition}
  \label{t--u--abstract}
  Let $\A$ and $\B$ be presentable $\infty$-categories, and let $\B_0$ be a relatively accessible full subcategory of $\B$. Then the following statements hold:
  \begin{enumerate}
  \item \label{t--u--abstract--extension}
    For any accessible functor $E : \B \to \A$, the restriction $E|_{\B_0}$ admits a left Kan extension $E' : \B \to \A$, and $E'$ is accessible.
  \item \label{t--u--abstract--initial}
    Suppose that $\A$ is pointed. Let $\Fun^\acc(\B,\A)$ be the $\infty$-category of accessible functors from $\B$ to $\A$ and let $\Fun^\acc_0(\B,\A)$ be the full subcategory thereof spanned by those accessible functors $E : \B \to \A$ such that $E|_{\B_0} \iso 0$. For $E \in \Fun^\acc(\B,\A)$, let $E' : \B \to \A$ denote the left Kan extension of \cref{t--u--abstract--extension}, $\psi : E' \to E$ the tautological natural transformation, and $\chi : E \to E''$ the cofiber of $\psi$. Then the inclusion functor $\Fun^\acc_0(\B,\A) \inj \Fun^\acc(\B,\A)$ admits a left adjoint $L$, given as follows: for $E \in \Fun^\acc(\B,\A)$, we have a natural identification $L(E) \iso E''$, under which the counit transformation $E \to L(E)$ identifies with $\chi$.
  \item \label{t--u--abstract--monoidal}
    Suppose that $\A$ is pointed and $\B$ is stable, that both are equipped with presentable symmetric monoidal structures, and that $\B_0$ is a thick tensor ideal in $\B$. Let $\Fun^{\acc,\laxdash{\otimes}}(\B,\A)$ be the $\infty$-category of accessible lax symmetric monoidal functors from $\B$ to $\A$ and let $\Fun^{\acc,\laxdash{\otimes}}_0(\B,\A)$ be the full subcategory thereof spanned by those accessible lax symmetric monoidal functors $E : \B \to \A$ such that $E|_{\B_0} \iso 0$. Then the inclusion functor $\Fun^{\acc,\laxdash{\otimes}}_0(\B,\A) \inj \Fun^{\acc,\laxdash{\otimes}}(\B,\A)$ admits a left adjoint $L^{\laxdash{\otimes}}$, and the diagram
    \[
      \begin{tikzcd}
        \Fun^{\acc,\laxdash{\otimes}}(\B,\A) \ar[r, "L^{\laxdash{\otimes}}"] \ar[d] &
        \Fun^{\acc,\laxdash{\otimes}}_0(\B,\A) \ar[d] \\
        \Fun^\acc(\B,\A) \ar[r, "L"] &
        \Fun^\acc_0(\B,\A)
      \end{tikzcd}
    \]
    commutes (the vertical arrows being the forgetful functors and $L$ being as in \cref{t--u--abstract--initial}) .
  \end{enumerate}
\end{proposition}

\begin{proof}
  Let $E : \B \to \A$ be an accessible functor. By \cref{t--u--relatively-accessible-meaning}, we may choose a regular cardinal $\kappa$ such that $E$ is $\kappa$-accessible, $\B$ is $\kappa$-accessible, and $E|_{\B_0}$ is left Kan extended from $\B^\kappa_0$. Then the left Kan extension $E' : \B \to \A$ of $E|_{\B^\kappa_0}$, which exists because $\B^\kappa_0$ is small and $\A$ is presentable, is also a left Kan extension of $E|_{\B_0}$. Moreover, $E'$ is left Kan extended from $\B^\kappa$, and hence is $\kappa$-accessible since $\B$ is $\kappa$-accessible. This proves \cref{t--u--abstract--extension}.

  We now move on to \cref{t--u--abstract--initial}, assuming that $\A$ is pointed. Note that the restriction of $\psi : E' \to E$ to $\B_0$ is an equivalence, and hence $E''|_{\B_0} \iso 0$. On the other hand, for any functor $F : \B \to \A$ such that $F|_{\B_0} \iso 0$, we have
  \[
    \Map_{\Fun(\B,\A)}(E',F) \iso \Map_{\Fun(\B_0,\A)}(E|_{\B_0},F|_{\B_0}) \iso *.
  \]
  These observations imply the claim.

  For \cref{t--u--abstract--monoidal}, we first give an alternative description of $E''$ under the additional assumptions that $\B$ is stable and $\B_0$ is a thick subcategory of $\B$. We begin with a digression: we have adjunctions of stable, presentable $\infty$-categories
  \[
    \begin{tikzcd}
      \Ind(\B^\kappa_0) \ar[r, "j", shift left=0.5ex] &
      \Ind(\B^\kappa) \ar[r, "q", shift left=0.5ex] \ar[l, "r", shift left=0.5ex] &
      \cat{Q}. \ar[l, "i", shift left=0.5ex]
    \end{tikzcd}
  \]
  In the left hand adjunction, $j$ is the $\Ind$-extension of the inclusion $\B^\kappa_0 \inj \B^\kappa$ and $r$ is its right adjoint. In the right hand adjunction, $\cat{Q}$ is the kernel of $r$, i.e. the full subcategory of $\Ind(\B^\kappa)$ spanned by those objects $B$ such that $\smash{\Map_{\Ind(\B^\kappa)}(B',B) \iso *}$ for all $B' \in \B^\kappa_0$; this is a localization of $\Ind(\B^\kappa)$, and $i$ denotes the inclusion functor and $q$ its left adjoint.\footnote{Another name for $\cat{Q}$ is $\Ind(\B^\kappa/\B^\kappa_0)$, where $\B^\kappa/\B^\kappa_0$ denotes the Verdier quotient.} The key point now is that the counit map for the left hand adjunction and the unit map for the right hand adjunction determine a cofiber sequence
  \[
    \smash{jr \lblto{\coun} \id_{\Ind(\B^\kappa)} \lblto{\un} iq.}
  \]
  To see this, let $B \in \Ind(\B^\kappa)$, and set $B' := jr(B)$ and $B'' := \cofib(\coun : B' \to B)$. We have $q(B') \iso 0$, implying that $q$ carries the canonical map $B \to B''$ to an isomorphism. It remains only to check that $B'' \in \cat{Q}$. For this we use our stability hypotheses, which imply that $r$ is an exact functor; since $r$ carries $\coun$ to an isomorphism, we deduce that $r(B'') \iso 0$, i.e. $B'' \in \cat{Q}$.

  Let us now return to our cofiber sequence of functors $E' \to E \to E''$. Each of these functors is $\kappa$-accessible, and hence left Kan extended from $\B^\kappa$; so let us focus on their restrictions to $\B^\kappa$. We may factor $E|_{\B^\kappa}$ as the composition
  \[
    \smash{\B^\kappa \lblto{e} \Ind(\B^\kappa) \lblto{\Ind(E)} \Ind(\A) \lblto{c} \A,}
  \]
  where $e$ is the tautological embedding and $c$ is the colimit functor, i.e. the $\Ind$-extension of $\id : \A \to \A$. Similarly, since $E'|_{\B^\kappa}$ is the left Kan extension of $E|_{\B^\kappa_0}$, it may be factored as the composition
  \[
    \smash{\B^\kappa \lblto{e} \Ind(\B^\kappa) \lblto{r} \Ind(\B^\kappa_0) \lblto{j} \Ind(\B^\kappa) \lblto{\Ind(E)} \Ind(\A) \lblto{c} \A,}
  \]
  (the composition $re$ is the left Kan extension of the tautological embedding $\B^\kappa_0 \inj \Ind(\B^\kappa_0)$). We then deduce from the cofiber sequence in the previous paragraph that $E''|_{\B^\kappa}$ factors as the composition
  \begin{equation}
    \label{t--u--abstract--composition}
    \smash{\B^\kappa \lblto{e} \Ind(\B^\kappa) \lblto{q} \cat{Q} \lblto{i} \Ind(\B^\kappa) \lblto{\Ind(E)} \Ind(\A) \lblto{c} \A,}
  \end{equation}
  with the map $\chi : E|_{\B^\kappa} \to E''|_{\B^\kappa}$ being induced by the unit map $\un : \id_{\Ind(\B^\kappa)} \to iq$.

  With that alternative description of $E''$ in hand, let us suppose finally that $\B$ and $\A$ are equipped with presentable symmetric monoidal structures and that $\B_0$ is a thick tensor ideal, and prove \cref{t--u--abstract--monoidal}. Enlarging $\kappa$ if necessary, we may assume that the tensor product functor $\otimes : \B \times \B \to \B$ preserves $\kappa$-compact objects, so that the symmetric monoidal structure of $\B$ restricts to the subcategory $\B^\kappa$ and $\B^\kappa_0$ is a tensor ideal in $\B^\kappa$. It follows that $\Ind(\B^\kappa_0)$ is a tensor ideal in $\Ind(\B^\kappa)$, hence that the localization $q : \Ind(\B^\kappa) \to \cat{Q}$ is compatible with the symmetric monoidal structure of $\Ind(\B^\kappa)$, and hence that the unit map $\un : \id_{\Ind(\B^\kappa)} \to iq$ is canonically a lax symmetric monoidal transformation. By the previous paragraph, this implies that the restriction to $\B^\kappa$ of the functor $E''$ and the natural transformation $\chi : E \to E''$ carry canonical lax symmetric monoidal structures. The same is then true before restriction, i.e. after left Kan extension along the inclusion $\B^\kappa \inj \B \iso \Ind_\kappa(\B^\kappa)$.

  To complete the argument, we must show that, equipped with this structure, the map $\chi : E \to E''$ exhibits $E''$ as the localization $L^{\laxdash{\otimes}}(E)$. Writing ``$\acck{\kappa}$'' in place of ``$\acc$'' to denote restriction to $\kappa$-accessible functors, note that our construction $E \mapsto E''$ defines a functor
  \[
    L^{\laxdash{\otimes}}_\kappa : \Fun^{\acck{\kappa},\laxdash{\otimes}}(\B,\A) \to \Fun^{\acck{\kappa},\laxdash{\otimes}}_0(\B,\A).
  \]
  As we are free to enlarge $\kappa$, it suffices to check that the map $\chi : E \to E''$ exhibits $L^{\laxdash{\otimes}}_\kappa$ as left adjoint to the inclusion functor. This is equivalent to the induced maps $L^{\laxdash{\otimes}}_\kappa(\chi) : L^{\laxdash{\otimes}}_\kappa(E) \to L^{\laxdash{\otimes}}_\kappa(E'')$ and $\chi : L^{\laxdash{\otimes}}_\kappa(E) \to L^{\laxdash{\otimes}}_\kappa(E)''$ being equivalences, which may be checked after forgetting lax symmetric monoidal structures, and hence follows from \cref{t--u--abstract--initial}.
\end{proof}

We now deduce \cref{t--u--main} from \cref{t--u--abstract}.

\begin{lemma}
  \label{t--u--nilpotent-relatively-accessible}
  Let $i : V \to T$ be a $\D$-$\sharp$-able map in $\S$. Then $\D(T)_{\nil{i}} \subseteq \D(T)$ is relatively accessible.
\end{lemma}

\begin{proof}
  Let $\kappa$ be a regular cardinal. Choose $\kappa' \ge \kappa$ such that $\D(V)$ is $\kappa'$-presentable. Let $\D(T)_0$ be the thick subcategory of $\D(T)$ generated by $i_\sharp(M)$ for $\smash{M \in \D(V)^{\kappa'}}$. Note that $i_\sharp$ preserves $\kappa'$-compactness, since its right adjoint $i^*$ preserves colimits, so that $\smash{\D(T)_0 \subseteq \D(T)^{\kappa'} \cap \D(T)_{\nil{i}}}$. Let $\D(T)_1$ be the full subcategory of $\D(T)_{\nil{i}}$ consisting of objects that can be written as $\kappa'$-filtered colimits in $\D(T)$ of objects in $\D(T)_0$. In other words, we have $\D(T)_1 = \D(T)_{\nil{i}} \cap \Ind_{\kappa'}(\D(T)_0)$, where we regard $\Ind_{\kappa'}(\D(T)_0)$ as a full subcategory of $\D(T)$ using the (fully faithful) $\Ind_{\kappa'}$-extension $\Ind_{\kappa'}(\D(T)_0) \inj \D(T)$ of the inclusion $\D(T)_0 \inj \D(T)$.

  To prove the claim, it suffices to show that $\D(T)_1 = \D(T)_{\nil{i}}$. Since $\D(V)$ is $\kappa'$-presentable, we know that $\D(T)_1$ contains $i_\sharp(M)$ for all $M \in \D(V)$. It remains to check that $\D(T)_1$ is thick. This follows from the fact that $\Ind_{\kappa'}(\D(T)_0)$ is thick.
\end{proof}

\begin{lemma}
  \label{t--u--conservative-monadic}
  Let $i : V \to T$ be a $\D$-$\sharp$-able map in $\S$. Then, if the functor $i^* : \D(T) \to \D(V)$ is conservative, it is monadic.
\end{lemma}

\begin{proof}
  This follows from the monadicity theorem, as $i^*$ preserves limits and colimits.
\end{proof}

\begin{lemma}
  \label{t--u--monadic-identity}
  Let $\B$ and $\C$ be $\infty$-categories, and let $\rho : \B \to \C$ be a monadic functor with left adjoint $\lambda : \C \to \B$. Let $\B' \subseteq \B$ be the essential image of $\lambda$. Then $\B'$ is dense in $\B$, i.e. the identity functor $\id : \B \to \B$ is left Kan extended from $\B'$.
\end{lemma}

\begin{proof}
  Since $\rho$ is monadic, the functor $\id : \B \to \B$ can be obtained as a geometric realization of the functors $(\lambda\rho)^n : \B \to \B$ for $n \ge 1$ \cite[Proof of Lemma 4.7.3.13]{lurie--HA}. As the property of being left Kan extended from $\B'$ is closed under colimits, it therefore suffices to show that this property is satisfied by $(\lambda\rho)^n$ for any $n \ge 1$. This follows from two facts: that any functor $\phi : \B \to \B$ obtained by left Kan extension along $\lambda : \C \to \B$ is left Kan extended from $\B'$, since $\lambda$ factors through $\B'$; and that left Kan extension along $\lambda$ is equivalent to precomposition with its right adjoint $\rho$.
\end{proof}

\begin{proof}[Proof of \cref{t--u--main}]
  Let $\gamma : f_* \to F$ be a natural transformation of functors $\D(T) \to \D(S)$ such that $F$ is accessible and $F|_{\D(T)_{\nil{i}}} \iso 0$. We will first show that, if $\fib(\gamma)$ preserves colimits, then $\gamma$ is initial among such natural transformations; this will prove \cref{t--u--main--initial} and \cref{t--u--main--fiber-criterion}, as $f_\t|_{\D(T)_{\nil{i}}} \iso 0$ by definition of $i$ being a $\D$-Tate cover of $f$, and $\fib(\can) \iso f_\sharp(- \otimes \poincaretwist_f)$ preserves colimits.

  By \cref{t--u--nilpotent-relatively-accessible}, we may apply \cref{t--u--abstract} apply here. \itemref{t--u--abstract}{initial} shows that we need only check that $\fib(\gamma)$ is left Kan extended from $\D(T)_{\nil{i}}$. Keeping in mind \cref{t--u--relatively-accessible-meaning}, this holds because $\fib(\gamma)$ preserves colimits and, by definition of $i$ being a $\D$-Tate cover of $f$ and \cref{t--u--conservative-monadic,t--u--monadic-identity}, the identity functor $\id_{\D(T)}$ is left Kan extended from $\D(T)_{\nil{i}}$.

  Finally, \itemref{t--u--abstract}{monoidal} gives statement \cref{t--u--main--monoidal} here, thanks to \cref{t--u--nilpotent-ideal}.
\end{proof}

We close this subsection with the following criterion for checking the last condition of \cref{t--u--cover}; it may be regarded as a generalization of a result of Klein concerning compact Lie groups \cite[Corollary 10.2]{klein--dualizing}.

\begin{proposition}
  \label{t--u--nilpotent-tate-vanishing}
  Let $f : T \to S$ be a $\D$-Tate map in $\S$, let $i : V \to T$ be a $\D$-Poincar\'e map in $\S$, and suppose that $fi : V \to S$ is $\D$-prim and $\D$-suave. Then $f_\t \circ i_\sharp : \D(V) \to \D(S)$ vanishes.
\end{proposition}

\begin{proof}
  We need to show that the norm map $\Nm_f : f_\sharp(i_\sharp(M) \otimes \poincaretwist_f) \to f_*(i_\sharp(M))$ is an equivalence for every $M \in \D(V)$. Since $\poincaretwist_i$ is invertible, we may replace $M$ by $M \otimes \poincaretwist_i$. \cref{t--c--poincare-comp} gives us a commutative diagram
  \[
    \begin{tikzcd}[row sep=small]
      f_\sharp i_\sharp (M \otimes \poincaretwist_i \otimes i^*(\poincaretwist_f)) \ar[r, "\pr^\sharp_i"] \ar[dd, "\sim", swap] &
      f_\sharp(i_\sharp(M \otimes \poincaretwist_i) \otimes \poincaretwist_f) \ar[d, "\Nm_f", swap] &
      \\
      &
      f_*(i_\sharp(M \otimes \poincaretwist_i)) \ar[r, "\Nm_i"] &
      f_*i_* \ar[d, "\sim"] \\
      (fi)_\sharp(M \otimes \poincaretwist_i \otimes i^*(\poincaretwist_f)) \ar[r, "\comp^\poincaretwist_{f,i}"]&
      (fi)_\sharp(M \otimes \poincaretwist_{fi}) \ar[r, "\Nm_{fi}"] &
      (fi)_*.
    \end{tikzcd}
  \]
  The projection map $\pr^\sharp_i$ is an isomorphism by \cref{t--p--pr-bc}, the map $\comp^\poincaretwist_{f,i}$ is an isomorphism by \cref{t--c--poincare-equiv}, the map $\Nm_i$ is an isomorphism because $i$ is Poincar\'e, and the map $\Nm_{fi}$ is an isomorphism because $fi$ is $\D$-prim and $\D$-suave. We deduce that the map $\Nm_f$ is an isomorphism too, as desired.
\end{proof}


\subsection{Functoriality of Tate cohomology}
\label{t--n}

In this subsection we again fix a commutative triangle
\[
  \begin{tikzcd}[row sep=small]
    U \ar[rr, "g"] \ar[dr, "h", swap] &
    &
    T \ar[dl, "f"] \\
    &
    S
  \end{tikzcd}
\]
in $\S$. Recall that cohomology has the following contravariant functoriality in this situation:

\begin{construction}
  \label{t--n--cohomology}
  Suppose that $f$ and $g$ are $\D$-$*$-able, and hence so is $h$. Let $\un^*_g : \id_{\D(T)} \to g_*g^*$ denote the unit transformation. We denote by $g^\zstar : f_* \to h_*g^*$ the natural transformation of functors from $\D(T)$ to $\D(S)$ given by the following composition:
  \[
    f_* \lblto{\un^*_g} f_*g_*g^* \iso  h_*g^*.
  \]
\end{construction}

The following construction shows that, under certain further assumptions, the above functoriality of cohomology extends compatibly to Tate cohomology. The subsequent result then explains how this interacts with the universal property of Tate cohomology discussed in the previous subsection.

\begin{construction}
  \label{t--n--tate}
  Suppose that $f$ is $\D$-Tate and that $g$ is $\D$-Poincar\'e. By \cref{t--c--poincare-equiv}, $h$ is $\D$-Tate and we have an isomorphism $\comp^\poincaretwist_{f,g} : \poincaretwist_g \otimes g^*(\poincaretwist_f) \isoto \poincaretwist_h$. We define the natural transformation $g^\zstar : f_\sharp(- \otimes \poincaretwist_f) \to h_\sharp(g^*(-) \otimes \poincaretwist_h)$ of functors from $\D(T)$ to $\D(S)$ to be the composition
  \begin{align*}
    f_\sharp(- \otimes \poincaretwist_f)
    \smash{{}\lblto{\un^*_g}{}}& f_\sharp(g_*g^*(-) \otimes \poincaretwist_f) \\[0.5ex]
    \smash{\lbliso{\Nm_g}{}}& f_\sharp(g_\sharp(g^*(-) \otimes \poincaretwist_g) \otimes \poincaretwist_f) \\[1ex]
    \smash{\lbliso{\pr^\sharp_g}{}}& f_\sharp g_\sharp(g^*(-) \otimes \poincaretwist_g \otimes g^*(\poincaretwist_f)) \\[0.5ex]
    \smash{\lblto{\comp^\poincaretwist_{f,g}}{}}& f_\sharp g_\sharp(g^*(-) \otimes \poincaretwist_h) \\
    \iso{}& h_\sharp(g^*(-) \otimes \poincaretwist_h);
  \end{align*}
  here $\eta^*_g$ is the unit map as in \cref{t--n--cohomology}, and the norm map $\Nm_g : g_\sharp(- \otimes \poincaretwist_g) \to g_*$ and projection map $\pr^\sharp_g$ are isomorphisms since $g$ is $\D$-Poincare (for the latter see \cref{t--p--pr-bc}). From the commutativity of the diagram in \cref{t--c--poincare-comp} we deduce the commutativity of the diagram
  \[
    \begin{tikzcd}
      f_\sharp(- \otimes \poincaretwist_f) \ar[r, "g^\zstar"] \ar[d, "\Nm_f", swap] &
      h_\sharp(g^*(-) \otimes \poincaretwist_h) \ar[d, "\Nm_h"] \\
      f_* \ar[r, "g^\zstar"] &
      h_*g^*.
    \end{tikzcd}
  \]
  Assuming that $\D(S)$ is pointed and admits cofibers, we may taking vertical cofibers to define a natural transformation $g^\zstar : f_\t \to h_\t g^*$ of functors from $\D(T)$ to $\D(S)$.
\end{construction}

\begin{theorem}
  \label{t--n--tate-uniqueness-monoidality}
  Suppose as in \cref{t--u} that $\D$ is stable and presentable, and as in \cref{t--n--tate} that $f$ is $\D$-Tate and $g$ is $\D$-Poincar\'e. Suppose also given a $\D$-Tate cover $i : V \to T$ of $f$ (\cref{t--u--cover}). Then the following statements hold:
  \begin{enumerate}
  \item \label{t--n--tate--map}
    The natural transformation $g^\zstar : f_\t \to h_\t g^*$ of \cref{t--n--tate} is the unique such making the following diagram of functors from $\D(T)$ to $\D(S)$ commute:
    \[
      \begin{tikzcd}
        f_* \ar[r, "g^\zstar"] \ar[d, "\can", swap] &
        h_*g^* \ar[d, "\can"] \\
        f_\t \ar[r, "g^\zstar"] &
        h_\t g^*.
      \end{tikzcd}
    \]

  \item \label{t--n--tate--monoidal}
    Let us regard the functor $f_\t$ and the map $\can : f_* \to f_\t$ as equipped with the lax symmetric monoidal structures of \cref{t--u--main}, and suppose also given lax symmetric monoidal structures on the functor $h_\t : \D(U) \to \D(S)$ and the map $\can : h_* \to h_\t$. Then there is a unique lax symmetric monoidal structure on the map $g^\zstar : f_\t \to h_t g^*$ making the commutative diagram in \cref{t--n--tate--map} one of lax symmetric monoidal functors.
  \end{enumerate}
\end{theorem}

\begin{proof}
  By \cref{t--u--main}, it suffices to show that the composition $h_\t g^* i_\sharp : \D(V) \to \D(S)$ vanishes. In other words, we have to show that, for $M \in \D(V)$, the norm map
  \[
    \Nm_h : h_\sharp(g^*i_\sharp(M) \otimes \poincaretwist_h) \to h_*g^*i_\sharp(M)
  \]
  is an isomorphism. By the commutative diagram in \cref{t--c--poincare-comp}, this map identifies with the composition
  \begin{align*}
    f_\sharp g_\sharp(g^*i_\sharp(M) \otimes \poincaretwist_g \otimes g^*(\poincaretwist_f))
    \smash{{}\lblto{\pr^\sharp_g}{}}& f_\sharp (g_\sharp(g^*i_\sharp(M) \otimes \poincaretwist_g) \otimes \poincaretwist_f) \\[0.5ex]
    \smash{\lblto{\Nm_f}{}}& f_* g_\sharp(g^*i_\sharp(M) \otimes \poincaretwist_g) \\[0.5ex]
    \smash{\lblto{\Nm_g}{}}& f_*g_*(g^*i_\sharp(M)).
  \end{align*}
  The projection map $\pr^\sharp_g$ and the norm map $\Nm_g$ are isomorphisms because $g$ is $\D$-Poincar\'e (for the former see \cref{t--p--pr-bc}). The norm map $\Nm_f$ is an isomorphism in this case because $i$ is a $\D$-Tate cover of $f$ and we have
  \[
    \smash{g_\sharp(g^*i_\sharp(M) \otimes \poincaretwist_g) \lbliso{\pr^\sharp_g} i_\sharp(M) \otimes g_\sharp(\poincaretwist_g) \lbliso{\pr^\sharp_i} i_\sharp(M \otimes i^*g_\sharp(\poincaretwist_g)) \in \D(T)_{\nil{i}}}
  \]
  (the $\D$-$\sharp$-projection formula holds also for $i$ by definition of $\D$-Tate cover).
\end{proof}


\section{Examples}
\label{e}

This section illustrates how the generalities of \cref{t} play out in more specific examples of three functor formalisms.

In \cref{e--l}, we will discuss local systems on spaces, rederiving what was overviewed in \cref{i--bg}, and generalizing by allowing the local systems to be valued in any stable, presentable symmetric monoidal $\infty$-category, rather than just the $\infty$-category of spectra. The generalization will allow us to formulate a new example of functoriality in Tate cohomology for the circle group, involving the ambidexterity phenomenon of Hopkins--Lurie \cite{hopkins-lurie--ambi}.

In \cref{e--q}, we will discuss quasicoherent sheaves on algebraic stacks, for a rather loose/general interpretation of this notion: we will take the affine objects to correspond to commutative algebra objects in an arbitrary stable, presentable symmetric monoidal $\infty$-category. Our main aim will be to present in a more systematic way the construction from \cite[\textsection 2.4]{raksit--hhdr} of Tate cohomology for modules over suitable bialgebras. The account here will also recover and generalize some of Rognes's theory of ``stably dualizable groups'' \cite{rognes--groups}.

\begin{remark}
  \label{e--absolute}
  In both contexts to be discussed in this section, the $\infty$-category $\S$ on which our three functor formalism $\D$ is defined (namely the $\infty$-category of spaces in \cref{e--l}, and the $\infty$-category of stacks in \cref{e--q}) will have a final object $*$. We will allow ourselves to abbreviate by identifying an object $T \in \S$ with the unique map $f : T \to *$. For instance, we will say that $T$ is $\D$-Poincar\'e if $f$ is. In addition, as in \cref{i--bg}, we will use the notation
  \[
    \Ch_*(T;-) := f_\sharp : \D(T) \to \D(*), \quad
    \Ch^*(T;-) := f_* : \D(T) \to \D(*), \quad
    \Ct^*(T;-) := f_\t : \D(T) \to \D(*)
  \]
  for the associated homology, cohomology, and Tate cohomology functors (when these are defined).
\end{remark}


\subsection{Local systems}
\label{e--l}

In this subsection, our three functor formalism will come from the theory of $\C$-valued local systems on spaces, where $\C$ is a fixed stable, presentable symmetric monoidal $\infty$-category. We begin with the functor
\begin{equation}
  \label{e--l--loc}
  \Loc_\C := \C^{(-)} : \Spc^\op \to \CAlg(\PrLst),
\end{equation}
carrying a space $S$ to the stable, presentable symmetric monoidal $\infty$-category $\Loc_\C(S) = \C^S$ of $\C$-valued local systems on $S$, i.e. functors from $S$ to $\C$ (the symmetric monoidal structure being the pointwise one), and carrying a map of spaces $f : T \to S$ to the restriction functor $f^* : \Loc_\C(S) \to \Loc_\C(T)$. We may extend $\Loc_\C$ to a three functor formalism as follows.

\begin{proposition}
  \label{e--l--hoco}
  Let $f : T \to S$ be a map of spaces. Then:
  \begin{enumerate}
  \item \label{e--l--hoco--hoco}
    the functor $f^* : \Loc_\C(S) \to \Loc_\C(T)$ admits right and left adjoints, $f_* : \Loc_\C(T) \to \Loc_\C(S)$ and $f_\sharp : \Loc_\C(T) \to \Loc_\C(S)$.
  \item \label{e--l--hoco--pr}
    the projection map $\pr^\sharp_f : f_\sharp(M \otimes f^*(L)) \to f_\sharp(M) \otimes L$ is an isomorphism for all $L \in \Loc_\C(S)$ and $M \in \Loc_\C(T)$;
  \item \label{e--l--hoco--bc}
    for any cartesian square of spaces $\sigma :$
    \[
      \begin{tikzcd}
        V \ar[r, "q"] \ar[d, "g", swap] &
        T \ar[d, "f"] \\
        U \ar[r, "p"] &
        S,
      \end{tikzcd}
    \]
    the base change maps $\bc^*_\sigma : p^*f_* \to g_*q^*$ and $\bc^\sharp_\sigma : g_\sharp q^* \to p^* f_\sharp$ of functors from $\Loc_\C(T)$ to $\Loc_\C(U)$ are isomorphisms.
  \end{enumerate}
\end{proposition}

\begin{proof}
  \cref{e--l--hoco--hoco} follows from the fact that $\C$ is presentable; the adjoints $f_*$ and $f_\sharp$ are given by right and left Kan extension along $f$, respectively. For \cref{e--l--hoco--bc}, see e.g. \cite[Proposition 4.3.3]{hopkins-lurie--ambi}. For \cref{e--l--hoco--pr}, it suffices to show for any point $s : * \to S$ that $\pr^\sharp_f$ becomes an equivalence after applying $s^*$. Applying \cref{e--l--hoco--bc}, this reduces us to the case that $S = *$. Then the projection map is the canonical map
  \[
    \colim_{t \in T} {(M(t) \otimes L)} \to (\colim_{t \in T} M(t)) \otimes L
  \]
  in $\C$, which again by presentability is an isomorphism.
\end{proof}

\begin{corollary}
  \label{e--l--D}
  The functor \cref{e--l--loc} extends canonically to a (stable, presentable) three functor formalism $\Loc_\C : \Span(\Spc) \to \PrLst$, where:
  \begin{enumerate}
  \item \label{e--l--D--able}
    every map of spaces is $!$-able, $*$-able, and $\sharp$-able, and satisfies the $\sharp$-projection formula and both $*$-- and $\sharp$--base change;
  \item \label{e--l--D--push}
    for any map of spaces $f : T \to S$, there is a canonical identification of functors $f_! \iso f_\sharp$, as well as a canonical identification of projection isomorphisms $\pr^!_f \iso \pr^\sharp_f$;
  \item \label{e--l--D--bc}
    for any cartesian square of spaces $\sigma :$
    \[
      \begin{tikzcd}
        V \ar[r, "q"] \ar[d, "g", swap] &
        T \ar[d, "f"] \\
        U \ar[r, "p"] &
        S,
      \end{tikzcd}
    \]
    there is a canonical identification of base change isomorphisms $\bc^!_\sigma \iso \bc^\sharp_\sigma$.
  \end{enumerate}
\end{corollary}

\begin{proof}
  See for instance \cite[Proposition 3.3.3]{heyer-mann--six} (or the work of Cnossen--Lenz--Linskens \cite{cll--six} for another approach to the same result). In terms of the notation there, we take in this case the class $I$ to consist of all maps of spaces and the class $P$ to consist only of isomorphisms.
\end{proof}

For the remainder of the subsection, we will work in the context of the three functor formalism $\Loc_\C$ of \cref{e--l--D}. To get going, we note the following.

\begin{proposition}
  \label{e--l--suave}
  Let $f : T \to S$ be a map of spaces. Then $f$ is $\Loc_\C$-suave, and there is a canonical isomorphism $\suavetwist_f \iso \unit_T$ in $\Loc_\C(T)$.
\end{proposition}

\begin{proof}
  This follows from the construction of $\Loc_\C$ as a three functor formalism: unravelling definitions shows that the $\sharp$-norm map associated to $f$ recovers the identification $f_\sharp \iso f_!$ of \itemref{e--l--D}{push}.
\end{proof}

It follows from \cref{e--l--suave} that every map of spaces $f : T \to S$ is $\Loc_\C$-Tate, and hence by \cref{t--d--norm} has an associated norm map $\Nm_f : f_\sharp(- \otimes \poincaretwist_f) \to f_*$ and Tate cohomology functor $f_\t = \cofib(\Nm_f) : \Loc_\C(T) \to \Loc_\C(S)$. Note that, in this setting, the norm map $\Nm_f$ is identical to the $*$-norm map $\Nm^*_f : f_!(- \otimes \primtwist_f) \to f_*$. In particular, a map of spaces is $\Loc_\C$-Poincar\'e if and only if it is $\Loc_\C$-terse, and hence every $\Loc_\C$-proper map is $\Loc_\C$-Poincar\'e; regarding these notions, we have the following basic results.

\begin{proposition}
  \label{e--l--prim-fiber}
  Let $f : T \to S$ be a map of spaces. Then the following conditions are equivalent:
  \begin{enumerate}
  \item \label{e--l--prim-fiber--prim}
    $f$ is $\Loc_\C$-prim (resp. $\Loc_\C$-Poincar\'e, $\Loc_\C$-proper);
  \item \label{e--l--prim-fiber--fiber}
    for every point $s : * \to S$, the fiber $T_s$, i.e. the base change of $f$ along $s$, is $\Loc_\C$-prim (resp. $\Loc_\C$-Poincar\'e, $\Loc_\C$-proper).
  \end{enumerate}
\end{proposition}

\begin{proof}
  That \cref{e--l--prim-fiber--prim} implies \cref{e--l--prim-fiber--fiber} follows from the closure of primness and Poincar\'eness under base change (\cref{t--b--closure}). The converse follows from the fact that these conditions are local with respect to $\Loc_\C$-$*$-covers of $S$ (\cref{t--b--descent}), as $\Loc_\C$ is a sheaf for the canonical topology on $\Spc$ (i.e. carries colimits of spaces to limits of $\infty$-categories) and the collection of points of $S$ generates a cover for the canonical topology.
\end{proof}

\begin{lemma}
  \label{e--l--kunneth}
  For any spaces $T_0$ and $T_1$, the canonical map
  \[
    \Loc_\C(T_0) \otimes_\C \Loc_\C(T_1) \to \Loc_\C(T_0 \times T_1)
  \]
  in $\CAlg(\PrL)$, induced by the pullback functors along the projection maps $p_i : T_0 \times T_1 \to T_i$, is an equivalence.
\end{lemma}

\begin{proof}
  The key point is that the tensor product $- \otimes_\C -$ preserves limits of space-indexed diagrams in each variable: this is a consequence of the fact that such a limit in $\PrL$ is naturally equivalent to the colimit of the same diagram \cite[Example 4.3.11]{hopkins-lurie--ambi}, together with the fact that the tensor product preserves colimits in each variable. It follows that both the source and target of the map in question carry colimits in each $T_i$ to limits of $\infty$-categories, which allows us to reduce to the case that $T_0 = T_1 = *$, where the claim is clear.
\end{proof}

\begin{proposition}
  \label{e--l--prim-absolute}
  Let $T$ be a space. Then the norm map $\Nm_T : \Ch_*(T;- \otimes \poincaretwist_T) \to \Ch^*(T;-)$ is a $\C$-linear assembly map: that is, it exhibits the source as the terminal colimit preserving $\C$-linear functor equipped with a lax $\C$-linear map to the target. In particular, $T$ is $\Loc_\C$-prim if and only if the cohomology functor $\Ch^*(T;-) : \Loc_\C(T) \to \C$ is $\C$-linear (i.e. satisfies the projection formula) and preserves colimits.
\end{proposition}

\begin{proof}
  This follows from \cref{e--l--kunneth,t--d--norm-assembly,t--p--assembly}.
\end{proof}

\begin{remark}
  It follows from \cref{e--l--prim-absolute} that, in the case $\C = \Spt$, our norm map $\Nm_T$ agrees with that of Klein and Nikolaus--Scholze: see \cite[Theorem B]{klein--axioms} and \cite[Theorem I.4.1]{nikolaus-scholze--TC}.
\end{remark}

\begin{remark}
  \label{e--l--Cmap}
  Suppose given another stable, presentable symmetric monoidal $\infty$-category $\C'$ and a colimit preserving symmetric monoidal functor $\alpha : \C \to \C'$. Then there is an induced map of three functor formalisms $\alpha : \Loc_\C \to \Loc_{\C'}$. Thus, by \cref{t--p--Dprime}, for $f : T \to S$ a map of spaces, if $f$ is $\Loc_\C$-prim (resp. $\Loc_\C$-Poincar\'e), then it is also $\Loc_{\C'}$-prim (resp. $\Loc_\C'$-Poincar\'e).
\end{remark}

We can now point out the most basic class of $\Loc_\C$-Poincar\'e objects, upon which our terminology is based.

\begin{example}
  \label{e--l--poincare}
  Let $T$ be a space. We claim that if $T$ is compact (resp. Poincar\'e, in the sense defined in \cref{i--bg}), then it is $\Loc_\C$-prim (resp. $\Loc_\C$-Poincar\'e). Given that there is a (unique) colimit preserving symmetric monoidal functor $\Spt \to \C$, and given \cref{e--l--Cmap}, it's enough to verify this claim in the case $\C = \Spt$. Moreover, we need only check that compactness implies $\Loc_\Spt$-primness (that Poincar\'eness implies $\Loc_\Spt$-Poincar\'eness is then immediate from the definitions); this is precisely the assertion made in \cref{i--bg} that, when $T$ is compact, the norm map $\Nm_T$ is an isomorphism.

  To prove this, it suffices by \cref{e--l--prim-absolute} to show that, for $T$ compact (i.e. finitely dominated), the cohomology functor $\Ch^*(T;-) : \Loc_T(\Spt) \to \Spt$ preserves colimits. This property is stable under retracts in $T$, so we may assume that $T$ is finite, in which case it holds because finite limits commute with colimits in $\Spt$.
\end{example}

From this we may recover the following motivating example from \cref{i--bg} of the results of \cref{t--h} and \cref{t--u}. (See \cref{e--q--poincare-hopf,e--q--poincare-hopf-B} in the next subsection for a generalization.)

\begin{example}
  \label{e--l--poincare-group}
  Let $G$ be a compact group space, and let $\clspc{G}$ be its classifying space. Let $f : \clspc{G} \to *$ and $g : G \to *$ denote the unique such maps and let $i : * \to \clspc{G}$ denote the canonical basepoint. By \cref{e--l--poincare}, $G$ is $\Loc_\C$-prim. By \cref{t--h--main} it follows that $G$ is $\Loc_\C$-Poincar\'e. In view of the cartesian square
  \[
    \begin{tikzcd}
      G \ar[r, "g"] \ar[d, "g", swap] &
      * \ar[d, "i"] \\
      * \ar[r, "i"] &
      \clspc{G}
    \end{tikzcd}
  \]
  and \cref{e--l--prim-fiber}, $i$ is then $\Loc_\C$-Poincar\'e. Applying \cref{t--u--nilpotent-tate-vanishing} to the composition
  \[
    \begin{tikzcd}
      * \ar[rr, "i"] \ar[dr, "\id", swap] &
      &
      \clspc{G} \ar[dl, "f"] \\
      &
      *,
    \end{tikzcd}
  \]
  we deduce that $i$ is a $\Loc_\C$-Tate cover of $f$. Hence, \cref{t--u--main} applies to $f$: in terms of the traditional notation
  \[
    (-)^{\h G} := f_* : \Loc_\C(\clspc{G}) \to \C,
    \qquad
    (-)^{\t G} := f_\t : \Loc_\C(\clspc{G}) \to \C,
  \]
  we get that the Tate cohomology functor $(-)^{\t G}$ and the canonical map $(-)^{\h G} \to (-)^{\t G}$ admit canonical lax symmetric monoidal structures.
\end{example}

We next consider the $\Loc_\C$-proper maps (which, as noted above, are in particular $\Loc_\C$-Poincar\'e). By definition, these recover the notion of ``ambidexterity'' introduced by Hopkins--Lurie \cite{hopkins-lurie--ambi}:

\begin{definition}
  \label{e--l--ambidextrous}
  Let $f : T \to S$ be a map of spaces. We say that $f$ is \emph{$\C$-ambidextrous} if it is $\Loc_\C$-proper.
\end{definition}

\begin{remark}
  \label{e--l--ambidextrous-fiber}
  Following the convention of \cref{e--absolute}, we say that a space $T$ is $\C$-ambidextrous if the unique map $T \to *$ is $\C$-ambidextrous. By \cref{e--l--prim-fiber}, a map of spaces $T \to S$ is $\C$-ambidextrous if and only if each of its fibers is $\C$-ambidextrous.
\end{remark}

\begin{definition}
  \label{e--l--pi-finite}
  We say that a space $S$ is \emph{$\pi$-finite} if it satisfies the following conditions:
  \begin{enumerate}
  \item the set $\pi_0(S)$ is finite;
  \item for each $s \in S$ and each integer $j \ge 1$, the group $\pi_j(S,s)$ is finite;
  \item $S$ is truncated, i.e. there exists an integer $k$ such that for each $s \in S$ and each integer $j \ge k$, the group $\pi_j(S,s)$ vanishes.
  \end{enumerate}
\end{definition}

\begin{definition}
  \label{e--l--semiadditive}
  For any integer $k \ge -2$, we say that $\C$ is \emph{$k$-semiadditive} if every $k$-truncated, $\pi$-finite space $S$ is $\C$-ambidextrous. We say that $\C$ is $\infty$-semiadditive if it is $k$-semiadditive for all $k \ge -2$.
\end{definition}

\begin{example}
  \label{e--l--0semiadditive}
  The stability of $\C$ guarantees that it is $0$-semiadditive: that is, every finite set is $\C$-ambidextrous (see \cite[Proposition 4.4.9]{hopkins-lurie--ambi}).
\end{example}

\begin{example}
  \label{e--l--Kn-local}
  Let $n$ be a nonnegative integer. By work of Carmeli--Schlank--Yanovski \cite{csy--ambi}, the $\infty$-category $\C = \Spt_{\mathrm{T}(n)}$ of $\mathrm{T}(n)$-local spectra is $\infty$-semiadditive. This implies that the $\infty$-category $\C = \Spt_{\mathrm{K}(n)}$ of $\mathrm{K}(n)$-local spectra is $\infty$-semiadditive, proved earlier by Hopkins--Lurie \cite{hopkins-lurie--ambi}. In particular, the classifying spaces of finite groups are $\C$-ambidextrous for these $\infty$-categories $\C$, a result from even earlier work of Hovey--Sadofsky \cite{hovey-sadofsky--blueshift} and Kuhn \cite{kuhn--blueshift}.
\end{example}

From \cref{e--l--0semiadditive} and the results of \cref{t--n} we obtain the classical functoriality of Tate cohomology for finite groups.

\begin{example}
  \label{e--l--fingrp}
  Let $\gamma : H \to G$ be an injective map of finite groups, and let $g : \clspc{H} \to \clspc{G}$ be the induced map on classifying spaces. Let $i : * \to \clspc{G}$ and $j : * \to \clspc{H}$ be the canonical basepoints. By \cref{e--l--poincare-group}, $i$ and $j$ are $\Loc_\C$-Tate covers of $\clspc{G}$ and $\clspc{H}$, respectively, so that \cref{t--u--main} endows the Tate cohomology functors $(-)^{\t G} : \Loc_\C(\clspc{G}) \to \C$ and $(-)^{\t H} : \Loc_\C(\clspc{H}) \to \C$, as well as the canonical maps $(-)^{\h G} \to (-)^{\t G}$ and $(-)^{\h H} \to (-)^{\t H}$, with lax symmetric monoidal structures.

  Since $\gamma$ is injective, $g$ is a finite covering map: its fiber is the finite set $G/H$. It thus follows from \cref{e--l--0semiadditive} that $g$ is $\C$-ambidextrous, so that we may apply \cref{t--n--tate,t--n--tate-uniqueness-monoidality} to the commutative diagram
  \[
    \begin{tikzcd}[row sep=small]
      \clspc{H} \ar[rr, "g"] \ar[dr] &
      &
      \clspc{G} \ar[dl] \\
      &
      *,
    \end{tikzcd}
  \]
  giving us a commutative diagram 
  \[
    \begin{tikzcd}
      (-)^{\h G} \ar[r, "g^\zstar"] \ar[d, "\can", swap] &
      (g^*(-))^{\h H} \ar[d, "\can"] \\
      (-)^{\t G} \ar[r, "g^\zstar"] &
      (g^*(-))^{\t H}
    \end{tikzcd}
  \]
  of lax symmetric monoidal functors from $\Loc_\C(\clspc{G})$ to $\C$.
\end{example}

Our work also affords the following more exotic example of functoriality of Tate cohomology, which is applied in \cite[\textsection 2.2]{devalapurkar-r--THHZ}, in the case $\C = \Spt_{\mathrm{K}(1)}$ (allowable by \cref{e--l--Kn-local}), as a component of an analysis of the topological cyclic homology of $\ZZ$.

\begin{example}
  \label{e--l--TCZ}
  Let $\cir$ denote the circle group, and $\clspc{\cir}$ its classifying space. Let $i : * \to \clspc{\cir}$ denote the canonical basepoint. By \cref{e--l--poincare-group}, $i$ is a $\Loc_\C$-Tate cover of $\clspc{\cir}$, so that \cref{t--u--main} endows the Tate cohomology functor $(-)^{\t\cir} : \Loc_\C(\clspc{\cir}) \to \C$ and the canonical map $(-)^{\h\cir} \to (-)^{\t\cir}$ with lax symmetric monoidal structures.

  Now suppose that $\C$ is $1$-semiadditive, let $n$ be an integer, and let $[n] : \clspc{\cir} \to \clspc{\cir}$ be the map given by multiplication by $n$ (equivalently, the map induced by the $n$-fold covering map from $\cir$ to itself). The fiber of $[n]$ is $\clspc{\ZZ/n}$, and hence is $\C$-ambidextrous by hypothesis. Thus, the map $[n]$ is $\C$-ambidextrous, and in particular $\Loc_\C$-Poincar\'e.  Consequently, we may apply \cref{t--n--tate,t--n--tate-uniqueness-monoidality} to the commutative diagram
  \[
    \begin{tikzcd}[row sep=small]
      \clspc{\cir} \ar[rr, "{[n]}"] \ar[dr] &
      &
      \clspc{\cir} \ar[dl] \\
      &
      *,
    \end{tikzcd}
  \]
  giving us a commutative diagram 
  \[
    \begin{tikzcd}
      (-)^{\h\cir} \ar[r, "{[n]^\zstar}"] \ar[d, "\can", swap] &
      ([n]^*(-))^{\h\cir} \ar[d, "\can"] \\
      (-)^{\t\cir} \ar[r, "{[n]^\zstar}"] &
      ([n]^*(-))^{\t\cir}
    \end{tikzcd}
  \]
  of lax symmetric monoidal functors from $\Loc_\C(\clspc{\cir})$ to $\C$.
\end{example}


\subsection{Quasicoherent sheaves}
\label{e--q}

We fix for this subsection, as in the previous subsection, a stable, presentable symmetric monoidal $\infty$-category $\C$. Here we will extract a different three functor formalism from $\C$. We set $\Aff_\C := \CAlg(\C)^\op$. For $A \in \CAlg(\C)$, we will write $\Spec(A)$ for the corresponding object in $\Aff_\C$. And we let
\begin{equation}
  \label{e--q--qcoh}
  \QCoh : \Aff_\C^\op \to \CAlg(\PrLst)
\end{equation}
denote the functor corresponding to $\Mod : \CAlg(\C) \to \CAlg(\PrLst)$: that is, $\QCoh$ carries $\Spec(A) \in \Aff_\C$ to the stable, presentable symmetric monoidal $\infty$-category $\QCoh(\Spec(A)) = \Mod_A(\C)$ of $A$-modules in $\C$, and carries a map $f : \Spec(B) \to \Spec(A)$, corresponding to a map $A \to B$ of commutative algebras in $\C$, to the base change functor $f^* : \QCoh(\Spec(A)) \to \QCoh(\Spec(B))$, i.e. the functor $B \otimes_A - : \Mod_A(\C) \to \Mod_B(\C)$ (see \cite[\textsection 4.8.3--4.8.5]{lurie--HA} for the full construction). We will extend $\QCoh$ to a three functor formalism on $\Aff_\C$ in a manner dual to what was done in the previous subsection for local systems on spaces (the thought that this choice might serve our purposes was motivated by the definition of the six functor formalism on analytic stacks by Clausen--Scholze \cite{clausen-scholze--anstk}; see also Scholze's discussion of quasicoherent sheaves in \cite[Lecture VIII]{scholze--six}).

\begin{proposition}
  \label{e--l--star}
  Let $f : T \to S$ be a map in $\Aff_\C$. Then:
  \begin{enumerate}
  \item \label{e--l--star--star}
    the functor $f^* : \QCoh(S) \to \QCoh(T)$ admits a right adjoint $f_* : \QCoh(T) \to \QCoh(S)$, which preserves colimits;
  \item \label{e--l--star--pr}
    the projection map $\pr^*_f : f_*(M) \otimes L \to f_*(M \otimes f^*(L))$ is an isomorphism for all $L \in \QCoh(S)$ and $M \in \QCoh(T)$;
  \item \label{e--l--star--bc}
    for any cartesian square $\sigma :$
    \[
      \begin{tikzcd}
        V \ar[r, "q"] \ar[d, "g", swap] &
        T \ar[d, "f"] \\
        U \ar[r, "p"] &
        S,
      \end{tikzcd}
    \]
    in $\Aff_\C$, the base change map $\bc^*_\sigma : p^*f_* \to g_*q^*$ of functors from $\QCoh(T) \to \QCoh(U)$ is an isomorphism.
  \end{enumerate}
\end{proposition}

\begin{proof}
  Write $S = \Spec(A)$ and $T = \Spec(B)$, so that $f$ corresponds to a map $A \to B$ in $\CAlg(\C)$. Then $f^*$ is given by the functor $B \otimes_A - : \Mod_A(\C) \to \Mod_B(\C)$, whose right adjoint is the restriction functor $\Mod_B(\C) \to \Mod_A(\C)$; this proves \cref{e--l--star--star}. Claim \cref{e--l--star--pr} is the assertion that, for $L \in \Mod_A(\C)$ and $M \in \Mod_B(\C)$, the canonical map
  \[
    M \otimes_A L \to M \otimes_B (B \otimes_A L)
  \]
  is an isomorphism, which follows from associativity of the relative tensor product \cite[Propositon 4.4.3.13]{lurie--HA}. For \cref{e--l--star--bc}, let
  \[
    \begin{tikzcd}
      D \ &
      B \ar[l] \\
      C \ar[u] &
      A \ar[l] \ar[u],
    \end{tikzcd}
  \]
  be the cocartesian square in $\CAlg(\C)$ corresponding to $\sigma$. The assertion is then that, for $M \in \Mod_B(\C)$, the canonical map
  \[
    C \otimes_A M \to D \otimes_B M
  \]
  is an isomorphism. Noting that cocartesianness of the above square means that we have $D \iso C \otimes_A B$ (by \cite[Proposition 3.2.4.7 and Corollary 3.4.1.7]{lurie--HA}), this again follows from associativity of relative tensor product.
\end{proof}

\begin{corollary}
  \label{e--q--D}
  The functor \cref{e--q--qcoh} extends canonically to a (stable, presentable) three functor formalism $\QCoh : \Span(\Aff_\C) \to \PrLst$, where:
  \begin{enumerate}
  \item \label{e--q--D--able}
    every map in $\Aff_\C$ is $!$-able and $*$-able, and satisfies the $*$-projection formula and $*$--base change;
  \item \label{e--q--D--push}
    for any map $f : T \to S$ in $\Aff_\C$, there is a canonical identification of functors $f_! \iso f_*$, as well as a canonical identification of projection isomorphisms $\pr^!_f \iso \pr^*_f$;
  \item \label{e--q--D--bc}
    for any cartesian square $\sigma :$
    \[
      \begin{tikzcd}
        V \ar[r, "q"] \ar[d, "g", swap] &
        T \ar[d, "f"] \\
        U \ar[r, "p"] &
        S,
      \end{tikzcd}
    \]
    in $\Aff_\C$, there is a canonical identification of base change isomorphisms $\bc^!_\sigma \iso \bc^*_\sigma$.
  \end{enumerate}
\end{corollary}

\begin{proof}
  As with the proof of \cref{e--l--D}, we may apply \cite[Proposition 3.3.3]{heyer-mann--six}; but this time, we take the class $P$ to consist of all maps in $\Aff_\C$ and the class $I$ to consist only of isomorphisms.
\end{proof}

Working now in the context of the three functor formalism of \cref{e--q--D}, let us note first that primness and suaveness have simple characterizations.

\begin{proposition}
  \label{e--q--prim}
  Let $f : T \to S$ be a map in $\Aff_\C$. Then $f$ is $\QCoh$-prim, and there is a canonical isomorphism $\primtwist_f \iso \unit_T$ in $\QCoh(T)$.
\end{proposition}

\begin{proof}
  This is similar/dual to \cref{e--l--suave} from the previous subsection: here, the $*$-norm map associated to $f$ recovers the identification $f_! \iso f_*$ of \itemref{e--q--D}{push}.
\end{proof}

\begin{lemma}
  \label{e--q--kunneth}
  The functor \cref{e--q--qcoh} preserves pushouts: that is, given a cospan $T_0 \to S \from T_1$ in $\Aff_\C$, the induced map
  \[
    \QCoh(T_0) \otimes_{\QCoh(S)} \QCoh(T_1) \to \QCoh(T_0 \times_S T_1)
  \]
  in $\CAlg(\PrL)$ is an equivalence.
\end{lemma}

\begin{proof}
  Write $S = \Spec(A)$. We may replace $\C$ by $\Mod_A(\C)$, and thereby assume that $A$ is the unit object of $\C$ (or equivalently that $S$ is the final object of $\Aff_\C$). Then the claim follows from \cite[Remark 4.8.5.19]{lurie--HA}.
\end{proof}

\begin{proposition}
  \label{e--q--suave}
  Let $f : T \to S$ be a map in $\Aff_\C$, corresponding to a map $A \to B$ in $\CAlg(\C)$. Then the following conditions are equivalent:
  \begin{enumerate}
  \item \label{e--q--suave--suave}
    $f$ is $\QCoh$-suave;
  \item \label{e--q--suave--pr}
    $f$ is $\QCoh$-$\sharp$-able and satisfies the $\QCoh$-$\sharp$-projection formula;
  \item \label{e--q--suave--dual}
    $B$ is dualizable as an $A$-module.
  \end{enumerate}
\end{proposition}

\begin{proof}
  That \cref{e--q--suave--suave} implies \cref{e--q--suave--pr} was explained in \cref{t--p--pr-bc}, and that \cref{e--q--suave--pr} implies \cref{e--q--suave--dual} follows from \cref{t--h--dualizable}.

  Assume that \cref{e--q--suave--dual} holds. For any $B$-module $M$, we may consider the condition that the left adjoint $f_\sharp$ of $f^* = B \otimes_A (-)$ is defined at $M$: that is, that the functor
  \[
    \Mod_A(\C) \to \Spc, \quad N \mapsto \Map_{\Mod_B(\C)}(M, B \otimes_A N)
  \]
  is corepresentable by an $A$-module $f_\sharp(M)$. We may also consider the following stronger condition:
  \begin{enumerate}[label=($\dagger$), ref=$\dagger$]
  \item \label{e--q--suave--pr-pointwise}
    the left adjoint $f_\sharp$ is defined at $M \otimes_B f^*(L) \iso M \otimes_A L$ for all $A$-modules $L$, and the associated projection map $\pr^\sharp_f : f_\sharp(M \otimes_A L) \to f_\sharp(M) \otimes_A L$ is an isomorphism for all $A$-modules $L$.
  \end{enumerate}
  Condition \cref{e--q--suave--pr} is equivalent to \cref{e--q--suave--pr-pointwise} holding for all $B$-modules $M$. Note that the collection of $B$-modules $M$ for which \cref{e--q--suave--pr-pointwise} holds is closed under colimits in $\Mod_B(\C)$, and that any $B$-module $M$ can be written as a colimit of modules of the form $B \otimes_A M_0$, for $M_0$ an $A$-module (by virtue of the bar resolution). Thus, to show that \cref{e--q--suave--pr} holds, it suffices to show that \cref{e--q--suave--pr-pointwise} holds in the case $M = B \otimes_A M_0$. In that case, it is straightforward to see that $f_\sharp(M)$ is given by $\dual{B} \otimes_A M_0$, where $\dual{B}$ denotes the dual of $B$ in $\Mod_A(\C)$, and hence that \cref{e--q--suave--pr-pointwise} does indeed hold.

  We have now shown that \cref{e--q--suave--pr} holds (under the assumption that \cref{e--q--suave--dual} holds). Since condition \cref{e--q--suave--dual} is stable under base change, this also implies that \cref{e--q--suave--pr} holds for $\lp{f}$ in place of $f$. In particular, $f$ is in fact $\QCoh$-$\sharp$-normed. We conclude finally that \cref{e--q--suave--suave} holds by applying \cref{e--q--kunneth,t--d--norm-assembly,t--p--assembly}.
\end{proof}

\begin{remark}
  \label{e--q--suave-concrete}
  In the situation of \cref{e--q--suave}, when the stated conditions hold, we have a canonical identification $\suavetwist_f \iso \dual{B}$. Here $\dual{B}$ denotes the dual of $B$ in $\Mod_A(\C)$, which promotes to a $B$-module by virtue of being the image of $A$ under the right adjoint $f^! : \Mod_A(\C) \to \Mod_B(\C)$ to the restriction functor $f_! \iso f_* : \Mod_B(\C) \to \Mod_A(\C)$ (see \cref{t--p--op}).
\end{remark}

Now we may apply the main result of \cref{t--h} to obtain the following class of examples of $\QCoh$-Poincar\'e objects.

\begin{example}
  \label{e--q--poincare-hopf}
  Let $\G$ be a group object in $\Aff_\C$. Writing $\G = \Spec(B)$, the group structure on $\G$ corresponds to the promotion of $B$ from a commutative algebra in $\C$ to a commutative Hopf algebra in $\C$. Assume that $B$ is dualizable as an object of $\C$, so that $\G$ is $\QCoh$-suave by \cref{e--q--suave}. We also have that $\QCoh$-terse by \cref{e--q--prim}. Applying \cref{t--h--main}, we find that $\G$ is $\QCoh$-Poincar\'e. More precisely, what we have, in light of \cref{e--q--suave-concrete} and \cref{t--h--constant}, are isomorphisms of $B$-modules
  \[
    \poincaretwist_\G \iso \dual{B} \iso B \otimes \o\poincaretwist_\G,
  \]
  where $\o\poincaretwist_\G$ is an invertible object of $\C$ (note that $\poincaretwist_\G \iso \suavetwist_\G$ since $\primtwist_\G \iso \unit_\G$ by \cref{e--q--prim}).
\end{example}

\begin{example}
  \label{e--q--dualizable-group}
  \cref{e--q--poincare-hopf} connects to work of Rognes \cite{rognes--groups}, as follows. Let $\unit \in \C$ denote the unit object, and let $\unit[-] : \Spc \to \C$ denote the unit map in $\CAlg(\PrL)$, whose underlying functor is uniquely characterized by preserving colimits and carrying $*$ to $\unit$. Let $\C^\diamond \subseteq \C$ denote the full subcategory spanned by the dualizable objects, and let $\Spc^{\diamond_\C} \subseteq \Spc$ denote the full subcategory spanned by those spaces $T$ such that $\unit[T] \in \C^\diamond$; we will say that a space is \emph{$\C$-dualizable} if it lies in $\Spc^{\diamond_\C}$.

  Now, $\unit[-]$ induces a functor on cocommutative coalgebra objects
  \[
    \Spc \iso \cCAlg(\Spc) \to \cCAlg(\C),
  \]
  the first equivalence owing to the fact that the symmetric monoidal structure on $\Spc$ is cartesian. Restricting this to $\C$-dualizable spaces and then composing with the equivalence
  \[
    \dual{(-)} : \cCAlg(\C^\diamond) \iso \CAlg(\C^\diamond)^\op
  \]
  determined by duality \cite[Corollary 3.2.5]{lurie--elliptic-i}, we obtain a finite product preserving functor
  \[
    \Spc^{\diamond_\C} \to \Aff_\C, \quad T \mapsto \Spec(\unit^T),
  \]
  where $\unit^T$ denotes the dual of $\unit[T]$.
  
  Finally, let $G$ be a $\C$-dualizable group space, i.e. a group object in $\Spc^{\diamond_\C}$. Then $\G := \Spec(\unit^G)$ is an example of a group object in $\Aff_\C$ of the type discussed in \cref{e--q--poincare-hopf}. From the conclusion there we get a canonical isomorphim $\unit[G] \iso \unit^G \otimes \o\poincaretwist_\G$, where $\o\poincaretwist_\G$ is an invertible object of $\C$; this may be compared with the results of \cite[\textsection 3]{rognes--groups} (which apply in the case that $\C$ is a Bousfield localization of the $\infty$-category of spectra).
\end{example}

As emphasized by \cref{e--q--dualizable-group}, \cref{e--q--poincare-hopf} is an algebraic analogue of the fact explained in \cref{e--l--poincare-group} that compact group spaces are $\Loc_\C$-Poincar\'e. In \cref{e--l--poincare-group}, we were also able to derive consequences for the Tate cohomology of representations of such group spaces, by passing to the associated classifying spaces. So that we may make analogous observations in the current setting, let us now recall that there is a canonical procedure for extending our three functor formalism on $\Aff_\C$ to the associated $\infty$-category of \emph{prestacks}, $\PStk_\C := \P(\Aff_\C)$; here by a prestack we mean an accessible functor from $\Aff_\C^\op \iso \CAlg(\C)$ to $\Spc$, so that $\PStk_\C$ is the free cocompletion of $\Aff_\C$. We will freely identify $\Aff_\C$ with its essential image in $\PStk_\C$ under the Yoneda embedding.

\begin{proposition}
  \label{e--q--D-extended}
  There is a minimum wide subcategory $\PStk_{\C,!}$ of $\PStk_\C$ for which the following properties are satisfied:
  \begin{enumerate}
  \item \label{e--q--D-extended--geometric}
    The pair $(\PStk_\C,\PStk_{\C,!})$ is a geometric setup.
  \item \label{e--q--D-extended--star-local}
    Let $f : T \to S$ be a map in $\PStk_\C$ such that, for every cartesian square $\sigma :$
    \[
      \begin{tikzcd}
        V \ar[r, "q"] \ar[d, "g", swap] &
        T \ar[d, "f"] \\
        U \ar[r, "p"] &
        S,
      \end{tikzcd}
    \]
    in $\PStk_\C$ where $U \in \Aff_\C$, the map $g$ lies in $\PStk_{\C,!}$. Then $f$ lies in $\PStk_{\C,!}$.
  \item \label{e--q--D-extended--D}
    There exists a unique three functor formalism $\QCoh : \Span(\PStk_\C,\PStk_{\C,!}) \to \PrLst$ that extends the three functor formalism $\QCoh : \Span(\Aff_\C) \to \PrLst$ of \cref{e--q--D} and such that the composition
    \[
      \PStk_\C^\op \lblto{\iota_0} \Span(\PStk_\C,\PStk_{\C,!}) \lblto{\QCoh} \PrLst
    \]
    preserves limits (i.e. is right Kan extended from $\QCoh : \Aff_\C^\op \to \PrLst$).
  \end{enumerate}
  For the final two properties, we use the notion of $!$-cover from \itemitemref{t--f}{D}{cover}, with respect to the three functor formalism of \cref{e--q--D-extended--D}.
  \begin{enumerate}[resume]
  \item \label{e--q--D-extended--shriek-local}
    Let $f : T \to S$ be a map in $\PStk_\C$ that $\QCoh$-$!$-locally on the source or on the target lies in $\PStk_{\C,!}$. Then $f$ lies in $\PStk_{\C,!}$.
  \item \label{e--q--D-extended--tame}
    Let $f : T \to S$ be a map in $\PStk_{\C,!}$ where $S \in \Aff_\C$. Then $T$ admits a $\QCoh$-$!$-cover by objects of $\Aff_\C$.
  \end{enumerate}
\end{proposition}

\begin{proof}
  See \cite[Theorem 3.4.11]{heyer-mann--six}.
\end{proof}

Passing to the context of the three functor formalism of \itemref{e--q--D-extended}{D}, we gain access to the following examples of the results of \cref{t--u}, with which we end our discussion.

\begin{example}
  \label{e--q--poincare-hopf-B}
  As in \cref{e--q--poincare-hopf}, let $B$ be a dualizable commutative Hopf algebra in $\C$, so that $\G := \Spec(B)$ is a $\QCoh$-Poincar\'e group object in $\Aff_\C$. Then we may form the classifying prestack $\clspc{\G} \in \PStk_\C$. We have a canonical effective epimorphism $i : * \to \clspc{\G}$ and cartesian diagram
  \[
    \begin{tikzcd}
      \G \ar[r] \ar[d] &
      * \ar[d, "i"] \\
      * \ar[r, "i"] &
      \clspc{\G}
    \end{tikzcd}
  \]
  in $\PStk_\C$. From this and \itemref{e--q--D-extended}{star-local} we see that $i$ is $\QCoh$-$!$-able, and then \itemref{e--q--D-extended}{shriek-local} implies that $\clspc{\G}$ is $\QCoh$-$!$-able. We may then apply \cref{t--b--descent,t--c--descent} to find moreover that $i$ is $\QCoh$-Poincar\'e and that $\clspc{\G}$ is $\QCoh$-smooth, in particular $\QCoh$-Tate. Thus, we have homology and cohomology functors
  \[
    (-)_{\h\G} := \Ch_*(\clspc{\G};-) : \QCoh(\clspc{\G}) \to \C, \quad
    (-)^{\h\G} := \Ch^*(\clspc{\G};-) : \QCoh(\clspc{\G}) \to \C,
  \]
  as well as a norm map
  \[
    \Nm_{\clspc{\G}} : (- \otimes \poincaretwist_{\clspc{\G}})_{\h\G} \to (-)^{\h\G}
  \]
  whose cofiber gives the Tate cohomology functor
  \[
    (-)^{\t\G} := \Ct^*(\clspc{\G};-) : \QCoh(\clspc{\G}) \to \C.
  \]
  \cref{t--u--nilpotent-tate-vanishing} gives that $i$ is a $\QCoh$-Tate cover of $\clspc{\G}$, so that we may apply \cref{t--u--main}, and in particular get canonical lax symmetric monoidal structures on the Tate cohomology functor $(-)^{\t\G}$ and the canonical map $(-)^{\h\G} \to (-)^{\t\G}$.

  Regarding the twisting object $\poincaretwist_{\clspc{\G}}$, we may apply \cref{t--c--poincare-equiv} to the commutative triangle
    \[
      \begin{tikzcd}[row sep=small]
        * \ar[rr, "i"] \ar[dr] &
        &
        \clspc{\G} \ar[dl] \\
        &
        *,
      \end{tikzcd}
    \]
    and we obtain an identification
    \[
      i^*(\poincaretwist_{\clspc{\G}}) \iso \poincaretwist_i^{-1} \iso (\o\poincaretwist_\G)^{-1}
    \]
    in $\C$, where $\o\poincaretwist_\G$ is as in \cref{e--q--poincare-hopf} (the isomorphism $\poincaretwist_i \iso \o\poincaretwist_\G$ follows from the cartesian square above and \cref{t--b--poincare-equiv}).
\end{example}

\begin{remark}
  \label{e--q--poincare-hopf-B-addendum}
  We make a few additional remarks on the situation in \cref{e--q--poincare-hopf-B}:
  \begin{enumerate}
  \item \label{e--q--poincare-hopf-B-addendum--mod}
    It follows from the map $i$ being $\QCoh$-Poincar\'e that the functor $i^* : \QCoh(\clspc{\G}) \to \C$ admits $\C$-linear right and left adjoints. Since $i^*$ is also conservative, we may apply the comonadicity and monadicity theorems and thereby identify $\QCoh(\clspc{\G})$ with the $\infty$-category of comodules over the Hopf algebra $B$ and with the $\infty$-category of modules over the dual Hopf algebra $\dual{B}$.
  \item \label{e--q--poincare-hopf-B-addendum--assembly}
    One can check that the canonical functor $\QCoh(\clspc{\G}) \otimes_\C \QCoh(\clspc{\G}) \to \QCoh(\clspc{\G} \times \clspc{\G})$ is an equivalence (for example by using \cref{e--q--poincare-hopf-B-addendum--mod}, which applies also for $\clspc{\G} \times \clspc{\G} \iso \clspc{(\G \times \G)}$). Thus, by \cref{t--d--norm-assembly}, the norm map $\Nm_{\clspc{\G}}$ is a $\C$-linear assembly map.
  \item We may consider the case $\G = \Spec(\unit^G)$ for $G$ a $\C$-dualizable group space, as in \cref{e--q--dualizable-group}. Then there is a canonical identification
    \[
      \QCoh(\clspc{\G}) \iso \Loc_\C(\clspc{G}),
    \]
    the right hand side being the $\infty$-category of $\C$-valued local systems on the classifying space of $G$, as studied in previous subsection (the comparison functor is induced by a canonical map of prestacks $\u{\clspc{G}} \to \clspc{\G}$, where the source denotes the constant prestack on the space $\clspc{G}$, and that this functor is an equivalence may be checked by a monadicity or comonadicity argument as in \cref{e--q--poincare-hopf-B-addendum--mod}). Under this identification, the norm map of \cref{e--q--poincare-hopf-B} agrees with the norm map for local systems considered in \cref{e--l}; this follows for instance from \cref{e--q--poincare-hopf-B-addendum--assembly} and \cref{e--l--prim-absolute}. In this way, \cref{e--q--poincare-hopf-B} generalizes \cref{e--l--poincare-group} from the case of compact group spaces to the case of $\C$-dualizable group spaces (cf. Rognes's work \cite[\textsection 5]{rognes--groups}, which applies when $\C$ is a Bousfield localization of $\Spt$).
  \end{enumerate}
\end{remark}

\begin{example}
  \label{e--q--poincare-hopf-map}
  We may generalize \cref{e--q--poincare-hopf-B} by considering a map $B \to B'$ of dualizable commutative Hopf algebras in $\C$. Let $\G' \to \G$ be the corresponding map of group objects in $\Aff_\C$, and let $f : \clspc{\G'} \to \clspc{\G}$ denote the induced map on classifying spaces. Let $i : * \to \clspc{\G}$ and $j : * \to \clspc{\G'}$ denote the canonical basepoints. Then we have a commutative triangle
  \[
    \begin{tikzcd}[row sep=small]
      * \ar[rr, "j"] \ar[dr, "i", swap] &
      &
      \clspc{\G'} \ar[dl, "f"] \\
      &
      \clspc{\G}.
    \end{tikzcd}
  \]
  From \cref{e--q--poincare-hopf-B} we know that $i$ and $j$ are $\QCoh$-Poincar\'e, in particular $\QCoh$-smooth. Since $j$ is also an effective epimorphism, hence a $\QCoh$-$*$-cover, we deduce using \cref{t--c--descent} that $f$ is $\QCoh$-smooth, in particular $\QCoh$-Tate. Thus, we have a norm map and Tate cohomology functor associated to $f$, i.e. a cofiber sequence
  \[
    f_\sharp(- \otimes \poincaretwist_f) \lblto{\Nm_f} f_* \lblto{\can} f_\t
  \]
  of functors from $\QCoh(\clspc{\G'}) \to \QCoh(\clspc{\G})$. Moreover, $j$ is a $\QCoh$-Tate cover of $f$ by \cref{t--u--nilpotent-tate-vanishing}, so that \cref{t--u--main} gives lax symmetric monoidal structures to $f_\t$ and $\can : f_* \to f_\t$.

  To calculate the twisting object $\poincaretwist_f$, we consider also the commutative triangle
  \[
    \begin{tikzcd}[row sep=small]
      \clspc{\G'} \ar[rr, "f"] \ar[dr] &
      &
      \clspc{\G} \ar[dl] \\
      &
      *.
    \end{tikzcd}
  \]
  Since all three maps are $\QCoh$-smooth, we have an isomorphism $\suavetwist_{\clspc{\G'}} \iso \suavetwist_f \otimes f^*(\suavetwist_{\clspc{\G}})$ by \cref{t--c--suave-equiv}, and thus we may apply \cref{t--c--poincare-comp} to obtain a comparison map
  \[
    \poincaretwist_f \otimes f^*(\poincaretwist_{\clspc{\G}}) \to \poincaretwist_{\clspc{\G}'}.
  \]
  We claim that this is an isomorphism, so that $\poincaretwist_f \iso \poincaretwist_{\clspc{\G'}} \otimes f^*(\poincaretwist_{\clspc{\G}})^{-1}$ (recall from \cref{e--q--poincare-hopf-B} that $\poincaretwist_{\clspc{\G}}$ and $\poincaretwist_{\clspc{\G'}}$ are invertible). It suffices to check that this holds after applying $j^*$. For this, we return to the first triangle above. Since $j$ is $\QCoh$-Poincar\'e, \cref{t--c--poincare-equiv} gives an isomorphism $\poincaretwist_i \iso \poincaretwist_j \otimes j^*(\poincaretwist_f)$. The desired claim follows from this and the identifications $i^*(\poincaretwist_{\clspc{\G}}) \iso \poincaretwist_i^{-1}$ and $j^*(\poincaretwist_{\clspc{\G'}}) \iso \poincaretwist_j^{-1}$ of \cref{e--q--poincare-hopf-B}.
\end{example}


\bigskip
{\footnotesize
  \printbibliography}

\end{document}